\documentclass[a4paper,10pt,reqno]{amsart}
\usepackage{amsmath,amsfonts,amsthm,bbm}
\usepackage{amsaddr}
\usepackage{mathtools} 
\usepackage[mathscr]{eucal}
\usepackage{rotating,float,enumerate}
\usepackage{multirow}
\usepackage[left=1in, right=1in, bottom=1.25in, top = 1.25in]{geometry}
\numberwithin{equation}{section}
\newtheorem{theorem}{Theorem}[section]
\newtheorem{lemma}[theorem]{Lemma}
\newtheorem{remark}[theorem]{Remark}
\newtheorem{proposition}{Proposition}
\newtheorem{corollary}[theorem]{Corollary}

\theoremstyle{definition}
\newtheorem{definition}[theorem]{Definition}
\newtheorem{example}[theorem]{Example}
\usepackage[colorlinks=true]{hyperref}
\hypersetup{urlcolor=blue, citecolor=red, linkcolor=blue}
\usepackage{color}
\DeclareMathOperator{\sign}{sign}
\begin{document}

\title[Perturbed Nonlinear Volterra Equations]
{Growth and fluctuation in perturbed nonlinear\\ Volterra equations}

\author{John A. D. Appleby}
\address{\vspace*{-8pt}School of Mathematical Sciences, Dublin City University} \email{john.appleby@dcu.ie}
\author{Denis D. Patterson}
\address{\vspace*{-8pt} Princeton Environmental Institute, Princeton University \\ Department of Mathematics, Brandeis University} \email{denispatterson@princeton.edu} 

\date{\today}

\begin{abstract}
We develop precise bounds on the growth rates and fluctuation sizes of unbounded solutions of deterministic and stochastic nonlinear Volterra equations perturbed by external forces. The equation is sublinear for large values of the state, in the sense that the state--dependence is negligible relative to linear functions. If an appropriate functional of the forcing term has a limit $L$ at infinity, the solution of the differential equation behaves asymptotically like the underlying unforced equation when $L=0$, like the forcing term when $L=+\infty$, and inherits properties of both the forcing term and underlying differential equation for values of $L\in (0,\infty)$. Our approach carries over in a natural way to stochastic equations with additive noise and we treat the illustrative cases of Brownian and L\'evy noise.
\end{abstract}

\maketitle

\tableofcontents
\newpage
\setlength\parskip{0.5\baselineskip}
\section{Introduction}
\subsection{Problem Overview}
We analyse the long--run dynamics of solutions to the scalar Volterra integro-differential equation
\begin{equation} \label{eq.xpert}
x'(t)= \int_{[0,t]} \mu(ds) f(x(t-s))+h(t), \quad t>0, \quad x(0)=\psi \in \mathbb{R}.
\end{equation}
In particular, we concentrate on the behaviour of unbounded but non-explosive solutions, i.e. $x \in C(\mathbb{R}^+;\mathbb{R})$ but $\limsup_{t\to\infty}|x(t)|=\infty$. The nonlinearity $f$ is assumed to be \emph{sublinear} in the state variable in the sense that $\lim_{x\to\infty}f(x)/x=0$, guaranteeing global existence of solutions. We draw a distinction between when solutions of \eqref{eq.xpert} grow, $\lim_{t\to\infty}x(t)=\infty$, and when solutions can be said to fluctuate asymptotically, $\liminf_{t\to\infty}x(t) = - \infty$ and $\limsup_{t\to\infty}x(t)=+\infty$. When solutions grow it is natural to ask at what rate they grow and when they fluctuate to ask if the size of these fluctuations can be captured in an appropriate sense; these are the primary goals of this paper.

Nonlinear equations with after effect, such as \eqref{eq.xpert}, appear naturally in a myriad of diverse applications from models of nuclear reactors to heat flow or even financial management~\cite{kolmanovskii2012applied,kolmanovskii2013introduction}. In the present work, we are especially motivated by economic applications, such as to vintage capital models; in this context the sublinear response of the system represents saturated growth or diminishing returns to scale.  This analogy further motivates our study of \eqref{eq.xpert} in the presence of random forcing as \eqref{eq.xpert} encompasses a broad class of continuous time nonlinear time-series models when $h$ is replaced by an appropriate stochastic process~\cite{brockwell2009existence,marquardt2007multivariate}. The measure $\mu$ can be thought of as weighting the contributions of capital of different ages to the current growth in the economy, as well as ``time-to-build lags'' and other delay inducing effects~\cite{benhabib1991vintage}. Hence we assume $\mu$ to be a finite measure so that the influence of older capital becomes negligible as time advances (see \cite{sublinear2015} for further discussion and motivation with regard to applications). Precisely, $\mu$ is a finite Borel--measure on $\mathbb{R^+}:=[0,\infty)$, i.e.
\begin{equation}\label{finite_measure}
	\mu(E) \geq 0 \mbox{ for all } E \in \mathcal{B}(\mathbb{R^+}),\quad \mu(\mathbb{R}^+) \in (0,\infty),
\end{equation}
where $\mathcal{B}(\mathbb{R}^+)$ denotes the $\sigma$-algebra generated by the open sets of $\mathbb{R}^+$. We also define $m(t) = \mu([0,t])$ so that $\lim_{t\to\infty}m(t)=\mu(\mathbb{R}^+)$ and $H(t):= \int_{[0,t]} h(s)\,ds$ for $t \geq 0$. 

In the framework outlined above, the following is a convenient sufficient condition to guarantee a positive, growing solution to \eqref{eq.xpert} (see \cite[Theorem 1]{sublinear2015}):
\begin{equation} \label{eq.usualdata_growth}
f\in C(\mathbb{R}^+;(0,\infty)), \quad H\in C(\mathbb{R}^+;\mathbb{R}^+).
\end{equation}
After developing results regarding the qualitative behaviour of solutions to \eqref{eq.xpert} we extend our deterministic analysis to consider the asymptotic behaviour of the related stochastic Volterra equation
\begin{equation} \label{eq.introsfde}
dX(t) = \int_{[0,t]} \mu(ds) f(X(t-s))\,dt +dZ(t), \quad t>0,
\end{equation}
where $Z$ is a semimartingale. We establish an appropriate existence and uniqueness theorem for equation \eqref{eq.introsfde} and specialise to the cases of Brownian and L\'evy noise in order to prove precise asymptotic results.

Equations \eqref{eq.xpert} and \eqref{eq.introsfde} can both be viewed as perturbations of the underlying  Volterra integro-differential equation
\begin{equation} \label{eq.unpert}
y'(t)=\int_{[0,t]} \mu(ds) f(y(t-s)), \quad t>0, \quad y(0)=\psi.
\end{equation}
When $f$ is positive and sublinear at infinity, the solution $y(t)$ of \eqref{eq.unpert} obeys $y(t)\to \infty$ as $t\to\infty$ and grows asymptotically like the solution of the ODE 
$z'(t) = \mu(\mathbb{R}^+)f(z(t))$, i.e. the corresponding ODE with the mass of the measure concentrated at zero~ \cite{sublinear2015}.  It is natural to ask how large the forcing terms $h$ in \eqref{eq.xpert} and $Z$ in \eqref{eq.introsfde} can become while the solutions $x$ of \eqref{eq.xpert} and $X$ of \eqref{eq.introsfde} continue to grow in the same manner as solutions to $z'(t) = \mu(\mathbb{R}^+)f(z(t))$. Furthermore, can we identify a new asymptotic regime or growth rate if the forcing terms exceed this critical rate? The main goal of this paper is to identify such critical rates of growth on $h$ and $Z$, and to determine precise estimates on the growth rate of solutions, or the rate of growth of the partial maxima when solutions fluctuate. 

The analysis in our paper~\cite{sublinear2015} deals with growth estimates for the unperturbed version of \eqref{eq.xpert} (i.e. $h\equiv 0$) and can be thought of as demonstrating the asymptotic sharpness of Bihari's inequality (and related retarded integral inequalities) for a large class of functional differential equations~\cite{bihari1956generalization,lipovan2006}. It is a well-established, but still active, area of work to determine the interplay between the properties of perturbations and fundamental solutions or unperturbed equations. In the context of linear and quasilinear integral equations, the following formulation has proven fruitful: If the perturbation belongs to a function class $\mathcal{C}$, and the equation (operator) is stable in an appropriate sense, then the solution also belongs to the class $\mathcal{C}$. This phenomenon is referred to as admissibility, and an excellent account by one the most important contributors and initiators of this theory is given in Corduneanu~\cite{corduneanu1991integral}. In the linear case, properties of the fundamental solution and variations of constants formulae are especially helpful in obtaining these admissibility results, but naturally nonlinear equations require entirely different methods. 

A central contribution of this paper is to prove results in the spirit of classical linear admissibility theory, exploiting methods more suited to dealing with nonlinear equations in which the state-dependence is \emph{sublinear}. We give precise estimates on the asymptotic behaviour of solutions whether perturbations are large or small. In doing this, we bear in mind that the great bulk of research focuses on bounded perturbations with different properties or exploits the theory of $L^p$ weighted spaces~\cite{GLS}. In the latter case, $p = + \infty$ (which is our focus here) tends to attract relatively little attention. Our work treats, in a unified manner, large or rapidly fluctuating perturbations which may be either deterministic or stochastic in nature. For many classes of linear differential equations with additive forcing (with or without memory) the asymptotic behaviour of solutions tends to be in one of two regimes. When the perturbation is sufficiently small, the solution tracks that of the unperturbed equation asymptotically; this can even extend to precise quantitative measures, such as Lyapunov exponents, being preserved  (see \cite[Ch. 10]{hartman2002ordinary} and the references therein, \cite{pituk1999hartman} and Proposition \ref{prop.linear} below). The second typical regime is when the perturbation becomes so large that the perturbed solution no longer behaves at all like the fundamental solution and instead the forcing term dominates the dynamics (cf. \cite{appleby2017growth}). The challenge is to characterise precisely the appropriate quantities that are preserved in the case of small perturbations and the critical perturbation size at which the regime shift occurs. There may be interesting dynamics at the critical transition between the small and large perturbative regimes which can sometimes also be characterised; we show that this can be achieved for \eqref{eq.xpert}, as well as it's stochastic counterpart \eqref{eq.introsfde}.   

There is an extensive literature regarding existence, uniqueness and regularity results for stochastic Volterra equations subject to additive Brownian forcing~\cite{berger1980volterra,bharucha1972random}. More recently, there has been considerable progress in developing stability theory for linear stochastic Volterra equations~\cite{appleby2005mean,appleby2007almost,mao2006mean,reynolds2008decay,wu2012lasalle} and stochastic delay differential equations~\cite{appleby2003non,mao2007stochastic,mohammed2003stable,mohammed2004stable}, particularly with an instantaneous diffusion term (i.e. delays are confined to the drift term). Moreover, several authors have extended the stability theory to certain classes of nonlinear equations, typically using Lyapunov methods~\cite{appleby2004asymptotic,kolmanovskii2013introduction,shaikhet2013lyapunov,wu2010robustness}. There has also been some limited exploration of global estimates on solutions in the absence of asymptotic stability for the small noise regime~\cite{nualart2000large,zhang2008euler}. In contrast, we explicitly explore both the case where the noise term is small relative to the nonlinear drift term in \eqref{eq.introsfde} and the situation in which the state-independent noise term dominates the dynamics. This is an important first step towards understanding the dynamics of \eqref{eq.introsfde} with state dependent noise, e.g. $dZ(t) = G(t,\,X(t))\,dB(t)$ where $B$ denotes Brownian motion. There are relatively few results in the extant literature concerning the qualitative or asymptotic behaviour of stochastic Volterra equations with non-Brownian noise, although some authors have considered linear equations with quite general noisy driving processes~\cite{mohammed1996lyapunov}. Hence the combination of nonlinearity, general stochastic noise and memory in \eqref{eq.introsfde} is especially novel; in particular, we will prove asymptotic results with $\alpha$-stable noise, in which case the solution to \eqref{eq.introsfde} is a non-Markovian jump process with nonlinear state dependence. 
\subsection{Outline \& Motivation}
Much of our analysis flows from the simple matter of integrating \eqref{eq.xpert} to obtain the forced Volterra integral equation  
\begin{equation}\label{eq.volterra}
x(t) =x(0) + \int_0^t m(t-s)f(x(s))\,ds + H(t), \quad t \geq 0.
\end{equation}
Since stochastic ``differential'' equations must be rigorously formulated in integral form, it is perhaps even more natural to treat \eqref{eq.introsfde} similarly; this leads us to consider 
\begin{equation} \label{eq.introsintde}
X(t)=X(0) + \int_0^t m(t-s)f(X(s))\,ds + Z(t), \quad t\geq 0. 
\end{equation}
Equation \eqref{eq.volterra} shows that the solution to \eqref{eq.xpert} is a functional of the aggregate behaviour of the forcing term $h$ purely through $H$ and hence it is natural to formulate asymptotic results in terms of $H$. When studying the asymptotic behaviour of many forced differential systems, hypotheses on the aggregate or average behaviour of the forcing terms are preferable to more restrictive pointwise conditions. When studying stochastic equations pointwise estimates become unrealistically restrictive and it is more natural, perhaps even necessary, to consider average behaviour. In this spirit, Proposition \ref{prop.linear} below illustrates how hypotheses on the averaged behaviour of perturbations can be used to classify the behaviour of the elementary perturbed linear ordinary differential equation (ODE) $x'(t) = ax(t) + h(t)$ shown in \eqref{eq.linearODE}; the proof is a simple matter of applying the variation of constants formula and integration by parts multiple times (see Section \ref{sec_deterministic_proofs} for details). The asymptotic behaviour of the solution to this linear ODE is characterized via the functional
\begin{equation}\label{eq.L_linear}
	 L = \lim_{t\to\infty}\frac{H(t)}{a\int_0^t H(s)\,ds}, \quad H(t) = \int_0^t h(s)\,ds, 
\end{equation}
quantifies size of the perturbation term. When $L = 0$, the perturbation is small (in some appropriate sense), and the solution to \eqref{eq.linearODE} is asymptotic to the unperturbed solution in case (i.). On the other hand, when perturbations grow more rapidly than some critical rate, the solution can track the perturbation asymptotically; this is case (iii.) of Proposition \ref{prop.linear} when $L = \infty$ and the solution is asymptotic to the perturbation term as $t\to\infty$. 
\begin{proposition}\label{prop.linear}
	Consider the nonautonomous linear ordinary differential equation given by
	\begin{equation}\label{eq.linearODE}
	x'(t) = a\,x(t) + h(t), \quad x(0) >0, \quad a >0.
	\end{equation}
	Suppose $H(t) :=  \int_0^t h(s)\,ds \geq 0$ for all $t \geq 0$ and let $L$ be defined as in \eqref{eq.L_linear}.
	\begin{enumerate}[(i.)]
		\setlength\itemsep{5pt}
		\item If $L = 0$, then $x(t)/e^{at} \to \xi^* \in (0,\infty)$ as $t\to\infty$.
		\item If $L \in (0,\infty)$, then
		\begin{itemize}
			\setlength\itemsep{5pt}
			\item $L<1$ implies $x(t)/e^{at} \to \xi^* \in (0,\infty)$ as $t\to\infty$,
			\item $L = 1$ implies $x(t)/H(t) \to \infty$ and $\log(x(t))/t \to a$ as $t\to\infty$,
			\item $L > 1$ implies $x(t)/H(t) \to L/(L-1)$ as $t\to\infty$. 
		\end{itemize}
		\item If $L = \infty$, then $x(t)/H(t) \to 1$ as $t\to\infty$.
	\end{enumerate}
\end{proposition}
When the perturbation is precisely at some ``critical size'', there is often an intermediate regime where the solution to a perturbed differential system inherits properties of both the unperturbed solution and the perturbation term (cf. \cite[Corollary 1]{appleby2017growth}). This regime requires a certain asymptotic balance in the sense that the perturbation term and the unperturbed solution should be of roughly the same order of magnitude. We observe this situation even for the simple linear ordinary differential equation \eqref{eq.linearODE} in Proposition \ref{prop.linear} case (ii.); we prove results of a similar character for \eqref{eq.xpert}, even for stochastic perturbations.
 
We now specify our hypotheses on the nonlinearity and outline typical results for the nonlinear Volterra equation \eqref{eq.xpert}. For solutions of the unperturbed Volterra equation \eqref{eq.unpert} to behave similarly to those of the corresponding nonlinear ordinary differential equation with the measure concentrated at zero, i.e.
\begin{equation}\label{eq.introode}
	z'(t) = \mu(\mathbb{R}^+)f(z(t)), \quad t > 0, \quad z(0)>0,
\end{equation}
it is important that $f$ be \emph{sublinear}. In previous work we showed that if $f$ is asymptotic to a $C^1$ function $\phi$ which is increasing and obeys $\phi'(x)\to 0$ as $x\to\infty$ (a hypothesis implying sublinearity of $f$), then
the solution to \eqref{eq.unpert} obeys
\begin{equation} \label{eq.unpertgrowth}
	\lim_{t\to\infty} \frac{F(y(t))}{\mu(\mathbb{R}^+)t}=1,
\end{equation}
where $F$ is the function defined by 
\begin{equation} \label{def.F}
	F(x)=\int_1^x \frac{1}{f(u)}\,du, \quad x>0
\end{equation}
(see \cite{sublinear2015} for further details). \emph{Sublinearity} is crucial to this result since, for example, the linear Volterra  equation of the form \eqref{eq.unpert} does not share the exponential rate of growth of the linear ODE with all of the mass of $\mu$ concentrated at zero (cf. \cite[Theorem 7.2.3]{GLS}). However, the distribution of the measure $\mu$, as opposed to simply it's total mass, can impact rates of asymptotic growth when $f$ is sublinear in the case where $m(t) \to \infty$ as $t\to\infty$~\cite{appleby2017memory}. We retain the aforementioned hypothesis on $f$ and occasionally strengthen it so that $\phi'(x)$ decays monotonically to $0$ as $x\to\infty$; the implications and technical motivations for such hypotheses are discussed in Section \ref{det.results}.

Before stating our main results precisely we give a heuristic argument as to their likely validity. In this discussion we consider the simple (deterministic) case in which both the solution and the perturbation are positive. If the unperturbed equation \eqref{eq.unpert} is integrated as above, $H\equiv 0$. In this case, 
the solution $y$ of the integral equation is roughly of order $F^{-1}(\mu(\mathbb{R}^+)t)$, like the solution $z(t) = F^{-1}\left(\mu(\mathbb{R}^+)t+F(\psi)\right)$ to the nonlinear ODE \eqref{eq.introode}. This leads to the naive idea that if $H$ is of smaller order than $y$ (i.e., smaller than $F^{-1}(\mu(\mathbb{R}^+)t)$), then $H$ on the right--hand side of  \eqref{eq.volterra} could be absorbed into $x$ on the left--hand side, without changing the leading order asymptotic behaviour of $x$. However, if $H$ dominates $y$, or is of comparable order, such an outcome is improbable and the asymptotic behaviour of $x$ is unlikely to be determined by $y$. Since the asymptotic behaviour of \eqref{eq.unpert} is described well by the asymptotic relation $F(y(t))/\mu(\mathbb{R}^+)t\to 1$ as $t\to\infty$, and $F^{-1}$ is increasing, it is natural to characterise the forcing term as ``small'' or ``large'' according as to whether $F(H(t))/\mu(\mathbb{R}^+)t$ tends to a small or large limit as $t\to\infty$ (if such a limit exists). Hence we define the dimensionless parameter $L\in [0,\infty]$ by 
\begin{equation} \label{eq.introL}
\lim_{t\to\infty}\frac{F(H(t))}{\mu(\mathbb{R}^+)t} = L.
\end{equation}
In some sense $L=1$ is critical; for $L<1$, $H$ is dominated by the solution of \eqref{eq.unpert}. But for $L>1$, $H$ dominates the solution of \eqref{eq.unpert}. The cases $L=0$ and $L=+\infty$ are especially decisive; in these cases it is very clear whether the solution of the unperturbed equation or the perturbation dominates. A condition which implies \eqref{eq.introL}, and turns out to be very useful in classifying asymptotic behaviour, is 
\begin{equation}\label{intromed_pert}
\lim_{t\to\infty}\frac{H(t)}{\mu(\mathbb{R}^+)\int_0^t f(H(s))ds} = L. 
\end{equation}
If $L=0$ in \eqref{intromed_pert}, then 
\[
\lim_{t\to\infty} \frac{F(x(t))}{\mu(\mathbb{R}^+)t}=1,  \quad \lim_{t\to\infty} \frac{x(t)}{H(t)}=+\infty,
\]
so small perturbations give rise to asymptotic behaviour as in \eqref{eq.unpert}, 
and the solution dominates the perturbation. If $L=+\infty$, then 
\[
\lim_{t\to\infty} \frac{x(t)}{H(t)}=1, \quad \lim_{t\to\infty} \frac{F(x(t))}{\mu(\mathbb{R^+})t}=+\infty, 
\]
so large perturbations cause the solution to grow at exactly the same rate as $H$, and the solution grows much faster than the original unperturbed Volterra equation. When the perturbation is of a scale comparable to the solution of \eqref{eq.unpert}, in the sense that $L\in (0,\infty)$, 
\begin{equation}\label{eq.intro.limits}
1 \leq \liminf_{t\to\infty} \frac{F(x(t))}{\mu(\mathbb{R}^+)t}\leq \limsup_{t\to\infty} \frac{F(x(t))}{\mu(\mathbb{R}^+)t}\leq 1+L, \quad \liminf_{t\to\infty} \frac{x(t)}{H(t)}\geq 1+\frac{1}{L}.
\end{equation}
Examples show that the limits in the first part of \eqref{eq.intro.limits} are not, in general, equal to $1$ or $1+L$. Further investigation for finite and positive $L$ leads to better estimates, especially when $L>1$. The critical character of the case when $L=1$ is 
demonstrated by the following result: if $L\in (1,\infty)$ then 
\begin{align}\label{x_over_H_intro}
1\leq \liminf_{t\to\infty} \frac{x(t)}{H(t)}\leq \limsup_{t\to\infty} \frac{x(t)}{H(t)}  \leq \frac{L}{L-1}.
\end{align}
This provides sharper estimates for large $L$ than the asymptotic bounds given for $L\in (0,\infty)$ above and identifies that $x$ is of order $H$. We also show by means of examples that when $L\in(0,1]$, the limit 
\[
\lim_{t\to\infty} \frac{x(t)}{H(t)}=+\infty
\]
can result, so that $x$ can only be expected to be exactly of the order of $H$ for $L>1$ (see Example \ref{example_L=1}). However, if $L\in (0,1]$, it is not necessarily the case that $x(t)/H(t)\to\infty$ as $t\to\infty$ (see Example \ref{golden}). As $L\to\infty$, equation \eqref{x_over_H_intro} correctly anticipates that $x(t)/H(t)\to 1$ as $t\to\infty$, which is what pertains when $L=+\infty$. To generalise the analysis above to stochastic equations, and for notational convenience, we define the following functional for later use:
\begin{equation}\label{L_functional}
L_f(\gamma) = \lim_{t\to\infty}\frac{\gamma(t)}{\mu(\mathbb{R}^+)\int_0^t f(\gamma(s))\,ds}, 
\end{equation} 
for all functions $f$ and $\gamma \in C(\mathbb{R}^+;(0,\infty))$ such that the above limit is well defined. For deterministic equations we will typically choose $\gamma=H$ in \eqref{L_functional} but other choices will prove advantageous when considering stochastic perturbations.

The rest of the paper is organized as follows: in Section \ref{det.results} we provide the mathematical framework for studying solutions to \eqref{eq.xpert}, state our main theorems for both growing and fluctuating solutions, and provide examples to illustrate the strengths and limitations of our results. Section \ref{stochastic} contains results for the stochastic Volterra equation \eqref{eq.introsintde}; we first prove a general existence theorem for solutions to \eqref{eq.introsintde} and then proceed to extend our deterministic results to cover Brownian and $\alpha$-stable L\'evy noise. We also present some applications and numerical simulations to illustrate our stochastic results. Section \ref{sec_deterministic_proofs} contains the proofs for Section \ref{det.results} on deterministic Volterra equations and Section \ref{sec_stoch_proofs} contains the proofs for Section \ref{stochastic} on stochastic Volterra equations. The interested reader can find detailed justification of all examples and an outline of the numerical scheme used to produce Figure \ref{fig_sharpness} in Appendix \ref{sec_examples}.
\section{Deterministic Volterra Equations}\label{det.results}
\subsection{Mathematical Preliminaries}
We briefly recall the relevant existence and uniqueness theory for the deterministic Volterra equation \eqref{eq.xpert} in order to keep the presentation self contained. 
\begin{definition}\label{defn.solution}
	A function solves the initial value problem \eqref{eq.xpert} if it obeys \eqref{eq.xpert} almost everywhere on an interval containing zero and is absolutely continuous on that interval. A solution which obeys \eqref{eq.xpert} for almost every $t \geq 0$ is called a \emph{global solution}. A solution which obeys \eqref{eq.xpert} for almost every $t \in [0,T]$ for some $T>0$ is called a \emph{local solution}. 
\end{definition}
The following result guarantees the existence of a local solution to \eqref{eq.xpert} in the sense of Definition \ref{defn.solution}~\cite[Corollary 12.3.2]{GLS}.
\begin{theorem}[Local Existence Theorem]\label{existence_local}
	Suppose $\mu$ is a Borel measure on $\mathbb{R}^+$, $h \in L^1_{loc}(\mathbb{R}^+;\mathbb{R})$, and $f \in C(\mathbb{R};\mathbb{R})$. Then, for each $\psi \in \mathbb{R}$, there exists a locally absolutely continuous solution $x$ to \eqref{eq.xpert} on an interval $[0,T]$ for some $T>0$. 	Moreover, every solution to \eqref{eq.xpert}, defined on some interval $[0,T]$, can be continued to a noncontinuable solution on $[0,T_{max})$ for some $T_{max}>T$. If $T_{max}<\infty$, then $\limsup_{t\to T_{max}^-}|x(t)|=\infty$.
\end{theorem}
We assume throughout that $\lim_{|x|\to\infty}f(x)/x = 0$ (or stronger hypotheses implying this) and hence $f$ will always obey a global linear bound of the form 
\[
|f(x)| \leq A + B|x|, \mbox{ for all }x \in \mathbb{R}, \mbox{ for some positive constants } A \mbox{ and }B.
\]
Therefore solutions to \eqref{eq.xpert} are always globally defined in the present framework. Moreover, Gripenberg at al. \cite[Theorem 13.5.1]{GLS} guarantees that solutions to equation \eqref{eq.xpert} are unique under the hypotheses of Theorem \ref{existence_local} if $f$ is locally Lipschitz continuous.

We next define a useful equivalence relation on the space of positive continuous functions; in essence, two functions are equivalent if they have the same leading order asymptotic behaviour.
\begin{definition}
$f,\,\phi \in C((0,\infty);(0,\infty))$ are asymptotically equivalent if $\lim_{x\to\infty}f(x)/\phi(x)=1$; written $f(x) \sim \phi(x)$ as $x\to\infty$ for short.
\end{definition}
$f(x) \sim \phi(x)$ implies $1/f(x) \sim 1/\phi(x)$ as $x\to\infty$ and $\lim_{x\to\infty}f(x)/x=0$ implies that $F$, defined by \eqref{def.F}, obeys $\lim_{x\to\infty}F(x)=\infty$. Hence the following convenient lemma can be proven immediately by asymptotic integration. 
\begin{lemma}\label{capF_sim_capPhi}
If $f,\,\phi \in C((0,\infty);(0,\infty))$ are asymptotically equivalent and obey \[
\lim_{x\to\infty}\frac{f(x)}{x}=\lim_{x\to\infty}\frac{\phi(x)}{x}=0,\quad\lim_{x\to\infty}f(x)=\lim_{x\to\infty}\phi(x) = \infty,
\]
then
$
F(x) \sim \Phi(x) \text{ as } x\to\infty,
$
where $F$ is defined by \eqref{def.F} and $\Phi(x)$ is defined by
\begin{align}\label{def.Phi}
\Phi(x)= \int_1^x \frac{1}{\phi(u)}du, \quad x >0.
\end{align}
\end{lemma} 
We impose the following sublinearity hypothesis on the nonlinear function $f$:
\begin{align}\label{fasym}
f \sim \phi \in C^1 \mbox{ such that } \lim_{|x|\to\infty}\phi(x)=\infty,\,\, \phi'(x) > 0 \mbox{ for all } x\in\mathbb{R}\mbox{ and }\phi'(x) \to 0 \mbox{ as } |x| \to \infty.
\end{align}
In many cases the following slightly stronger hypothesis is necessary
\begin{align}\label{fasym2}
f \sim \phi \in C^1 \mbox{ such that } \lim_{|x|\to\infty}\phi(x)=\infty,\,\, \phi'(x) > 0 \mbox{ for all } x\in\mathbb{R} \mbox{ and } \phi'(x) \downarrow 0 \mbox{ as } |x| \to \infty.
\end{align}
If $f$ is an increasing, sublinear function, then $\liminf_{x\to\infty}f'(x)=0$ but it is still possible that\\ $\limsup_{x\to\infty}f'(x)=\infty$ in the ``worst'' case. In previous work we provided an example of such a pathological $f$ but such nonlinearities are unlikely to arise naturally in applications so condition \eqref{fasym} is a relatively mild strengthening of sublinearity in this context~\cite{sublinear2015}. Assuming further that $\phi'$ tends to zero monotonically, as in \eqref{fasym2}, one can establish the following lemmata which often prove crucial in the asymptotic analysis of \eqref{eq.xpert} and \eqref{eq.introsfde}. 
\begin{lemma}\label{phi_props}
If \eqref{fasym2} holds, then $\phi$ obeys
\begin{align}\label{props}
\limsup_{x\to\infty}\frac{x\,\phi'(x)}{\phi(x)} \leq 1, \quad \limsup_{x\to\infty}\frac{\phi(\Lambda x)}{\phi(x)} \leq \Lambda, \quad \Lambda \in [1,\infty).
\end{align}
\end{lemma}
The conclusions of Lemma \ref{props} are remarkably close to some of the key properties enjoyed by the class of regularly varying functions with unit index, denoted  $RV_\infty(1)$. Namely, $\phi \in RV_\infty(1)$ implies $\lim_{x\to\infty}\phi(\Lambda x)/\phi(x) = \Lambda$ for all $\Lambda > 0$ and $\lim_{x\to\infty}x\,\phi'(x)/\phi(x)=1$~\cite{BGT}. The next lemma shows that the auxiliary function $\phi$ from \eqref{fasym2} preserves asymptotic equivalence (cf. \cite[Ch. 3, Problem 2]{dieudonne1971infinitesimal}) and hence $L_f(\gamma) = L_\phi(\gamma)$, if the limit exists.
\begin{lemma}\label{asym_equiv}
If \eqref{fasym2} holds, then the function $\phi$ preserves asymptotic equivalence, i.e. if $x,y \in C(\mathbb{R^+},(0,\infty))$ obey $\lim_{t\to\infty}x(t)=\lim_{t\to\infty}y(t) = \infty$, and $x(t) \sim y(t)$ as $t\to\infty$, then $\phi(x(t)) \sim \phi(y(t))$ as $t\to\infty$.
\end{lemma}
The connection between the ``natural'' size hypothesis on $H$, \eqref{eq.introL}, and the functional condition, \eqref{L_functional}, is supplied by the following result.
\begin{proposition}\label{PhiafterH}
Suppose $\phi\in C(\mathbb{R}^+;(0,\infty))$ is increasing and continuous with $\Phi$ defined by \eqref{def.Phi}. Let $\gamma \in C(\mathbb{R}^+;(0,\infty))$. If $L_{\phi}(\gamma)$ from \eqref{L_functional} is well defined, then
\[\lim_{t\to\infty}\frac{\Phi(\gamma(t))}{\mu(\mathbb{R}^+)t} = L_{\phi}(\gamma).
\]
\end{proposition}
Occasionally, we employ the standard Landau ``O'' and ``o'' notation. For $a$ and $b$ in $ C(\mathbb{R}^+;\mathbb{R})$, $b\text{ is }O(a)$ if $|b(t)|\leq K a(t)$ for some $K\in (0,\infty)$ and $t$ sufficiently large, and $b \text{ is }o(a)$ if $b(t)/a(t)\to 0$ as $t\to\infty$. 

\subsection{Growth Results}
With the requisite preliminaries in hand, we now turn to the computation of explicit growth estimates for solutions to \eqref{eq.xpert}. Suppose that \eqref{eq.usualdata_growth} holds so that $0< x(t) \to \infty$ as $t\to\infty$, subject to a positive initial condition. Our first result provides an easy to check sufficient condition on $H$ which guarantees solutions of \eqref{eq.xpert} retain the rate of growth of solutions to the ordinary differential equation \eqref{eq.introode}. This sufficient condition is of a different character to conditions involving the functional $L_f(\cdot)$ and expresses more explicitly the idea that the perturbation term, H, should be small relative to the solution of \eqref{eq.introode}.
\begin{theorem}\label{Thm.det.epsilon}
Suppose \eqref{finite_measure}, \eqref{eq.usualdata_growth}, and \eqref{fasym} hold and $\psi>0$. If
\begin{equation}\label{suff_eps}
\lim_{t\to\infty}\frac{H(t)}{F^{-1}\left(\mu(\mathbb{R}^+)(1+\epsilon)t \right)}=0 \mbox{ for each } \epsilon \in (0,1),
\end{equation}
then solutions of \eqref{eq.xpert} obey
\begin{align}\label{eq.thm.Finv}
\lim_{t\to\infty}\frac{F(x(t))}{\mu(\mathbb{R}^+)t} = 1, \quad \lim_{t\to\infty}\frac{x(t)}{H(t)}=\infty.
\end{align}
\end{theorem}
\noindent Now we formulate a sufficient condition for $\lim_{t\to\infty}F(x(t))/\mu(\mathbb{R}^+)t = 1$ to hold in terms of $L_f(\cdot)$. We also prove that when the solution of \eqref{eq.xpert} retains the growth rate of solutions of \eqref{eq.introode} it is of a strictly larger order of magnitude than the perturbation term, H.
\begin{theorem}\label{Thm.det.zero}
Suppose \eqref{finite_measure}, \eqref{eq.usualdata_growth}, and \eqref{fasym} hold and $\psi>0$. If $L_{f}(H)=0$, then solutions of \eqref{eq.xpert} obey
\begin{equation} \label{eq.thm21FxtMxHinf}
\lim_{t\to\infty}\frac{F(x(t))}{\mu(\mathbb{R}^+)t} = 1, \quad \lim_{t\to\infty} \frac{x(t)}{H(t)}=\infty.
\end{equation}
\end{theorem}
Notably, we do not assume that $H(t)\to\infty$ as $t\to\infty$ in Theorem~\ref{Thm.det.zero}; this is in the case where $L_f(H)=0$. However, if $L_f(H) \in (0,\infty]$, then $\lim_{t\to\infty} H(t)=\infty$. The rationale is as follows in the case $L_f(H)\in (0,\infty)$, with the case of $L_f(H)=\infty$ being similar. By hypothesis $H(t)>0$ for $t>0$ and as $f$ is a positive function, $t\mapsto \int_0^t f(H(s))ds$ is increasing. Therefore, $H$ either tends to $\infty$ or to a finite limit. In the former case, $H(t)\to\infty$ as $t\to\infty$ automatically. If, to the contrary, the limit is finite, then $H(t)$ tends to a finite positive limit as $t\to\infty$. But this forces $\int_{0}^t f(H(s))\,ds\to\infty$ as $t\to\infty$, a contradiction.   

When $L_f(H)$ is nonzero but finite we expect the solution of \eqref{eq.xpert} to inherit properties of both the ODE \eqref{eq.introode} and the perturbation term. Our next theorem investigates results of the type \eqref{eq.unpertgrowth} when $L_f(H)\in (0,\infty)$; we show that the growth of solutions to \eqref{eq.xpert} is at least as fast as that of solutions to the ODE \eqref{eq.introode} and we prove an upper bound on the growth rate. The resulting upper bound is linear in $L_f(H)$ and this is intuitively appealing as a ``larger'' $H$ should speed up growth. However, this upper estimate on the growth rate is not sharp in general. Without additional hypotheses this upper bound is hard to improve but can be shown to be suboptimal for specific classes of nonlinearity, for example when $f$ is regularly varying with less than unit index. We will demonstrate this possible improvement in further work.
\begin{theorem}\label{Thm.det.L.F}
Suppose \eqref{finite_measure}, \eqref{eq.usualdata_growth}, and \eqref{fasym} hold and $\psi>0$. If $L_f(H) \in (0,\infty)$, then solutions of \eqref{eq.xpert} obey
\[
1 \leq \liminf_{t\to\infty}\frac{F(x(t))}{\mu(\mathbb{R}^+)t} \leq \limsup_{t\to\infty}\frac{F(x(t))}{\mu(\mathbb{R}^+)t} \leq 1+L_f(H).
\]
If \eqref{fasym} is strengthened to \eqref{fasym2}, solutions of \eqref{eq.xpert} also obey
\[
\liminf_{t\to\infty}\frac{x(t)}{H(t)} \geq 1 + \frac{1}{L_f(H)}.
\]
\end{theorem}
The asymptotic lower bound on the quantity $x(t)/H(t)$ in the result above agrees with Theorem \ref{Thm.det.zero} as $L_f(H) \to 0^+$ since it correctly predicts that $\lim_{t\to\infty}x(t)/H(t)=\infty$ when $L_f(H)=0$. In some sense, the case where $L_f(H) \in (0,\infty)$ is special since the perturbation term is approximately the same order of magnitude as the solution to the unperturbed equation \eqref{eq.unpert}. More precisely, $L_f(H) \in (0,\infty)$ implies that $\lim_{t\to\infty}F(H(t))/\mu(\mathbb{R}^+)t = L_f(H)$. Furthermore, if $F^{-1}$ preserves asymptotic equivalence (see Lemma \ref{asym_equiv}) and $L_f(H) \in (0,\infty)$, then $H(t) \sim F^{-1}(\mu(\mathbb{R}^+)L_f(H)t)$ as $t\to\infty$ where $F^{-1}(\mu(\mathbb{R^+})t)$ is the leading order behaviour of the solution to \eqref{eq.unpert}. In many practical cases of interest, such as when $f \in \text{RV}_\infty(\beta)$ for $\beta\in(0,1)$, $F^{-1}$ inherits the asymptotic preserving property from $f$ (see \cite[Theorem 5]{sublinear2015}). 

When $L_f(H)>1$, we can additionally provide upper bounds on the quantity $x(t)/H(t)$. Moreover, these bounds are sharp in the limit as $L_f(H)\to\infty$.
\begin{theorem}\label{Thm.det.L.H}
	Suppose \eqref{finite_measure}, \eqref{eq.usualdata_growth}, \eqref{fasym2} hold and that $\psi>0$. Let $x$ denote a solution of \eqref{eq.xpert}.
	\begin{itemize}
		\item[(a.)] If $L_f(H) \in (1,\infty)$, then
		\[
		G_L := 1+\frac{1}{L_f(H)} \leq \liminf_{t\to\infty}\frac{x(t)}{H(t)} \leq \limsup_{t\to\infty}\frac{x(t)}{H(t)} \leq \frac{L_f(H)}{L_f(H)-1} =: G_U.
		\]
		\item[(b.)] If $L_f(H) =\infty$, then
		\begin{align}\label{eq.thm25FxtinfxH1}
		\lim_{t\to\infty}\frac{x(t)}{H(t)} = 1, \quad \lim_{t\to\infty}\frac{F(x(t))}{\mu(\mathbb{R}^+)t} = \infty.
		\end{align}
	\end{itemize}
\end{theorem}
Albeit under stronger hypotheses, Theorem \ref{Thm.det.L.H} provides more refined conclusions than Theorems \ref{Thm.det.gamma} and \ref{Thm.det.gamma.pm} (see Section \ref{sec.fluctuation} on fluctuation results). In particular, case $(a.)$ establishes bounds which demonstrate that $x$ will closely track the asymptotic behaviour of $H$ and case $(b.)$ establishes that when the forcing term, $H$, is sufficiently large $x(t) \sim H(t)$ as $t\to\infty$. Furthermore, when $x(t) \sim H(t)$ as $t\to\infty$, $x$ is of a strictly larger order of magnitude than the solution of the corresponding ODE \eqref{eq.introode}. This result allows us to pick up fluctuations in the solution even when $H$ is nonnegative. Even though the solution grows to infinity, it may not do so monotonically and the conclusion of Theorem \ref{Thm.det.L.H} identifies upper and lower rates of growth of the solution ($G_L H(t)$ and $G_U H(t)$ respectively) when $L_f(\gamma)\in(1,\infty)$. When $L_f(\gamma)=\infty$, the fluctuations are entirely determined by $H$.

The main results of this section are all proven by comparison arguments and the careful asymptotic analysis of the resulting differential inequalities. Since we assume positivity of $H$ to ensure asymptotic growth of solutions, it is straightforward to show that $\liminf_{t\to\infty}F(x(t))/\mu(\mathbb{R}^+) t \geq 1$; this is proven by a translation argument and appealing to \cite[Corollary 2]{sublinear2015}. The proof of the corresponding upper bound, $\limsup_{t\to\infty}F(x(t))/\mu(\mathbb{R}^+)t < \infty$, is more involved but can be roughly summarized as follows:

\begin{itemize}
				\setlength\itemsep{5pt}
	\item[Step 1:] Use monotonicity and finiteness of the measure to construct the crude upper inequality
	\begin{equation}\label{upper_discuss}
	x(t) < H_\epsilon(t) + (1+\epsilon)\mu(\mathbb{R}^+)\int_{T}^t \phi( x(s) )\,ds, \quad t \geq T,
	\end{equation}
	where $H_\epsilon$ includes constants and lower order terms, $\phi$ is a monotone function asymptotic to $f$ and we define $I_\epsilon(t) = \int_{T}^t \phi( x(s) )\,ds$ for $ t \geq T$.
	\item[Step 2:] Using hypotheses on the size of the perturbation term, show that $H_\epsilon$ is $o(I_\epsilon)$ or $O(I_\epsilon)$. 
	\item[Step 3:] Conclude the argument via a variation on Bihari's inequality.
\end{itemize}
Results in this section can be restated with positivity assumptions on $f$ and $H$ replaced by \eqref{asym_odd} below and 
\begin{equation}\label{H_pm}
H \in C(\mathbb{R}^+;\mathbb{R}).
\end{equation}
With this modification one obtains upper bounds on the rate of growth of solutions to \eqref{eq.xpert}, i.e. results of the type $\limsup_{t\to\infty}F(|x(t)|)/\mu(\mathbb{R}^+)t < \infty$.
\subsection{Fluctuation Results}\label{sec.fluctuation}
The existence of the limit $L_f(H)$ (even when it takes the value $+\infty$) is too strong a condition if we hope to apply our deterministic arguments to related equations with stochastic perturbations. We weaken the hypothesis $L_f(H) \in (0,\infty)$ as follows: assume that there exists a function $\gamma$ such that
\begin{align}\label{gamma_H}
\gamma \in C((0,\infty);(0,\infty)) \mbox{ is increasing and obeys } \lim_{t\to\infty}\gamma(t)=\infty \mbox{ and } \limsup_{t\to\infty}\frac{|H(t)|}{\gamma(t)} = 1.
\end{align}
As we no longer restrict ourselves to positive solutions, we ask for a degree of symmetry in the problem to simplify the analysis. We ask for ``asymptotic oddness'' of the nonlinearity in the following sense:
\begin{align}\label{asym_odd}
f \in C(\mathbb{R};\mathbb{R}) \mbox{ and }\lim_{|x|\to\infty}\frac{|f(x)|}{\phi(|x|)} = 1 \mbox{ for some } \phi \in C^1(\mathbb{R}^+;(0,\infty)).
\end{align} 
We now make hypotheses on $L_f(\gamma)$, as opposed to $L_f(H)$. We take $\limsup_{t\to\infty}|H(t)|/\gamma(t) = 1$, rather than positive and finite since we can always normalise this quantity while keeping the properties of $\gamma$ unchanged. Since $L_f(\gamma) \in (0,\infty)$ forces $\gamma$ to be eventually increasing, we suppose that $\gamma$ is always increasing for ease of exposition but there is strictly no need to make this assumption. Under \eqref{gamma_H} we can permit highly irregular behaviour in $H$ as long as we can capture some underlying regularity in the asymptotics of $H$ via a well-behaved auxiliary function, $\gamma$. For example, in applications to stochastic equations, $H$ could be a stochastic process whose partial maxima are described in terms of a deterministic function; this is the case for classes of processes obeying so-called iterated logarithm laws for instance. The following result illustrates the immediate utility of the hypothesis \eqref{gamma_H} for deterministic equations and furthermore details how this hypothesis carries over to the case when $L_f(\gamma)=\infty$.
\begin{theorem}\label{Thm.det.gamma}
Suppose \eqref{finite_measure}, \eqref{asym_odd}, \eqref{H_pm}, \eqref{fasym2} and \eqref{gamma_H} hold. Let $x$ denote the solution to \eqref{eq.xpert}.
\begin{itemize}
\item[(a.)] If $L_f(\gamma) \in (1,\infty)$, then
\[
\limsup_{t\to\infty}\frac{|x(t)|}{\gamma(t)} \in \left[0,\frac{L_f(H)}{L_f(H)-1}\right).
\]
\item[(b.)] If $L_f(\gamma) =\infty$, then
\[
\limsup_{t\to\infty}\frac{|x(t)|}{\gamma(t)} = 1, \quad \lim_{t\to\infty}\frac{x(t)-H(t)}{\gamma(t)} = 0.
\]
\end{itemize}
\end{theorem}
Case $(a.)$ of the result above indicates that when the perturbation is of intermediate size, in the sense that $L_f(\gamma) \in (1,\infty)$, solutions of \eqref{eq.xpert} are at most the same order of magnitude as $H$, modulo a multiplier. In case $(b.)$, when the perturbation is so large that $L_f(\gamma)=\infty$, solutions of \eqref{eq.xpert} have partial maxima of exactly the same order as those of $H$. This conclusion is strongly hinted at in case $(a.)$ of Theorem \ref{Thm.det.gamma} if one lets $L_f(\gamma) \to \infty$ in that result.

The restriction $L_f(\gamma) > 1$ is crucial to the proof of Theorem \ref{Thm.det.L.H} and cannot be relaxed within the framework of the current argument; we make this comment precise at the relevant moment during the proof itself (see Remark \ref{remark_L>1}). In fact, $L_f(\gamma)>1$ is not a purely technical contrivance but is also essential to the validity of our result. In Example \ref{example_L=1} we demonstrate that when $L_f(\gamma) \in (0,1]$ it is possible to have $\lim_{t\to\infty}|x(t)|/\gamma(t)=\infty$.

If $\limsup_{t\to\infty}|H(t)|/\gamma(t)=0$ in \eqref{gamma_H} we can use the following hypothesis and the arguments from Theorem \ref{Thm.det.gamma} to extend the scope of the result above:
\begin{align}\label{gamma_pm}
\limsup_{t\to\infty}\frac{|H(t)|}{\gamma_+(t)} =0, \quad \limsup_{t\to\infty}\frac{|H(t)|}{\gamma_-(t)} =\infty.
\end{align}
\begin{theorem}\label{Thm.det.gamma.pm}
Suppose \eqref{finite_measure}, \eqref{asym_odd}, \eqref{H_pm} and \eqref{fasym2} hold. Furthermore, suppose there exist increasing functions $\gamma_\pm \in C((0,\infty);(0,\infty))$ obeying $\lim_{t\to\infty}\gamma_\pm(t)=\infty$ such that \eqref{gamma_pm} holds and let $x$ denote the solution of \eqref{eq.xpert}.
\begin{itemize}
\item[(a.)] If $L_f(\gamma_\pm) \in (1,\infty]$, then
\[
\limsup_{t\to\infty}\frac{|x(t)|}{\gamma_+(t)} \in \left[0,\frac{1}{L_f(\gamma_+)}\right], \quad \limsup_{t\to\infty}\frac{|x(t)|}{\gamma_-(t)} = \infty.
\]
\item[(b.)] If $L_f(\gamma_\pm) = \infty$, then
\begin{align}\label{x_over_gammma_pm}
\lim_{t\to\infty}\frac{|x(t)|}{\gamma_+(t)} = 0, \quad \limsup_{t\to\infty}\frac{|x(t)|}{\gamma_-(t)} = \infty,
\end{align}
\end{itemize}
where it is understood that $1/L_f(\gamma_+) = 0$ if $L_f(\gamma_+)=\infty$.
\end{theorem}
In the presence of limited information about the behaviour of $H$, in the sense that \eqref{gamma_pm} holds, the result above tells us that the solution of \eqref{eq.xpert} is roughly the same order of magnitude as $H$, in the sense that $x$ also obeys \eqref{gamma_pm}, when $L_f(\gamma_\pm) = \infty$. When $L_f(\gamma_\pm)\in (1,\infty]$ we are still left with a weak conclusion and we are tempted to ask if this is an artifact of the method of proof. Example \ref{gamma_+} shows that we cannot expect to conclude that $\limsup_{t\to\infty}|x(t)|/\gamma_+(t)=0$ in general in this case. However, in attempting to apply this result it is likely that the user would actually seek to refine their choice of $\gamma_+$ in order to obtain a $\gamma_+$ obeying $L_f(\gamma_+)=\infty$ and hence make the stronger conclusion that $x \text{ is } o(\gamma_+)$.

Theorem \ref{Thm.det.gamma.pm} could equally well be stated as follows: $L_f(\gamma_+) \in (1,\infty]$ implies that \\ $\limsup_{t\to\infty}|x(t)|/\gamma_+(t) \leq 1/L_f(\gamma_+)$ and $L_f(\gamma_-) \in (1,\infty]$ implies $\limsup_{t\to\infty}|x(t)|/\gamma_+(t) = \infty$. These two statements are proved independently of one another but we chose to present them as part of a single result as we feel this is the manner in which they would prove most useful in practice; choosing $\gamma_+$ and $\gamma_-$ ``close together'' can give useful bounds on the size of the solution but using either bound in isolation only gives very crude information (see Example \ref{stoch_eg_1} for an illustration of this comment). 

The results of this section are proven via the usual machinery of comparison and asymptotic analysis but also rely crucially on the construction of a \emph{linear} differential inequality to achieve sharp results. The key steps in the argument can be understood as follows:
\begin{itemize}
\item[Step 1:] Using \eqref{upper_discuss}, derive the nonlinear differential inequality
\[
I_\epsilon '(t) < \phi\left(H_\epsilon(t) + \mu(\mathbb{R}^+)(1+\epsilon)I_\epsilon(t) \right),\quad t \geq T,
\]
where $I_\epsilon(t) = \int_T^t \phi(x(s))\,ds$. 
\item[Step 2:] Use \eqref{fasym2} to derive the \emph{linear} differential inequality
\begin{equation}\label{linear_inequality_discuss}
I_\epsilon '(t) < \phi(H_\epsilon(t)) + \frac{\phi(H_\epsilon(t))}{H_\epsilon(t)}\mu(\mathbb{R}^+)(1+\epsilon)^2 I_\epsilon(t), \quad t \geq T_1 > T.
\end{equation}
Since we can solve this inequality directly, there is no additional loss of sharpness here.
\item[Step 3:] Careful asymptotic analysis of the solution to the inequality \eqref{linear_inequality_discuss} using hypotheses on $L_f(H)$ yield upper bounds on the size of the solution to \eqref{eq.xpert}.
\item[Step 4:] The upper bounds achieved in Step 3 are recycled and further estimation yields the conclusions shown in the results above. 
\end{itemize}
Essentially the same steps outlined above are successful with random forcing, as we demonstrate below.
\subsection{Deterministic Examples} The following examples illustrate both the limitations and sharpness of some of the results outlined above. Consider the Volterra integro-differential equation given by
\[
x'(t) = \int_0^t e^{-(t-s)}f(x(s))\,ds + h(t), \quad t >0; \quad x(0) = \psi > 0.
\]
In the notation of \eqref{eq.volterra}, $m(t) = \int_0^t e^{-s}ds = 1 - e^{-t}$ and $\mu(\mathbb{R}^+)=1$. Hence
\begin{align}\label{example_1}
H(t) = x(t) -x(0) - \int_0^t f(x(s))\,ds + \int_0^t e^{-(t-s)}f(x(s))\,ds, \quad t \geq 0.
\end{align}
We construct examples by choosing a solution $x$, up to asymptotic equivalence, and then using \eqref{example_1} to figure out how large the perturbation term, $H$, must have been to generate a solution of this size. The calculations relevant to this section can be found in Appendix \ref{sec_examples}. For simplicity we forego any mention of hypotheses of the form \eqref{gamma_H} in this section and concentrate on the special case $\gamma = H$ with $H$ positive.
\begin{example}\label{golden}
The limits in Theorem \ref{Thm.det.L.F} are not always equal to $1$ or $1+L_f(H)$ and furthermore $L_f(H) \in (0,1]$ does not in general imply that $\lim_{t\to\infty}x(t)/H(t) = \infty$. If $f(x) = x^\beta$, $\beta \in (0,1)$, then
	\begin{align}\label{F_F_inv_2}
	F(x) \sim \frac{1}{1-\beta}x^{1-\beta} \text{ and } F^{-1}(x) \sim \left[(1-\beta)x\right]^{\tfrac{1}{1-\beta}}, \text{ as } x\to\infty.
	\end{align}
	Suppose $A \in [1,\infty)$ and take $x(t) = A[(1-\beta)t]^{\tfrac{1}{1-\beta}}$, for all $t \geq 0$. Thus $H(t) \sim (A-A^\beta)[(1-\beta)t]^{\tfrac{1}{1-\beta}}$ as $t \to \infty$. If $H(t) \sim [L_f(H)(1-\beta)t]^{\tfrac{1}{1-\beta}}$ as $t\to\infty$, then
	\[
	\lim_{t\to\infty}\frac{H(t)}{\mu(\mathbb{R}^+)\int_0^t f(H(s))ds} = L_f(H).
	\]
	Now suppose that $A-A^\beta = L_f(H)^{\tfrac{1}{1-\beta}}$ so we can choose an advantageous value of $L_f(H)$. For the purposes of this example it is sufficient to take $L_f(H)=1$ and $\beta = 1/2$. With these choices
	\[
	1 < \lim_{t\to\infty}\frac{F(x(t))}{\mu(\mathbb{R}^+)t} = \frac{1+\sqrt{5}}{2} \approx 1.618 < 2 = 1 + L_f(H),
	\]
	and the reader can compare this with the conclusion of Theorem \ref{Thm.det.L.F}. Finally, note that
	\[
	\lim_{t\to\infty}\frac{x(t)}{H(t)} = A \in (0,\infty).
	\]
\end{example}
\begin{example}\label{example_L=1}
If $f(x) = (x+e) / \log(x+e)$, then
\begin{align}\label{F_F_inv}
F(x) \sim \frac{1}{2}\log^2(x+e) \text{ and } F^{-1}(x) \sim e^{\sqrt{2x}} \text{, as } x \to\infty.
\end{align}
This example highlights the potential problems that emerge when one attempts to address the case $L_f(H)\in (0,1]$ (resp. $L_f(\gamma)$) in the context of Theorem \ref{Thm.det.gamma}. In particular, one cannot extend the conclusion of Theorem \ref{Thm.det.gamma} to cover $L_f(H) \in (0,1]$ without additional hypotheses because when $L_f(H) \in (0,1]$ it is possible to have  $\lim_{t\to\infty}x(t)/H(t)=\infty$.

Choose $x(t) = \exp\left( \lambda(t) + \sqrt{2(t+1)} \right)-e = \exp(P(t))-e$ for $t \geq 0$ and let $\lambda(t) = (1+t)^\alpha$ for some $\alpha \in (0,1/2)$. In this case $H(t) \sim K\,P(t)^{2\alpha-1}\exp(P(t))$. Furthermore, $H$ obeys $L_f(H)=1$ and by construction $\lim_{t\to\infty}x(t)/H(t) = \infty$. However, we still have $\lim_{t\to\infty}F(x(t))/\mu(\mathbb{R}^+)t = 1.$ 
\end{example}
\begin{example}\label{example_L>1}
The bounds on $\lim_{t\to\infty}x(t)/H(t)$ and $\liminf_{t\to\infty}F(x(t))/\mu(\mathbb{R}^+)t$ obtained in Theorems \ref{Thm.det.L.F} and \ref{Thm.det.L.H} can be attained. Once more suppose that $f(x) = (x+e) / \log(x+e)$.

Suppose $L_f(H) \in (1,\infty)$ and choose $x(t) = \exp\left( \sqrt{2L_f(H)(t+1)} \right) - e$ for $t \geq 0$. This gives 
$
H(t) \sim ((L_f(H)-1)/L_f(H))\exp\left(\sqrt{2L_f(H)(t+1)}\right)
$
as $t\to\infty$ and 
\[
\lim_{t\to\infty}\frac{H(t)}{\mu(\mathbb{R}^+)\int_0^t f(H(s))ds} = L_f(H) \in (1,\infty).
\]
Hence $\lim_{t\to\infty}x(t)/H(t) = L_f(H)/(L_f(H)-1)$, achieving the upper bound predicted by Theorem \ref{Thm.det.L.H}. Futhermore, a simple calculation reveals that $\lim_{t\to\infty}F(x(t))/\mu(\mathbb{R}^+)t = 1$, achieving the lower bound from Theorem \ref{Thm.det.L.F}.
\end{example}
\begin{example}\label{example_L=infty}
We present a simple example illustrating the case when the solution to \eqref{eq.xpert} is asymptotic to $H$ and the functional $L_f(H)$ takes the value $+\infty$. Let $f(x) = (x+e) / \log(x+e)$.	
	
Suppose $x(t) = \exp\left([2(t+1)]^\alpha \right)-e$, $\alpha \in \left(1/2,1\right)$, $t \geq 0$. It follows easily from equation \eqref{example_1} that $x(t) \sim H(t)$ as $t\to\infty$ and hence that $L_f(H)=\infty$ as in Theorems \ref{Thm.det.gamma} and \ref{Thm.det.L.H}, case $(b.)$.
\end{example}
\begin{example}\label{gamma_+}
In case $(a.)$ of Theorem \ref{Thm.det.gamma.pm} it is possible to have $\limsup_{t\to\infty}H(t)/\gamma_+(t)=0$ but\\ $\limsup_{t\to\infty}x(t)/\gamma_+(t)>0$ when $L_f(\gamma_+) \in (1,\infty)$. Hence there is no straightforward improvement of the conclusion of Theorem \ref{Thm.det.gamma.pm} when $L_f(\gamma_+) \in (1,\infty)$.

Let $f(x)=x^\beta$ with $\beta\in(0,1)$, $H=0$, and $\gamma_+(t)=F^{-1}(\alpha\, \mu(\mathbb{R}^+) t)$ with $\alpha \in (1,\infty)$. This implies that $x(t) \sim F^{-1}(\mu(\mathbb{R}^+)t)$ as $t\to\infty$, where the asymptotics of $F^{-1}$ are given by \eqref{F_F_inv_2}, and hence $\lim_{t\to\infty}x(t)/\gamma_+(t) = \alpha^{-1/(1-\beta)} > 0$, as required. It is straightforward to verify that $L_f(\gamma_+)=\alpha \in (1,\infty)$. 
\end{example}
\section{Stochastic Volterra Equations}\label{stochastic}
We now study the pathwise asymptotic behaviour of solutions to \eqref{eq.introsfde}. Our approach is to treat \eqref{eq.introsfde} as a perturbed version of \eqref{eq.xpert} where the forcing term is now stochastic and hence to leverage our deterministic results as much as possible. After proving a strong existence theorem for solutions to \eqref{eq.introsfde}, we use the pathwise asymptotic theory for continuous Brownian martingales and $\alpha$--stable L\'evy processes to show that the main results from the previous section are sufficiently general that we can extend them to provide asymptotic estimates on the pathwise growth and fluctuation of solutions. We also explain how our results provide a programme for establishing similar pathwise bounds for broader classes of admissible stochastic noise.

We henceforth work on a given probability space $(\Omega, \mathcal{F},\mathbb{P},(\mathcal{F}_t)_{t\geq 0})$ which is complete and has a right continuous filtration. We ask that the nonlinear function $f: \mathbb{R} \mapsto \mathbb{R}$ obeys the following local Lipschitz condition: for each $d>0$ there exists $K_d>0$ such that
\begin{equation}\label{f_local_lipschitz}
|f(x)-f(y)| \leq K_d\, |x-y|, \mbox{ for each } x \mbox{ and }y \in [-d,d], 
\end{equation}
and that $f$ obeys a global linear bound of the form
\begin{equation}\label{f_global_linear}
|f(x)| \leq K + \eta \,|x|, \mbox{ for each } x \in \mathbb{R},
\end{equation}
where $K$ and $\eta$ are positive constants.

In order to leverage the framework of M\'etivier and Pellaumail \cite{metivier1980stopped} we make a slight modification to the formulation of \eqref{eq.introsfde} and consider the stochastic integral equation 
\begin{equation}\label{volterra_semi_1}
X(t) = X(0) + \int_0^t\left(\int_{(0,s]} \mu(du)f(X(s-u)) + \mu(\{0\})f(X(s-))\right)ds + Z(t), \quad t \geq 0.
\end{equation}
By applying Fubini's Theorem and making a suitable change of variable, \eqref{volterra_semi_1} can be written as
\begin{equation}\label{volterra_semi_2}
X(t) = X(0) + \mu(\{0\})\int_{[0,t]} f(X(s-))ds + \int_{[0,t)} m_{-}(t-s) f(X(s))\,ds + Z(t), \quad t \geq 0,
\end{equation}
where $X(t-) = X(\lim_{s\uparrow t})$ and $m_-(t) = \mu((0,t])$. This adjustment is necessary for the functional
\begin{equation}\label{volterra_functional}
a(s,\omega,X)=\int_{(0,s]} \mu(du)f(X(s-u)) + \mu(\{0\})f(X(s-)), \quad s \geq 0,
\end{equation}
to define a predictable process (measurable with respect to the filtration generated by adapted, left continuous processes) and hence be integrable with respect to general semimartingales (see Protter \cite{protter2004stochastic} for details).

In order to define the notion of a strong solution for stochastic equations such as \eqref{volterra_semi_2}, we recall some standard terminology from the theory of stochastic processes: a regular process is one which is adapted and has right continuous paths with left hand limits (RCLL). A process $X$ is called $\mathbb{P}$--null if $\mathbb{P}$ almost surely the paths $t \mapsto X(t)$ are identically zero functions.
\begin{definition}\label{defn_strong_solution}
	A process $X$ defined on $[0,\tau)$ is said to be a strong solution to equation \eqref{volterra_semi_2} on $[0,\tau)$ with initial value $X(0)$ if the process \[
	\mu(\{0\})\int_{[0,t]} f(X(s-))ds + \int_{[0,t)} m_{-}(t-s) f(X(s))\,ds + Z(t)
	\] is well--defined on $[0,\tau)$ as a regular process and differs from $X(t)-X(0)$ by a $\mathbb{P}$--null process.
\end{definition}
The solution to \eqref{volterra_semi_2} is unique if for any two processes $X$ and $Y$ obeying Definition \ref{defn_strong_solution}, $X-Y$ is a $\mathbb{P}$--null process.
\begin{theorem}\label{Thm.existence}
Let \eqref{finite_measure} hold and let $Z$ be a cadlag semimartingale. If $f: \mathbb{R}\mapsto \mathbb{\mathbb{R}}$ is measurable and obeys \eqref{f_local_lipschitz}, and \eqref{f_global_linear}, then there exists a unique, strong solution to  \eqref{volterra_semi_2}. Moreover, this solution exists for all $t \geq 0$ as a real-valued process with probability 1.
\end{theorem}
\begin{proof}
This theorem is a natural specialisation of a result of M\'etivier and Pellaumail \cite[Theorem 5]{metivier1980stopped}. In order to apply the aforementioned result we must check that the functional from \eqref{volterra_functional} and also the constant functional $a(s,\omega,X) = 1$ obey the following pair of conditions:
firstly for any regular processes (adapted with cadlag paths) $X$ and $Y$, for each $d>0$ there exists a constant $L_d>0$ such that
\begin{align}\label{functional_local_lipschitz}
|a(t,\omega,X) - a(t,\omega,Y)| \leq L_d \sup_{0 \leq s < t}| X(s)-Y(s) |
\end{align} 
for each $t \in \mathbb{R}^+$, $\sup_{0 \leq s < t}| X(s)| \leq d$ and $\sup_{0 \leq s < t}| Y(s)| \leq d$. Secondly, for any regular process $X$ there exists $C>0$ such that
\begin{equation}\label{global_linear_bound}
|a(t,\omega,X)| \leq C \sup_{0 \leq s < t}\left(| X(s)| + 1\right)
\end{equation}
for each $t \in \mathbb{R}^+$.

When the functional $a$ is constant the conditions above are trivially satisfied so suppose now that $a$ is given by \eqref{volterra_functional} and proceed to verify condition \eqref{functional_local_lipschitz}. Let $X$ and $Y$ be any two regular processes satisfying $\sup_{0 \leq s < t}|X(s)| \leq d$ (resp. $Y$), fix $t \in \mathbb{R}^+$ and estimate as follows:
\begin{align*}
|a(t,\omega,X) - a(t,\omega,Y)| &\leq \mu(\{0\})|f(X(t-)) - f(Y(t-))| + \int_{(0,t]} \mu(ds) |f(X(t-s)) - f(Y(t-s))|\\
&\leq \mu(\mathbb{R}^+)\,K_d, \sup_{0 \leq s < t}|X(s) - Y(s)|, 
\end{align*}
where we have used both \eqref{finite_measure} and \eqref{f_local_lipschitz}. Now check \eqref{global_linear_bound}; assume $X$ is a regular process and fix $t \in \mathbb{R}^+$. The following inequality is a straightforward consequence of \eqref{finite_measure} and \eqref{f_global_linear}:
\begin{align*}
|a(t,\omega,X)| \leq \mu(\{0\})|f(X(t-))| + \int_{(0,t]} \mu(ds) |f(X(t-s))|
\leq C^*\,\sup_{0 \leq s < t}\left(  |X(s)| + 1  \right),
\end{align*}
where $C^* = \mu(\mathbb{R}^+)K$. Thus there exists a unique, strong solution to \eqref{volterra_semi_2}, at least locally in time (i.e. we could have $\tau<\infty$ with positive probability in Definition \ref{defn_strong_solution}). Moreover, \eqref{f_global_linear} can be used to bound the solution to \eqref{volterra_semi_2} below the solution to the corresponding linear equation, which exists for all $t \geq 0$ with probability 1. Hence finite-time blow-up occurs with probability zero and the second part of the claim is proven.
\end{proof}
\begin{remark}
The condition \eqref{f_global_linear} will always be satisfied in this section since the hypotheses \eqref{asym_odd} and \eqref{fasym2} will be imposed throughout. The assumption \eqref{finite_measure} is also present throughout so the only additional hypothesis imposed by Theorem \ref{Thm.existence} is that of local Lipschitz continuity on the nonlinear function $f$ (which is only required to guarantee uniqueness of solutions).
\end{remark}
We pause now to consider the method by which the results of this section are proven and to illustrate that this presents a framework for generating similar pathwise asymptotic results for a wide range of suitable stochastic forcing terms. Our method of proof relies principally on building appropriate comparison equations so we are not concerned about the pathwise regularity of the solution to \eqref{eq.introsfde} and hence can treat quite irregular forcing processes.

If $Z = \{Z(t)\,:\,t\geq 0\}$ denotes the forcing term in \eqref{eq.introsfde}, then our general approach is as follows:
\begin{itemize}
				\setlength\itemsep{5pt}
	\item[Step 1:] Prove pathwise asymptotic bounds on the size of the process $Z$ in terms of a well-behaved deterministic function, say $\gamma$, on which we can formulate functional hypotheses in terms of $L_f(\cdot)$. These bounds should be in the spirit of \eqref{gamma_H} or \eqref{gamma_pm}. 
	\item[Step 2:] Construct an upper comparison solution (pathwise) in terms of $\gamma$ which majorizes the solution to the \eqref{eq.introsfde}; this essentially reduces the stochastic problem to a deterministic one.
	\item[Step 3:] Conclude the argument using suitable hypotheses on $L_f(\gamma)$ and the results of Section \ref{det.results}.
\end{itemize}
The steps outlined above also highlight how one can approach the important, and more challenging, case of state-dependent noise, i.e. equations of the form
\[
dX(t) = \int_{[0,t]} \mu(ds) f(X(t-s))\,dt +G\left(t,\,X(t)\right)dZ(t), \quad t>0.
\]
In particular, hypotheses involving the functional $L_f(\cdot)$ will still allow completion of Step 1, but now the function $\gamma$ will involve the process $X$ itself. Thus Step 2 will now necessitate further analysis of how this state-dependent perturbation interacts with the drift term $\int_{[0,t]} \mu(ds) f(X(t-s))$.
\subsection{Brownian Noise}\label{sec.brownian}
Throughout this section $X$ denotes the unique, strong solution to \eqref{volterra_semi_2},
\begin{equation}\label{sigma_con}
Z(t) = \int_0^t \sigma(s)\,dB(s), \mbox{ where }B \mbox{ is a standard 1-D Brownian motion} \mbox{ and }  \sigma \in C(\mathbb{R}^+, \mathbb{R}), 
\end{equation}
and 
\begin{align*}
\Sigma(t) = \sqrt{ 2 \left(\int_0^t \sigma^2 (s)\, ds\right) \log\log\left( \int_0^t \sigma^2 (s) \,ds \right)}, \quad t \geq 0.
\end{align*}
Analogous to the deterministic case, we classify the behaviour of solutions to \eqref{eq.introsfde} according to whether the number $L_f(\Sigma)$ is zero, finite or infinite.

The existence and uniqueness of solutions of \eqref{eq.introsfde} is naturally simpler in the case of Brownian noise. In particular, there is a unique, continuous (strong) solution to \eqref{eq.introsfde} with Brownian noise if \eqref{finite_measure} holds and the nonlinearity is locally Lipschitz continuous with a global linear bound (see Mao \cite[Ch. 5]{mao2007stochastic}).

\begin{figure}
\centering
\includegraphics[width=\textwidth]{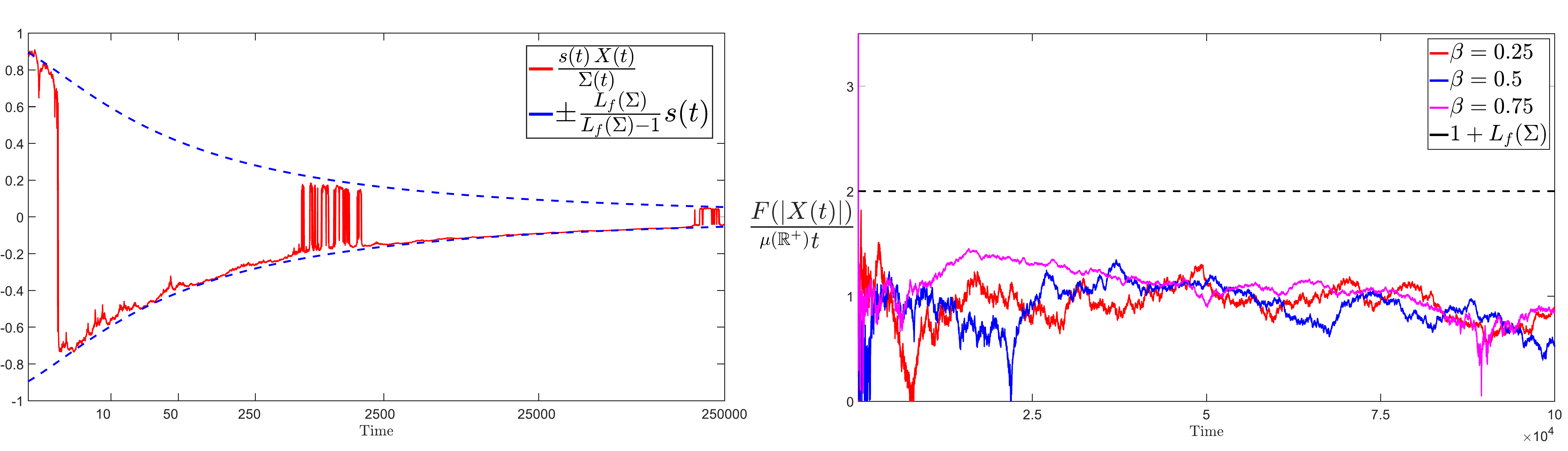}
\caption{\textbf{Left:} $s(t) = (1+t)^{-3}$ and thus according to Theorem \ref{Thm.stoch.L.H}, the quantity $s(t)X(t)/\Sigma(t)$ fluctuates between the bounds $\pm s(t) L_f(\Sigma)/(L_f(\Sigma)-1)$ as it tends to zero a.s. (for $L_f(\Sigma) >1$). This plot numerically illustrates these fluctuations and shows that the quantitative bounds given on the fluctuations appear to be approximately sharp. \textbf{Right:} Theorem \ref{Thm.stoch.L.F} provides an almost sure bound on the quantity $F(|X(t)|)/\mu(\mathbb{R}^+)t$ for large $t$. In this plot we observe that the linear bound in terms of $L_f(\Sigma)$ indeed holds but does not appear to be sharp in general. See Appendix \ref{sec_examples} for details.}
\label{fig_sharpness}
\end{figure}

When formulating functional conditions on \eqref{eq.introsfde} to preserve growth of the type \eqref{eq.unpertgrowth} it is necessary to distinguish between the cases $\sigma \in L^2(0,\infty)$ and $\sigma \notin L^2(0,\infty)$. When $\sigma \in L^2(0,\infty)$ the martingale term in \eqref{eq.introsintde}, $\int_0^t \sigma(s)dB(s)$, will tend to an a.s. finite random variable and in this case we clearly expect to retain the growth rate of solutions to \eqref{eq.introode}. However, when $\sigma \notin L^2(0,\infty)$ the martingale term is recurrent on $\mathbb{R}$ and has large fluctuations of order $\Sigma(t)$~\cite{revuzandyor}. Our first result shows that when $\sigma \notin L^2(0,\infty)$ and $L_f(\Sigma)=0$, the solution to \eqref{eq.introsfde} cannot grow faster than that of the ordinary differential equation \eqref{eq.introode}.
\begin{theorem}\label{Thm.stoch.epsilon}
Let \eqref{finite_measure}, \eqref{asym_odd}, \eqref{fasym}, and \eqref{sigma_con} hold. Suppose additionally that $\lim_{x\to\infty}f(x)=\infty$. If
\begin{align}\label{small_pert_stochastic_F}
\lim_{t\to\infty}\frac{\Sigma(t)}{F^{-1}\left(\mu(\mathbb{R}^+)(1+\epsilon)t\right)} = 0,\,\, \mbox{for each } \epsilon \in (0,1),
 \end{align}
then
\[
\limsup_{t\to\infty}\frac{F(|X(t)|)}{\mu(\mathbb{R}^+)t} \leq 1 \,\, a.s.
\]
\end{theorem}
Our next result mirrors the conclusion of Theorem \ref{Thm.stoch.epsilon} but uses hypotheses on the functional $L_f(\Sigma)$ instead of the condition \eqref{small_pert_stochastic_F}.
\begin{theorem}\label{Thm.stoch.zero}
Let \eqref{finite_measure}, \eqref{asym_odd}, \eqref{fasym}, and \eqref{sigma_con} hold. Suppose additionally that $\lim_{x\to\infty}f(x)=\infty$.
\begin{itemize}
\item[(a)] If $\sigma \notin L^2(0,\infty)$ and $L_f(\Sigma)=0$, then
\[
\limsup_{t\to\infty}\frac{F(|X(t)|)}{\mu(\mathbb{R}^+)t} \leq 1 \,\, a.s.
\]
\item[(b)] If $\sigma\in L^2(0,\infty)$, then $L_f(\Sigma)=0$ and 
\[
\limsup_{t\to\infty}\frac{F(|X(t)|)}{\mu(\mathbb{R}^+)t} \leq 1 \,\, a.s.
\]
\end{itemize}
\end{theorem}
An interesting special case of Theorem \ref{Thm.stoch.zero}, which is likely to be important in applications, is when the function $\sigma$ is a nonzero constant. In this case, solutions to \eqref{eq.introsfde} are unbounded with probability one. 
\begin{corollary}\label{sigma_const}
Let \eqref{finite_measure}, \eqref{asym_odd}, \eqref{fasym}, and \eqref{sigma_con} hold. Suppose additionally that $\lim_{x\to\infty}f(x)=\infty$. If $\sigma(t) = \sigma \in \mathbb{R}/\{0\}$ for all $t \geq 0$, then
\[
\limsup_{t\to\infty}|X(t)| = \infty\,\, a.s. \mbox { and }\limsup_{t\to\infty}\frac{F(|X(t)|)}{\mu(\mathbb{R}^+)t} \leq 1 \,\, a.s.
\]
\end{corollary}
As in the deterministic case when the perturbation is of intermediate or critical magnitude, in the sense that $L_f(\Sigma) \in (0,\infty)$, we expect the solution to inherit characteristics of both the perturbation and the ordinary differential equation \eqref{eq.introode}. Our next result demonstrates that this is indeed the case by showing that if the solution to \eqref{eq.introsfde} grows, then its growth rate is at most the same order of magnitude as the solution to \eqref{eq.introode}.
\begin{theorem}\label{Thm.stoch.L.F}
Let \eqref{finite_measure}, \eqref{asym_odd}, \eqref{fasym2} and \eqref{sigma_con} hold. Suppose additionally that $\lim_{x\to\infty}f(x)=\infty$ and $\sigma \notin L^2(0,\infty)$. If  $L_f(\Sigma) \in (0,\infty)$, then
\[
\limsup_{t\to\infty}\frac{F(|X(t)|)}{\mu(\mathbb{R}^+)t} \leq 1+L_f(\Sigma) \,\, a.s.
\]
\end{theorem}
When $L_f(\Sigma)\in(1,\infty)$ we show that if the the solution to \eqref{eq.introsfde} fluctuates, then these fluctuations are at most of order $\Sigma(t)$ times a multiplier which we can bound in terms of $L_f(\Sigma)$. As in Theorem \ref{Thm.det.L.H} we are unable to extend this argument to $L_f(\Sigma) \in (0,1)$ for technical reasons which become apparent in the relevant construction. While this bound is practically useful, it appears that it is not sharp in general (see Figure \ref{fig_sharpness}, right).

Nonnegativity of the measure $\mu$ no longer plays an important role in the results above; primarily because we are reduced to proving upper bounds on the growth rate of solutions once solutions are no longer necessarily of one sign. For ease of exposition we have left the hypothesis \eqref{finite_measure} in place but it could equally well be replaced by the hypothesis that $\mu$ is a Borel measure with finite total variation norm, i.e. $|\mu| = \mu(\mathbb{R}^+) \in (0,\infty)$, with the results above unchanged. 
\begin{theorem}\label{Thm.stoch.L.H}
Let \eqref{finite_measure}, \eqref{asym_odd}, \eqref{fasym2} and \eqref{sigma_con} hold. Suppose additionally that $\lim_{x\to\infty}f(x)=\infty$ and $\sigma \notin L^2(0,\infty)$. If $L_f(\Sigma) \in (1,\infty)$, then
\[
 \frac{-L_f(\Sigma)}{L_f(\Sigma)-1} \leq \liminf_{t\to\infty}\frac{X(t)}{\Sigma(t)} \leq \limsup_{t\to\infty}\frac{X(t)}{\Sigma(t)} \leq \frac{L_f(\Sigma)}{L_f(\Sigma)-1} \,\, a.s.
\]
\end{theorem}
Under the hypotheses of Theorem \ref{Thm.stoch.L.H} we can additionally conclude that 
\[
\liminf_{t\to\infty}\frac{X(t)}{\Sigma(t)} \leq \frac{2-L_f(\Sigma)}{L_f(\Sigma)-1}\,\, a.s., \quad \limsup_{t\to\infty}\frac{X(t)}{\Sigma(t)} \geq \frac{L_f(\Sigma)-2}{L_f(\Sigma)-1} \,\, a.s.
\]
Hence, when $L_f(\Sigma)>2$, $X(t)$ is recurrent on $\mathbb{R}$. This leaves open the question of recurrence, or in other words, whether or not the process actually fluctuates, for $L_f(\Sigma)\in(1,2)$. In Figure \ref{fig_sharpness} (left) we show  simulations of the process in which (a scaled version of) the quantity $X(t)/\Sigma(t)$ fluctuates between the quantitative bounds predicted by Theorem \ref{Thm.stoch.L.H}. These numerical experiments further illustrate that our bounds appear to be approximately sharp (at least for some classes of examples).

Finally, when the perturbation term is so large that $L_f(\Sigma)=\infty$ we expect this exogenous force to dominate the system and this intuition is confirmed by our next result. In particular, we prove that the solution to \eqref{eq.introsfde} is recurrent on $\mathbb{R}$ and that its fluctuations are precisely of order $\Sigma$.
\begin{theorem}\label{Thm.stoch.infty}
Let \eqref{finite_measure}, \eqref{asym_odd}, \eqref{fasym2} and \eqref{sigma_con} hold. Suppose additionally that $\lim_{x\to\infty}f(x)=\infty$ and $\sigma \notin L^2(0,\infty)$. If $L_f(\Sigma)=\infty$, then
\[
\liminf_{t\to\infty}\frac{X(t)}{\Sigma(t)} = -1 \,\, a.s. \mbox{ and }\limsup_{t\to\infty}\frac{X(t)}{\Sigma(t)} = 1 \,\, a.s.,
\]
and furthermore
\begin{align}\label{X_over_sigma_to_zero}
\lim_{t\to\infty}\frac{X(t) - \int_0^t \sigma(s)dB(s)}{\Sigma(t)} = 0 \,\, a.s.
\end{align}
\end{theorem}
\subsection{L\'evy Noise}
We now assume that the semimartingale $Z$ in \eqref{eq.introsfde} is an $\alpha$--stable L\'evy process; the results which follow further emphasize the fact that our methods do not rely on the path continuity of the process in any essential way. For the readers convenience we recall the relevant definitions from the theory of L\'evy processes.
\begin{definition}
If $Z = \{Z(t) \,:\, t \geq 0\}$ is a L\'evy process, it's characteristic function $\mathscr{F}_Z$ is given by
\[
\mathscr{F}_Z(\lambda) = e^{- \Psi(\lambda)}, \quad \lambda \in \mathbb{R},
\]
where $\Psi : \mathbb{R} \mapsto \mathbb{C}$ is of the form
\begin{equation}\label{char_exp}
\Psi(\lambda) = i a \lambda + \frac{1}{2}\sigma^2 \lambda^2 + \int_{\mathbb{R}}\left(1 - e^{i x  \lambda } + i x \lambda \mathbbm{1}_{\{|x|<1\}} \right)\Pi(dx), 
\end{equation}
with $a \in \mathbb{R}$, $\sigma \in \mathbb{R^+}$ and $\Pi$ 
a measure on $\mathbb{R} /\{0\}$ satisfying $\int_{\mathbb{R}}(1 \wedge |x|^2)\Pi(dx) < \infty$.

$\Psi$ is called the \emph{characteristic exponent} of the process $Z$.
\end{definition}
The number $a$ in \eqref{char_exp} corresponds to the linear ``drift'' coefficient of the L\'evy process in question, $\sigma$ is called the Gaussian coefficient and corresponds to the Brownian or continuous random component; $\Pi$ is called the L\'evy measure and represents the pure jump part of the process. A L\'evy process is uniquely specified by the triple $(a,\sigma,\Pi)$.
\begin{definition}
	For each $\alpha \in (0,2]$, a L\'evy process with characteristic exponent $\Psi$ is called a stable process with index $\alpha$ ($\alpha$--stable for short) if $\Psi(k\lambda) = k^{\alpha}\Psi(\lambda)$ for each $k>0,\, \lambda \in \mathbb{R}^d$.
\end{definition}
Stable processes are closely related to the class of stable distributions which gain their importance as ``attractors'' for normalised sums of independent and identically distributed random variables. In particular, a sum of random variables with power law decay in the tails, proportional to $|x|^{-1-\alpha}$, will tend to a stable distribution if $0<\alpha<2$ and to a normal distribution if $\alpha \geq 2$. Integrability of the L\'evy measure forces us to consider $\alpha \in (0,2]$ and in this section we also ignore the case $\alpha = 2$ since this corresponds to the case of Brownian noise (which was considered in detail in Section \ref{sec.brownian}). We tacitly exclude the degenerate case when $Z$ is a pure drift process and assume for the remainder of this section that
\begin{equation}\label{alpha_stable}
Z \mbox{ is an $\alpha$--stable process with }\alpha \in (0,2).
\end{equation}
Let $X$ denote the unique, strong solution to \eqref{eq.introsfde} throughout.

Our first result is a stochastic analogue of Theorem \ref{Thm.det.zero} and provides a sufficient condition to retain growth to infinity no faster than the solution of \eqref{eq.introode} in the presence of $\alpha$--stable noise.
\begin{theorem}\label{Thm.Levy0}
Let \eqref{finite_measure}, \eqref{asym_odd}, \eqref{fasym} and \eqref{alpha_stable} hold. If  $\lim_{x\to\infty}f(x)=\infty$ and there exists an increasing function $\gamma \in C ((0,\infty);(0,\infty))$ such that $L_f(\gamma) =0$ and $\int_0^\infty {\gamma(s)}^{-\alpha}ds < \infty$, then
	\[
	\limsup_{t\to\infty}\frac{F(|X(t)|)}{\mu(\mathbb{R}^+)t} \leq 1 \mbox{ a.s.}
	\]
\end{theorem}

The next results provides a direct stochastic analogue of Theorem \ref{Thm.det.gamma.pm}.
\begin{theorem}\label{Thm.Levy1}
Let \eqref{finite_measure}, \eqref{asym_odd}, \eqref{fasym2} and \eqref{alpha_stable} hold. Suppose additionally that $\lim_{x\to\infty}f(x)=\infty$ and $\gamma \in C ((0,\infty);(0,\infty))$ is an increasing function such that $L_f(\gamma) \in (1,\infty]$. If $\int_0^\infty {\gamma(s)}^{-\alpha}ds < \infty$, then
\[
\limsup_{t\to\infty}\frac{|X(t)|}{\gamma(t)} \leq \frac{1}{L_f(\gamma)} \mbox{ a.s.},
\]
where we interpret $1/L_f(\gamma) = 0$ if $L_f(\gamma) = \infty$. If $\int_0^\infty {\gamma(s)}^{-\alpha}ds = \infty$, then 
\[
\limsup_{t\to\infty}\frac{|X(t)|}{\gamma(t)} = \infty \mbox{ a.s.}
\]
\end{theorem}
\subsection{Stochastic Examples}\label{sec.stoch_examples}
\begin{example}\label{stoch_eg_1}
To illustrate the practical application of the results in Section \ref{sec.brownian} we present an example with power type nonlinearity and Brownian noise, i.e. $Z(t) = \int_0^t \sigma(s)\,dB(s)$. Suppose 
\[
f(x) = \sign(x)|x|^\beta, \quad x \in \mathbb{R}, \quad \beta \in (0,1),
\]
$\sigma(t) = t^\alpha$, $t \geq 0$, for some $\alpha>0$, and $\mu$ is a measure obeying \eqref{finite_measure}. In this framework 
\begin{equation}\label{Sigma_power_eg}
\Sigma(t) \sim t^{\alpha + 1/2} A(t,\alpha) \mbox{ as }t\to\infty, \mbox{ where } A(t,\alpha)=\sqrt{\frac{2 \log\log t}{2\alpha+1} },
\end{equation}
and 
\[
F(x) \sim \frac{1}{1-\beta}x^{1-\beta} \mbox{ as }x\to\infty.
\]
Clearly, $\Sigma(t) \to \infty$ as $t \to\infty$ and therefore
$
L_f(\Sigma) = \lim_{t\to\infty}\Sigma'(t)/\mu(\mathbb{R}^+)f(\Sigma(t))$.
It is straightforward to show that
\[
\Sigma'(t) = t^{\alpha - 1/2} \left(\frac{2}{2\alpha+1}\right)^{-1/2} \left(\log\log\left(\frac{t^{2\alpha+1}}{2\alpha+1}\right)\right)^{1/2} \left\{ 1 + \left(\log\left(\frac{t^{2\alpha+1}}{2\alpha+1}\right) \log\log\left(\frac{t^{2\alpha+1}}{2\alpha+1}\right) \right)^{-1} \right\},
\]
for $t \geq 0$ and hence
\[
L_f(\Sigma) = 
\begin{cases}
0, \quad &0 < \alpha < (1+\beta)/2(1-\beta),\\
\infty,  &\alpha \geq (1+\beta)/2(1-\beta).
\end{cases}
\]
By Theorem \ref{Thm.stoch.zero} we conclude that the unique, strong solution of \eqref{eq.introsfde} obeys
\[
\limsup_{t\to\infty}\frac{F(|X(t)|)}{\mu(\mathbb{R}^+)t} = \limsup_{t\to\infty}\frac{|X(t)|^{1-\beta}}{\mu(\mathbb{R}^+)(1-\beta)t} \leq 1 \mbox{ a.s. for }0 < \alpha < \frac{1+\beta}{2(1-\beta)}.
\]
Similarly, by Theorem \ref{Thm.stoch.infty}, 
\[
\liminf_{t\to\infty}\frac{X(t)}{A(t,\alpha)\, t^{\alpha + 1/2}} = -1 \mbox{ a.s. and } \limsup_{t\to\infty}\frac{X(t)}{A(t,\alpha)\, t^{\alpha + 1/2}} = 1 \mbox{ a.s. for } \alpha \geq\frac{1+\beta}{2(1-\beta)},
\]
where the function $A(t,\alpha)$ is given by \eqref{Sigma_power_eg}. 
\end{example}
\begin{example}\label{stoch_eg_2}
Let $Z$ be an $\alpha$--stable process with index $\alpha \in (0,2)$ and, as in the previous example, suppose we have a power--type nonlinearity given by
\[
f(x) = \sign(x)|x|^\beta, \quad x \in \mathbb{R}, \quad \beta \in (0,1).
\]
Let $\mu$ be a measure obeying \eqref{finite_measure} and suppose the function $\gamma_+$ is given by 
\[
\gamma_+(t) = (1+t)^{\epsilon}, \quad t \geq 0, \quad \epsilon> \frac{1}{\alpha} > 0.
\]
By construction, $\gamma_+$ is increasing, positive and satisfies $\int_0^\infty {\gamma_+(t)}^{-\alpha}\,dt<\infty$. Furthermore,
\[
L_f(\gamma_+) = 
\begin{cases}
0, \quad &1/\alpha < \epsilon < 1/(1-\beta),\\
\epsilon/\mu(\mathbb{R}^+), & 1/\alpha < \epsilon = 1/(1-\beta),\\
\infty, &\epsilon > \max\left( 1/\alpha, 1/(1-\beta) \right).
\end{cases}
\]
If the interval $(1/\alpha,\,1/(1-\beta))$ is nonempty, then we can take $\gamma$ in the statement of Theorem \ref{Thm.Levy0} to be $\gamma_+$ with $\epsilon \in (1/\alpha,\,1/(1-\beta))$. Hence the solution of \eqref{eq.introsfde} obeys
\[
\limsup_{t\to\infty}\frac{F(|X(t)|)}{\mu(\mathbb{R}^+)t} \leq 1 \mbox{ a.s., when }\beta > 1 - \alpha.
\]
This essentially means that if the nonlinearity is sufficiently strong we cannot experience growth in the solution of \eqref{eq.introsfde} faster than that seen in \eqref{eq.introode} (with positive probability). The restriction $\beta > 1 - \alpha$ is intuitive in the following sense: the smaller $\alpha$ is, the more mass there is in the tail of the L\'evy measure associated with $Z$ and hence the partial maxima of $Z$ will tend to grow faster the smaller the value of $\alpha$; when $\alpha$ is small we require a stronger nonlinearity (larger value of $\beta$) to retain the unperturbed growth rate. When $\alpha \geq 1$ we always retain the growth rate of the unperturbed equation.

If we take $\epsilon = 1/(1-\beta)$, then $L_f(\gamma_+)=1/\mu(\mathbb{R}^+)(1-\beta)$ and we can apply Theorem \ref{Thm.Levy1} to yield
\begin{equation}\label{example_est_1}
\limsup_{t\to\infty}\frac{|X(t)|}{t^{1/(1-\beta)}} \leq \mu(\mathbb{R}^+)(1-\beta) \mbox{ a.s., when }\beta> \max\left(1-\alpha,\, \frac{\mu(\mathbb{R}^+)-1}{\mu(\mathbb{R}^+)}\right),
\end{equation}
where we require $\beta > (\mu(\mathbb{R}^+)-1)/\mu(\mathbb{R}^+)$ to ensure that $L_f(\gamma_+)>1$. Theorem \ref{Thm.Levy1} also yields
\[
\limsup_{t\to\infty}\frac{|X(t)|}{t^\epsilon} = 0 \mbox{ a.s., for each }\epsilon> \max\left(\frac{1}{\alpha},\, \frac{1}{1-\beta}\right).
\] 
In other words, the solution of \eqref{eq.introsfde} is $o(t^\epsilon)$ with probability one for $\epsilon$ sufficiently large (in terms of both the noise and nonlinearity). Define the function $\gamma_-$ by
\[
\gamma_-(t) = (1+t)^\delta, \quad t \geq 0, \quad 0 < \delta \leq \frac{1}{\alpha}.
\]
Note that $\gamma_-$ is positive, increasing and obeys $\int_0^\infty {\gamma_-(t)}^{-\alpha}\,dt = \infty$. Since we aim to apply Theorem \ref{Thm.Levy1} we are only interested in the case $L_f(\gamma_-) \in (1,\infty]$. It is straightforward to show that
\[
L_f(\gamma_-) = 
\begin{cases}
\delta/\mu(\mathbb{R}^+), & \delta = 1/(1-\beta) \leq 1/\alpha,\\
\infty, & 1/(1-\beta) < \delta \leq 1/\alpha.
\end{cases}
\]
Hence Theorem \ref{Thm.Levy1} yields
\begin{equation}\label{example_est_2}
\limsup_{t\to\infty}\frac{|X(t)|}{t^{1/(1-\beta)}} = \infty \mbox{ a.s., when } \frac{\mu(\mathbb{R}^+)-1}{\mu(\mathbb{R}^+)} < \beta \leq 1 - \alpha, \quad \mbox{i.e. }\frac{1}{1-\beta}=\delta \leq \frac{1}{\alpha},
\end{equation}
and
\[
\limsup_{t\to\infty}\frac{|X(t)|}{t^\delta} =\infty \mbox{ a.s. for each $\delta$ such that } \frac{1}{1-\beta} < \delta \leq \frac{1}{\alpha}.
\]
This example highlights a limitation of Theorem \ref{Thm.Levy1} (and it's deterministic counterpart Theorem \ref{Thm.det.gamma.pm}). By comparing \eqref{example_est_1} and \eqref{example_est_2} the reader can see that it is not possible to have both $L_f(\gamma_+) \in (1,\infty)$ and $L_f(\gamma_-) \in (1,\infty)$ simultaneously in this example; indeed this case is difficult to engineer and only possible in limited circumstances (such as when the nonlinearity is regularly varying with unit index). 
\end{example}
\section{Proofs of Results for Deterministic Volterra Equations}\label{sec_deterministic_proofs}
To improve readability of the proofs, we let $\bar{m} := \lim_{t\to\infty}m(t) = \mu(\mathbb{R}^+)$ in the following sections. 

\begin{proof}[Proof of Proposition \ref{prop.linear}] Let $H_1(t) = \int_0^t H(s)\,ds$ for $t \geq 0$ so that \[L = \lim_{t\to\infty}\frac{H(t)}{a\, H_1(t)} = \lim_{t\to\infty} \frac{H_1'(t)}{a\, H_1(t)}.\] 
	
	Case (iii.): If $L = \infty$, then $H(t) \to \infty$ and $H_1(t) \to \infty$ as $t \to \infty$ due to positivity. Furthermore, $H_1'(t)/H(t) \to \infty$ as $t \to \infty$ and asymptotic integration of that limit shows that
	\[
	\lim_{t\to\infty}\frac{1}{t}\log H_1(t)
	= \infty.\]  
	Moreover, since $H_1(t)$ is o$(H(t))$ as $t\to\infty$, $\lim_{t\to\infty}\log H(t)/t = \infty$ as well. Applying the variation of parameters formula to the linear ODE \eqref{eq.linearODE} gives
	\[
	x(t) = x(0)e^{at} + e^{at} \int_0^t e^{-as} h(s)\,ds, \quad t \geq 0.
	\]
	Let $J(t) = e^{at}\int_0^t e^{-as} H(s)\,ds$ and $K(t) = e^{at}\int_0^t e^{-as}H_1(s)\,ds$. Now apply integration by parts to show that 
	\begin{equation}\label{eq.identity1}
		x(t) = x(0)e^{at} + H(t) + ae^{at}\int_0^t e^{-as} H(s)\,ds = x(0)e^{at} + H(t) +aJ(t), \quad t \geq 0.	
	\end{equation}
	Note that $\lim_{t\to\infty}J(t) = \infty$ since $\lim_{t\to\infty}H(t) = \infty$. Another application of integration by parts yields
	\begin{equation}\label{eq.indentity2}
		J(t) = H_1(t) + a K(t).	
	\end{equation}
	$L = \infty$ implies that $\lim_{t\to\infty}H_1(t)/H(t) = 0$ and hence there exists $T(\epsilon)>0$ such that
	\[
	H_1(t) < \epsilon H(t), \quad t \geq T(\epsilon).
	\]
	Thus
	\begin{align}\label{eq.est1_linear}
		K(t) &\leq e^{at} \int_0^T e^{-as} H_1(s)\,ds + e^{at}\int_T^t e^{-as} \epsilon\,H(s)\,ds \leq e^{at}\int_0^T e^{-as} H_1(s)\,ds + \epsilon J(t), \quad t > T(\epsilon).	
	\end{align}
	Divide by $J(t)$ to obtain
	\[
	\frac{K(t)}{J(t)} \leq \epsilon + \frac{\int_0^T e^{-as}H_1(s)\,ds }{\int_0^t e^{-as}H_1(s)\,ds}
	\]
	and since $\lim_{t\to\infty}\log H_1(t)/t = \infty$, it follows that $\lim_{t\to\infty}K(t)/J(t) = 0$. Thus from \eqref{eq.est1_linear} we have that $J(t) \sim H_1(t)$ as $t\to\infty$ and $\lim_{t\to\infty}J(t)/H(t) = 0$. Combining these limits with \eqref{eq.identity1} shows that $x(t) \sim H(t)$ as $t\to\infty$. 
	
	Case (i.): Suppose $L = 0$ so that 
	\[
	\lim_{t\to\infty}\frac{H(t)}{a H_1(t)} = \lim_{t\to\infty}\frac{H_1'(t)}{a H_1(t)} = 0. 
	\]  
	Apply integration by parts to \eqref{eq.identity1} to obtain
	\begin{equation}\label{eq.L=0identity}
		x(t) = x(0)e^{at} + H(t) + a H_1(t) + a^2 e^{at}\int_0^t e^{-as} H_1(s)\,ds, \quad t \geq 0. 	
	\end{equation}
	Since $L = 0$ we have
	\[
	\lim_{t\to\infty}\frac{\log H_1(t)}{t} = \lim_{t\to\infty}\frac{\log H(t)}{t} = 0.
	\]
	It is thus clear from \eqref{eq.L=0identity} that
	\begin{align}\label{eq.final_linear}
		\lim_{t\to\infty}e^{-at}x(t) &= x(0) + \lim_{t\to\infty} a^2 \int_0^t e^{-as}H_1(s)\,ds \nonumber \\
		&= x(0) + \lim_{t\to\infty} a^2 \int_0^t e^{-as} \int_0^s H(u)\,du\,ds=: \xi^* \in (0,\infty),
	\end{align}
	where $t \mapsto e^{-at}H_1(t) \in L^1(0,\infty)$ because $\lim_{t\to\infty}\log H_1(t)/t = 0$.
	
	Case (ii.): If $L \in (0,\infty)$, then $H(t)\to\infty$ as $t\to\infty$ and \[\lim_{t\to\infty}\frac{\log H(t)}{t} = \lim_{t\to\infty}\frac{\log H_1(t)}{t} = L.\] 
	If $L<1$, then 
	\[
	\lim_{t\to\infty}e^{-at}\left( H(t) + a H_1(t) \right) = 0 \mbox{ and } t \mapsto e^{-at} H_1(t) \in L^1(0,\infty). 
	\]
	Hence we can let $t\to\infty$ in \eqref{eq.final_linear} to once more show that $\lim_{t\to\infty}x(t)/e^{at} = \xi^* \in (0,\infty)$. 
	
	If $L>1$, then $H_1(t) \sim H(t)/aL$ as $t\to\infty$ and hence
	\[
	K(t) \sim e^{at}\int_0^t e^{-as} \frac{H(s)}{aL}\,ds =  \frac{J(t)}{aL}
	\]
	as $t\to\infty$. Now since $J(t) = H_1(t) + aK(t)$, it follows that 
	\[
	H_1(t) = J(t) - aK(t) \sim \left(1- \frac{1}{L}\right)J(t) \mbox{ as }t\to\infty.
	\]
	Use the asymptotic relation between $H$ and $H_1$ to show that $\lim_{t\to\infty} aJ(t)/H(t) = 1/(L-1)$. Furthermore, since $L>1$, $\lim_{t\to\infty}e^{at}/H(t) = 0$. Now from \eqref{eq.identity1} we have
	\[
	\frac{x(t)}{H(t)} = \frac{x(0)e^{at}}{H(t)} + 1 + \frac{aJ(t)}{H(t)} \to \frac{L}{L-1} \mbox{ as }t\to\infty,
	\]
	as claimed. 
	
	Finally, if $L = 1$, $H_1(t) \sim H(t)/a$ and $aK(t) \sim J(t)$ as $t\to\infty$. It follows that $\lim_{t\to\infty}H_1(t)/J(t) = 0$ and hence that $\lim_{t\to\infty}J(t)/H_1(t) = \infty$. Thus using a lower estimate from \eqref{eq.identity1} shows that 
	\[
	\lim_{t\to\infty}\frac{x(t)}{H_1(t)} \geq \lim_{t\to\infty}\frac{aJ(t)}{H_1(t)} = \infty
	\] 
	which in turn establishes that $\lim_{t\to\infty}x(t)/H(t) = \infty$. It follows that $\lim_{t\to\infty}x'(t)/x(t) = a$ and hence that $\lim_{t\to\infty}\log x(t)/t = a$. 
\end{proof}
\begin{proof}[Proof of Lemma \ref{phi_props}]
	Suppose that $x \geq a >0$. $\phi(x) - \phi(a) = \int_a^x \phi'(u)du \geq \phi'(x)(x-a)$. Thus 
	\begin{align}
		\limsup_{x\to\infty}\frac{\phi'(x)x}{\phi(x)} = \limsup_{x\to\infty}\frac{\phi'(x)(x-a)}{\phi(x)}\frac{x}{x-a} \leq \limsup_{x\to\infty}\frac{\phi(x)-\phi(a)}{\phi(x)} = 1,
	\end{align}
	establishing the first part of \eqref{props}. To prove the second claim estimate as follows
	\begin{align*}
		\frac{\phi(\Lambda x)}{\phi(x)} &= \frac{\int_a^{\Lambda x} \phi'(u)du + \phi(a)}{\phi(x)} = \frac{\int_a^{x} \phi'(u)du + \int_x^{\Lambda x} \phi'(u)du + \phi(a)}{\phi(x)} = 1 + \frac{\int_x^{\Lambda x} \phi'(u)du}{\phi(x)} \\ &\leq 1 + (\Lambda-1)\frac{\phi'(x)\,x}{\phi(x)}.
	\end{align*}
	Taking the limsup and using the first claim completes the proof.
\end{proof}
\begin{proof}[Proof of Lemma \ref{asym_equiv}]
	By hypothesis, for all $\epsilon>0$ there exists $T(\epsilon)>0$ such that for all $t \geq T(\epsilon)$
	\[
	(1-\epsilon)y(t) < x(t) < (1+\epsilon)y(t).
	\]
	Monotonicity of $\phi$ immediately yields
	\[
	\frac{\phi((1-\epsilon)y(t))}{\phi(y(t))} < \frac{\phi(x(t))}{\phi(y(t))} < \frac{\phi((1+\epsilon)y(t))}{\phi(y(t))}, \quad t \geq T.
	\]
	By Lemma \ref{phi_props}, and the divergence of $y$, there exists $T' > T$ such that $\phi((1+\epsilon)y(t)) < (1+\epsilon)^2 \phi(y(t))$ for all $t \geq T'$. Hence $\limsup_{t\to\infty}\phi(x(t))/\phi(y(t)) \leq 1.$ Reversing the roles of $x$ and $y$ in the argument above we have that $\limsup_{t\to\infty}\phi(y(t))/\phi(x(t)) \leq 1$, or equivalently, $\liminf_{t\to\infty}\phi(x(t))/\phi(y(t)) \geq 1$, completing the proof.
\end{proof}
\begin{proof}[Proof of Proposition \ref{PhiafterH}]
	Define $J(t)= \int_0^t \phi(\gamma(s))ds, \, t \geq 0$. Then, because $\phi$ is increasing and invertible, $J'(t) = \phi(\gamma(t))$ and $\gamma(t) = \phi^{-1}(J'(t))$. We begin by considering the case $L_\phi(\gamma) \in (0,\infty)$, so
	\[
	\lim_{t\to\infty}\frac{\phi^{-1}(J'(t))}{J(t)} = L_\phi(\gamma)\bar{m}.
	\]
	Thus for any $\epsilon \in (0,1)$ there exists $T(\epsilon)>0$ such that
	\[L_\phi(\gamma)\bar{m}(1-\epsilon) < \frac{\phi^{-1}(J'(t))}{J(t)} < L_\phi(\gamma)\bar{m}(1+\epsilon), \quad t \geq T.\] Now since $\phi$ is increasing 
	\begin{subequations}
		\begin{gather}
			\phi\left(L_\phi(\gamma)\bar{m}(1-\epsilon)J(t) \right) < J'(t) < \phi\left(L_\phi(\gamma)\bar{m}(1+\epsilon)J(t) \right), \label{est_1}\\
			L_\phi(\gamma)\bar{m}(1-\epsilon)J(t) < \gamma(t) < L_\phi(\gamma)\bar{m}(1+\epsilon)J(t), \label{est_2}
		\end{gather}
	\end{subequations}
	for all $t \geq T(\epsilon)$. From integrating \eqref{est_1} we obtain
	\[
	\int_T^t \frac{J'(s)ds}{\phi\left(L_\phi(\gamma)\bar{m}(1-\epsilon)J(s)\right)}\geq t-T; \, \int_T^t \frac{J'(s)ds}{\phi\left(L_\phi(\gamma)\bar{m}(1+\epsilon)J(s)\right)}\leq t-T,
	\]
	for all $t \geq T(\epsilon)$. If $a$ is a positive constant then
	\begin{align*}
		\int_T^t \frac{J'(s)ds}{\phi(aJ(s))} &= \int_{aJ(T)}^{aJ(t)}\frac{du}{a\phi(u)}
		=\frac{1}{a}\left\{\Phi(aJ(t)) -  \Phi(aJ(T))\right\}.
	\end{align*}
	With $a = L_\phi(\gamma)\bar{m}(1\pm \epsilon)$ this yields
	\begin{align*}
		&\frac{1}{L_\phi(\gamma)\bar{m}(1-\epsilon)}\left\{\Phi(L_\phi(\gamma)\bar{m}(1-\epsilon)J(t)) -  \Phi(L_\phi(\gamma)\bar{m}(1-\epsilon)J(T))\right\} \geq t-T, \\ &\frac{1}{L_\phi(\gamma)\bar{m}(1+\epsilon)}\left\{\Phi(L_\phi(\gamma)\bar{m}(1+\epsilon)J(t)) -  \Phi(L_\phi(\gamma)\bar{m}(1+\epsilon)J(T))\right\} \leq t-T.
	\end{align*}
	Thus for $t \geq T$
	\begin{align*}
		&\Phi(L_\phi(\gamma)\bar{m}(1-\epsilon)J(t)) \geq L_\phi(\gamma)\bar{m}(1-\epsilon)(t-T)+ \Phi(L_\phi(\gamma)\bar{m}(1-\epsilon)J(T)), \\ &\Phi(L_\phi(\gamma)\bar{m}(1+\epsilon)J(t)) \leq L_\phi(\gamma)\bar{m}(1+\epsilon)(t-T)+ \Phi(L_\phi(\gamma)\bar{m}(1+\epsilon)J(T)).
	\end{align*}
	Applying the monotone function $\Phi$ to \eqref{est_2}, for $t \geq T$, we have
	\begin{align*}
		&\Phi(\gamma(t)) > L_\phi(\gamma)\bar{m}(1-\epsilon)(t-T)+ \Phi(L_\phi(\gamma)\bar{m}(1-\epsilon)J(T)), \\ &\Phi(\gamma(t)) < L_\phi(\gamma)\bar{m}(1+\epsilon)(t-T)+ \Phi(L_\phi(\gamma)\bar{m}(1+\epsilon)J(T)).
	\end{align*}
	Taking limits across the final two sets of inequalities above we obtain
	\[
	\liminf_{t\to\infty}\frac{\Phi(\gamma(t))}{t} \geq \bar{m}L_\phi(\gamma)(1-\epsilon); \,\, \limsup_{t\to\infty}\frac{\Phi(\gamma(t))}{t} \leq L_\phi(\gamma)\bar{m}(1+\epsilon).
	\]
	Letting $\epsilon \to 0^+$ gives the desired result. When $L_\phi(\gamma) = 0$ we will have
	\begin{align*}
		\gamma(t) = \phi^{-1}(J'(t)) < \epsilon J(t), \quad t \geq T_1(\epsilon).
	\end{align*}
	Thus
	$J'(t) < \phi(\epsilon J(t))$ for all $t \geq T_1(\epsilon)$. Integrating we obtain
	\[
	\Phi(\epsilon J(t)) < \epsilon(t-T_1) + \Phi(\epsilon J(T_1)), \quad t \geq T_1.
	\]
	Hence
	\[
	\limsup_{t\to\infty}\frac{\Phi(\gamma(t))}{t} \leq \limsup_{t\to\infty}\frac{\Phi(\epsilon J(t))}{t} \leq \epsilon.
	\]
	It follows immediately that $\lim_{t\to\infty}\Phi(\gamma(t))/t = 0$. Similarly, when $L_\phi(\gamma) = \infty$, we have
	\[
	\gamma(t) = \phi^{-1}(J'(t)) > N J(t), \,\, t \geq T_2(N), \quad N \in \mathbb{R}^+.
	\]
	Integrating by substitution yields
	$
	\Phi(N J(t)) \geq N(t-T_1) - \Phi(N J(T_1)),\,\, t \geq T_1.
	$ 
	Hence
	\[
	\liminf_{t\to\infty}\frac{\Phi(\gamma(t))}{t} \geq \liminf_{t\to\infty}\frac{\Phi(N J(t)}{t} \geq N,
	\]
	and letting $N \to \infty$ completes the proof that $\lim_{t\to\infty}\Phi(\gamma(t))/t = \infty$.
\end{proof}
\begin{proof}[Proof of Theorem \ref{Thm.det.epsilon}]
	With $\Phi$ defined by \eqref{def.Phi}, condition \eqref{fasym} and Lemma \ref{capF_sim_capPhi} imply $F(x) \sim \Phi(x)$ as $x\to\infty$. Therefore, for every $\epsilon \in (0,1)$, there exists $x_1(\epsilon)$ such that
	\[
	\frac{1}{1+\epsilon}\Phi(x) < F(x) < (1+\epsilon)\Phi(x), \,\, x > x_1(\epsilon).
	\]
	Thus $F^{-1}(x)>x_1(\epsilon)$ implies $\tfrac{1}{1+\epsilon}\Phi(F^{-1}(x))< x$ or $x > F(x_1(\epsilon))= x_2(\epsilon)$ implies $F^{-1}(x) < \Phi^{-1}((1+\epsilon)x)$. By hypothesis, for every $\epsilon \in (0,1)$ and $\eta \in (0,1)$, there is $T(\epsilon,\eta)$ such that
	\[
	H(t) < \eta F^{-1}(\bar{m}(1+\epsilon)t), \,\, t \geq T(\epsilon,\eta).
	\]
	Define $T_1(\epsilon)=T(\epsilon,\epsilon)$. For $t \geq T_1(\epsilon)$, $H(t) < \epsilon F^{-1}(\bar{m}(1+\epsilon)t)$. Now let $T_2(\epsilon) = x_2(\epsilon)/(\bar{m}(1+\epsilon))$ and $T_3 = T_1+T_2$. Hence
	\[
	F^{-1}(\bar{m}(1+\epsilon)t)< \Phi^{-1}(\bar{m}(1+\epsilon)^2 t), \,\, t \geq T_3.
	\] 
	But since $t\geq T_3 \geq T_1$, we also have $H(t) < \epsilon \Phi^{-1}(\bar{m}(1+\epsilon)^2 t) < \epsilon \Phi^{-1}(\bar{m}(1+3\epsilon) t)$. Next, because $f(x) \sim \phi(x)$ as $x\to\infty$, there exists $x_3(\epsilon)>0$ such that 
	\[
	\frac{1}{1+4\epsilon}< \frac{f(x)}{\phi(x)}<1+4\epsilon, \,\, x > x_4(\epsilon).
	\]
	Since $\lim_{t\to\infty}x(t) = \infty$, there is $T_4(\epsilon)>0$, so $x(t) > x_3(\epsilon)$ for $t \geq T_4$. If $T^* = T_4 + T_3$, then
	\begin{align}\label{x_est_upper}
		x(t) &= x(0) + H(t) + \int_0^{T^*} \bar{m}(t-s)f(x(s))ds + \int_{T^*}^t \bar{m}(t-s)f(x(s))ds\nonumber\\
		&\leq x(0)+ H(t) + \bar{m} \int_0^{T^*}f(x(s))ds + (1+4\epsilon)\bar{m}\int_{T^*}^t \phi(x(s))ds\nonumber\\
		&\leq x(0) + \epsilon \Phi^{-1}(\bar{m}(1+3\epsilon) t) + x_*(\epsilon) + (1+4\epsilon)\bar{m}\int_{T^*}^t \phi(x(s))ds,\quad t \geq T^*,
	\end{align}
	where $x_*(\epsilon) = \bar{m} \int_0^{T^*}f(x(s))ds$. For $t \geq T^*$, define the function $z_\epsilon$ by
	\[
	z_\epsilon(t) = 1 + x_*(\epsilon) + \epsilon \Phi^{-1}(\bar{m}(1+3\epsilon) t) + (1+4\epsilon)\bar{m}\int_{T^*}^t \phi(z_\epsilon(s))ds.
	\]
	By construction $x(t) < z_\epsilon(t)$ for all $t \geq T^*$. Since $z_\epsilon$ is differentiable we have
	\begin{align*}
		z_\epsilon'(t) &= \epsilon \bar{m}(1+3\epsilon)\phi(\Phi^{-1}(\bar{m}(1+3\epsilon)t))+(1+4\epsilon)\,\bar{m}\, \phi(z_\epsilon(t)), \quad t \geq T^*,\\
		z_\epsilon(T^*) &= 1 + x_*(\epsilon) + \epsilon \Phi^{-1}(\bar{m}(1+3\epsilon) T^*) = z_*(\epsilon).
	\end{align*}
	Define 
	\[
	z_+(t) = \Phi^{-1}(A(\epsilon) +\bar{m}(1+8\epsilon)(t-T^*)), \quad t \geq T^*,
	\]
	where $A(\epsilon) > \Phi(z_*(\epsilon))+ \bar{m}(1+8\epsilon)T^*$. Then $z_+'(t) = \bar{m}(1+8\epsilon)\phi(z_+(t))$ for $t \geq T^*$ or $z_+'(t) = \bar{m}(1+4\epsilon)\phi(z_+(t))+ 4 \bar{m} \epsilon\phi(z_+(t))$. Choosing $\epsilon \in (0,1)$ guarantees that
	\begin{align*}
		4 \bar{m} \epsilon \phi(z_+(t)) &> 4 \bar{m} \epsilon\phi(\Phi^{-1}(\bar{m}(1+7 \epsilon)t))>4 \bar{m} \epsilon\phi(\Phi^{-1}(\bar{m}(1+3\epsilon)t))> \epsilon \bar{m}(1+3 \epsilon)\phi(\Phi^{-1}(\bar{m}(1+3\epsilon)t)).
	\end{align*}
	Hence
	\[
	z_+'(t) > \bar{m}(1+4\epsilon)\phi(z_+(t)) + \epsilon \bar{m} (1+3\epsilon)\phi(\Phi^{-1}(\bar{m}(1+3\epsilon)t)), \quad t \geq T^*,
	\]
	and $z_+(T^*) = \Phi^{-1}(A(\epsilon))> z_*(\epsilon) = z(T^*)$. From the preceding construction it follows that $z_+(t) > z_\epsilon(t) > x(t)$ for all $t \geq T^*$. Hence, from the definition of $z_+$,
	\[
	\Phi(x(t)) < A(\epsilon) + \bar{m}(1+8\epsilon)(t-T^*), \quad t \geq T^*.
	\]
	It follows that $\limsup_{t\to\infty}\Phi(x(t))/t \leq \bar{m}(1+8\epsilon)$ and letting $\epsilon\to 0^+$ shows that
	\[
	\limsup_{t\to\infty}\frac{\Phi(x(t))}{\bar{m}t} \leq 1.
	\]
	The lower bound is proved similarly and we refer the reader to Theorem \ref{Thm.det.zero}. Since $F \sim \Phi$, we will have $\lim_{t\to\infty}F(x(t))/\bar{m}t = 1$, as claimed.
	
	We now establish the second part of \eqref{eq.thm.Finv}, namely that $\lim_{t\to\infty}x(t)/H(t)=\infty$.
	By hypothesis and the first part of \eqref{eq.thm.Finv}, for an arbitrary $\epsilon \in (0,1)$ (chosen so small that $\bar{m}(1-\epsilon)/\epsilon > 1$), there exists $T_0(\epsilon)>0$ such that 
	\[
	F(x(t)) > \bar{m}(1-\epsilon)t, \quad F(H(t)) < \epsilon t, \quad t \geq T_0(\epsilon).
	\]
	Therefore, for $t \geq T_0(\epsilon)$,
	\[
	\frac{x(t)}{H(t)} > \frac{F^{-1}(\bar{m}(1-\epsilon)t)}{F^{-1}(\epsilon t)}.
	\]
	Hence with $K = K(\epsilon) = \bar{m}(1-\epsilon)/\epsilon > 1$, and with $y$ defined by $y'(t) = f(y(t))$ for $t>0$ and $y(0)=1$,
	\begin{align*}
		\liminf_{t\to\infty}\frac{x(t)}{H(t)} \geq \liminf_{t\to\infty}\frac{F^{-1}(\bar{m}(1-\epsilon)t)}{F^{-1}(\epsilon t)} = \liminf_{\tau\to\infty}\frac{F^{-1}(K \tau)}{F^{-1}(\tau)} = \liminf_{\tau\to\infty}\frac{y(K \tau)}{y(\tau)}.
	\end{align*}
	We show momentarily that 
	\begin{equation} \label{eq.yNtyt}
		\liminf_{\tau\to\infty}\frac{y(N \tau)}{y(\tau)}\geq N, \quad\text{ for any } N \geq 1.
	\end{equation}
	Using \eqref{eq.yNtyt} yields
	\[
	\liminf_{t\to\infty}\frac{x(t)}{H(t)} \geq \liminf_{\tau\to\infty}\frac{y(K \tau)}{y(\tau)}\geq K=\frac{\bar{m}(1-\epsilon)}{\epsilon}.
	\]
	Since $\epsilon$ was chosen arbitrarily, letting $\epsilon \to 0$ yields $\liminf_{t\to\infty} x(t)/H(t)=+\infty$, as required.
	
	Now we return to the proof of \eqref{eq.yNtyt}. Clearly, $\lim_{t\to\infty}y(t)=\infty$ and therefore there exists $T_1(\epsilon)>0$ such that $f(y(t)) > (1-\epsilon)\phi(y(t))$ for all $t\geq T_1(\epsilon)$. Let $t \geq T_1(\epsilon)$ and $N>1$, then by using the monotonicity of $\phi$ we obtain
	\begin{align*}
		y(Nt) = y(t) + \int_t^{Nt} f(y(s))\,ds \geq y(t)+\int_t^{Nt} (1-\epsilon)\phi(y(s))\,ds\geq y(t)+(N-1)t (1-\epsilon)\phi(y(t)).
	\end{align*}
	Since $y(t)=F^{-1}(t)$ for $t\geq 0$, we have for $t\geq T_1(\epsilon)$
	\[
	\frac{y(Nt)}{y(t)}\geq 1 + (1-\epsilon)(N-1) \frac{t\,\phi(F^{-1}(t))}{F^{-1}(t)}.
	\]
	Letting $t\to\infty$ yields
	\[
	\liminf_{t\to\infty} \frac{y(Nt)}{y(t)}\geq 1 + (1-\epsilon)(N-1) \liminf_{x\to\infty}\frac{F(x)\phi(x)}{x}
	=1 + (1-\epsilon)(N-1) \liminf_{x\to\infty}\frac{\Phi(x)\phi(x)}{x},
	\]
	since $\Phi(x)\sim F(x)$ as $x\to\infty$. Finally, as $\phi$ is increasing
	$
	\Phi(x)=\int_{1}^x \frac{1}{\phi(u)}\,du \geq x-1/\phi(x),
	$
	so
	\[
	\liminf_{t\to\infty} \frac{y(Nt)}{y(t)}\geq 1 + (1-\epsilon)(N-1).
	\]
	Since $\epsilon\in (0,1)$ was chosen arbitrarily, letting it tend to zero gives the desired bound \eqref{eq.yNtyt}.
\end{proof}
\begin{proof}[Proof of Theorem \ref{Thm.det.zero}]
	Firstly, with $\epsilon\in (0,1)$ arbitrary, rewrite \eqref{eq.xpert} as follows
	\begin{align*}
		x(t) &\leq x(0) + H(t) + \bar{m}\int_0^T f(x(s))ds + \bar{m}\int_T^t f(x(s))ds\\
		&\leq H_\epsilon(t) + (1+\epsilon)\bar{m}\int_T^t \phi(x(s))ds,\quad t \geq T,
	\end{align*}
	where $H_\epsilon(t) = x(0) + H(t) + \bar{m}\int_0^T f(x(s))ds$. Define $I_\epsilon(t)= \int_T^t \phi(x(s))ds$ for $t \geq T$, so
	\begin{align}\label{eq.x.est}
		x(t) \leq H_\epsilon(t) + (1+\epsilon)\bar{m} I_\epsilon(t) , \quad t \geq T.
	\end{align}
	Hence
	\begin{equation} \label{eq.Iprest}
		I_\epsilon '(t) = \phi(x(t)) < \phi\left( H_\epsilon(t) +\bar{m}(1+\epsilon)I_\epsilon(t) \right), \quad t \geq T.
	\end{equation}
	Note that $\lim_{t\to\infty}I_\epsilon(t) = \infty$. We claim
	\begin{equation} \label{eq.IepsHeps0}
		\lim_{t\to\infty}\frac{H_\epsilon(t)}{I_\epsilon(t)} =0.
	\end{equation}
	Suppose first that $\limsup_{t\to\infty} H(t)<\infty$. In this case $\limsup_{t\to\infty} H_\epsilon(t)<\infty$, but $\lim_{t\to\infty}I_\epsilon(t) = \infty$, and \eqref{eq.IepsHeps0} holds.
	
	Suppose next that $\limsup_{t\to\infty} H(t)=+\infty$. Since $f(x)\sim \phi(x)$ as $x\to\infty$, there is $x_1(\epsilon)>0$ such that $f(x)<(1+\epsilon)\phi(x)$ 
	for all $x\geq x_1(\epsilon)$. By the continuity of $f$ and $\phi$ the number $K=K_0(\epsilon)$ given by 
	$
	K_0(\epsilon)=\inf_{x\in (0,x_1(\epsilon))} \phi(x)/f(x)
	$
	is well--defined, and in $(0,\infty)$, even in the case when $f(0)=0$. Therefore, with $K_1(\epsilon)=\min(K_0(\epsilon),1/(1+\epsilon))$, 
	we have $\phi(x)\geq K_1(\epsilon) f(x)$ for all $x>0$. Since $H(t)>0$ for $t>0$, the estimate 
	\[
	\int_T^t \phi(H(s))\,ds \geq K_1(\epsilon)\int_T^t f(H(s))\,ds
	\]
	holds for $t\geq T$. Therefore,
	\begin{equation} \label{eq.HoverintphiH0}
		\frac{H(t)}{\int_T^t \phi(H(s))\,ds}
		\leq \frac{1}{ K_1(\epsilon)}\cdot \frac{H(t)}{\int_0^t f(H(s))\,ds}\cdot\frac{\int_0^t f(H(s))\,ds}{\int_T^t f(H(s))\,ds}, \quad t\geq T.
	\end{equation}
	Since $f$ and $H$ are positive, $t\mapsto \int_0^t f(H(s))\,ds$ tends to some $L \in(0,\infty)$ or infinity as $t\to\infty$. Suppose the former pertains. Then, because $L_f(H)=0$, $H(t)\to 0$ as $t\to\infty$, contradicting the hypothesis that 
	$\limsup_{t\to\infty} H(t)=\infty$. Thus, $\int_0^t f(H(s))\,ds\to\infty$ as $t\to\infty$, and the last quotient on the righthand side of \eqref{eq.HoverintphiH0} is an indeterminate limit as $t\to\infty$. But by l'H\^opital's rule, and because $L_f(H)=0$,
	\[
	\lim_{t\to\infty} \frac{H(t)}{\int_T^t \phi(H(s))\,ds}=0.
	\]
	To complete the proof of \eqref{eq.IepsHeps0} note that positivity of $H$ implies $\phi(x(t)) > \phi(x(0) + H(t))> \phi(H(t))$. Thus $I_\epsilon(t)=\int_T^t \phi(x(s))ds \geq \int_T^t \phi(H(s))ds$. Hence, because $I_\epsilon(t)\to \infty$ as $t\to\infty$,
	\begin{align*}
		\limsup_{t\to\infty}\frac{H_\epsilon(t)}{I_\epsilon(t)} 
		= \limsup_{t\to\infty} \left\{\frac{x(0) + \bar{m}\int_0^T f(x(s))ds}{I_\epsilon(t)} + \frac{H(t)}{I_\epsilon(t)}\right\}
		\leq \limsup_{t\to\infty} \frac{H(t)}{\int_T^t \phi(H(s))ds} =0,
	\end{align*}
	and \eqref{eq.IepsHeps0} holds.
	
	Equation \eqref{eq.IepsHeps0} implies that for every $\eta \in (0,1)$ there is $T'(\eta,\epsilon)>0$ such that $H_\epsilon(t) < \eta I_\epsilon(t)$ for all $t \geq T'(\eta,\epsilon)$. Hence for $t \geq T'(\epsilon,\epsilon)$, $H_\epsilon(t) < \bar{m}\epsilon I_\epsilon(t)$. Then for $t \geq T_2 = T + T'$,
	\begin{align*}
		I_\epsilon'(t) &< \phi(H_\epsilon(t) + \bar{m}(1+\epsilon)I_\epsilon(t)) < \phi(\bar{m}(1+2\epsilon)I_\epsilon(t)).
	\end{align*}
	Integrating we obtain
	\[
	\int_{T_2}^t \frac{I_\epsilon '(s) ds}{\phi(\bar{m}(1+2\epsilon)I_\epsilon(t))} \leq t -T_2, \quad t \geq T_2.
	\]
	Integrating by substitution with $u = \bar{m}(1+2\epsilon)I_\epsilon(s)$
	\[
	\Phi\left( \bar{m}(1+2\epsilon)I_\epsilon(t) \right) - \Phi\left( \bar{m}(1+2\epsilon)I_\epsilon(T_2) \right) \leq \bar{m}(1+2\epsilon)(t-T_2), \quad t \geq T_2.
	\]
	Letting $\Phi_\epsilon=\Phi\left( \bar{m}(1+2\epsilon)I_\epsilon(T_2) \right)$, we have
	\[
	I_\epsilon(t) \leq \frac{1}{\bar{m}(1+2\epsilon)}\Phi^{-1}(\Phi_\epsilon + \bar{m}(1+2\epsilon)(t-T_2)), \quad t \geq T_2.
	\]
	From \eqref{eq.x.est} we have $x(t) \leq H_{\epsilon}(t) + \bar{m}(1+\epsilon)I_\epsilon(t)$ for $t \geq T$ and for $t \geq T'$ we have $H_\epsilon(t) < \bar{m}\epsilon I_\epsilon(t)$. Hence
	\begin{align*}
		x(t) \leq \bar{m}\epsilon I_\epsilon(t) + \bar{m}(1+\epsilon)I_\epsilon(t) = \bar{m}(1+2\epsilon)I_\epsilon(t) \leq \Phi^{-1}(\Phi_\epsilon + \bar{m}(1+2\epsilon)(t-T_2)), \quad t \geq T_2.
	\end{align*}
	Therefore $\Phi(x(t)) < \Phi_\epsilon + \bar{m}(1+2\epsilon)(t-T_2)$ and hence $\limsup_{t\to\infty}\Phi(x(t))/t \leq \bar{m}(1+2\epsilon)$. Letting $\epsilon \to 0^+$ we have $\Phi(x(t))/\bar{m}t \leq 1$ and, since $F(x) \sim \Phi(x)$ as $x\to\infty$ by Lemma \ref{capF_sim_capPhi}, this implies 
	\[
	\limsup_{t\to\infty}\frac{F(x(t))}{\bar{m}t} \leq 1.
	\]
	We now proceed to compute the corresponding lower bound. Since $\lim_{t\to\infty}m(t) = \bar{m} < \infty$, there exists $T_3>0$ such that $m(t) >\bar{m}(1-\epsilon)$, for all $t \geq T_3$, with $\epsilon \in (0,1)$ arbitrary. For $t \geq 2T_3$ 
	\begin{align*}
		x(t) &\geq x(0) + \int_0^{T_3} m(t-s)f(x(s))ds + \int_{T_3}^t m(t-s)f(x(s))ds\\ &\geq x(0) + (1-\epsilon)\int_{T_3}^t m(t-s)\phi(x(s))ds \geq x(0) + (1-\epsilon)^2\bar{m}\int_{T_3}^t \phi(x(s))ds.
	\end{align*}
	Letting $y(t) = x(t+T)$ for $t \geq 2T_3$, it is straightforward to show that 
	\[
	y(t) \geq x(0) + \bar{m}(1-\epsilon)^2\int_0^{t-T_3}\phi(y(u))du, \quad t \geq T_3.
	\]
	Now define the lower comparison solution
	\begin{align*}
		z(t) &= z^* + \bar{m}(1-\epsilon)^2\int_0^{t-T_3}\phi(z(u))du, \quad t \geq T_3,
	\end{align*}
	and $z(t) = z^* = \tfrac{1}{2}\min_{t\in[0,2T_3]}x(t), \,\, t \in [0,T_3].$
	Thus for $t \in [0,T_3]$,\\ $y(t) = x(t+T_3) > z^* = z(t)$ and $z^* < x(0)$. Now suppose that $y(t)>z(t)$ for $t \in [0,\bar{T})$, $\bar{T}>T_3$, but $y(\bar{T}) = z(\bar{T)}$. Then $s \in [0,\bar{T}-T_3]$ implies $\phi(y(s)) > \phi(z(s))$ and
	$
	\int_0^{\bar{T}-T_3}\phi(y(s))ds \geq \int_0^{\bar{T}-T_3}\phi(z(s))ds. 
	$
	Therefore 
	\begin{align*}
		y(\bar{T}) &\geq x(0) + \bar{m}(1-\epsilon)^2 \int_0^{\bar{T}-T_3}\phi(y(s))ds 
		\geq x(0) + \bar{m}(1-\epsilon)^2 \int_0^{\bar{T}-T_3}\phi(z(s))ds \\
		&> z^* + \bar{m}(1-\epsilon)^2 \int_0^{\bar{T}-T_3}\phi(z(s))ds = z(\bar{T}) = y(\bar{T}),
	\end{align*}
	a contradiction. Hence $x(t+T_3) = y(t) > z(t)$ for all $t \geq 0$. For $t \geq T_3$, $z'(t) = \bar{m}(1-\epsilon)^2\phi(z(t-T_3))$ and thus by \cite[Corollary 2]{sublinear2015}, 
	$
	\lim_{t\to\infty}\Phi(z(t))/t = \bar{m}(1-\epsilon)^2,
	$
	under \eqref{fasym}. Hence 
	\[
	\liminf_{t\to\infty}\frac{\Phi(x(t+T_3))}{t} \geq \liminf_{t\to\infty}\frac{\Phi(z(t))}{t} \geq \bar{m}(1-\epsilon)^2.
	\]
	Thus
	\[
	\bar{m}(1-\epsilon)^2 \leq \liminf_{t\to\infty}\frac{\Phi(x(t))}{t-T_3} = \liminf_{t\to\infty}\frac{\Phi(x(t))}{t}.
	\]
	Recall Lemma \ref{capF_sim_capPhi} and let $\epsilon \to 0^+$ to obtain 
	$
	\liminf_{t\to\infty}F(x(t))/\bar{m}t \geq 1,
	$ 
	proving the first limit in \eqref{eq.thm21FxtMxHinf}.
	
	\noindent The proof of the second limit in \eqref{eq.thm21FxtMxHinf} is identical to the proof of the same statement in Theorem \ref{Thm.det.epsilon}.
\end{proof}
\begin{proof}[Proof of Theorem \ref{Thm.det.L.F}]
	The required lower bound, $\liminf_{t\to\infty}F(x(t))/\bar{m}t \geq 1$, can be derived exactly as in Theorem \ref{Thm.det.zero}. For the upper bound, recall the estimate \eqref{eq.Iprest} from the proof of Theorem \ref{Thm.det.zero}:
	\begin{equation*} 
		I_\epsilon '(t) < \phi\left( H_\epsilon(t) +\bar{m}(1+\epsilon)I_\epsilon(t) \right), \quad t \geq T,
	\end{equation*}
	where $I_\epsilon(t) = \int_T^t \phi(x(s))\,ds$ for $t \geq T$ and $H_\epsilon(t) = x(0) + H(t) + \bar{m} \int_0^T f(x(s))\,ds$.
	\begin{remark}
		The stronger hypothesis \eqref{fasym2} can be used to improve the estimate above . We state this improvement here for convenience. Using the mean value theorem, \eqref{fasym2} and the first part of Lemma \ref{phi_props}, estimate as follows: 
		\begin{align}\label{I'main_est}
			I_\epsilon '(t) &\leq \phi(H_\epsilon(t) + \bar{m}(1+\epsilon)I_\epsilon(t))
			= \phi(H_\epsilon(t)) + \phi'(H_\epsilon(t) + \bar{m}(1+\epsilon)I_\epsilon(t)\theta_t )\bar{m}(1+\epsilon)I_\epsilon(t)\nonumber\\
			&\leq \phi(H_\epsilon(t)) + \phi'(H_\epsilon(t))\bar{m}(1+\epsilon)I_\epsilon(t)
			\leq \phi(H_\epsilon(t)) + \frac{\phi(H_\epsilon(t))}{H_\epsilon(t)}\bar{m}(1+\epsilon)^2 I_\epsilon(t),
		\end{align}
		where $\theta_t \in [0,1]$ results from using the Mean Value Theorem. The differential inequality above  is now linear in $I_\epsilon(t)$ and can be solved explicitly; we will return to this estimate frequently.
	\end{remark}
	Next, since $x(t) > H(t)$, $\phi(x(t)) > \phi(H(t))$ and
	\[
	\frac{H_\epsilon(t)}{\bar{m} I_\epsilon(t)} = \frac{H_\epsilon(t)}{H(t)}\frac{H(t)}{\bar{m} \int_T^t \phi(x(s))ds} \leq \frac{H_\epsilon(t)}{H(t)}\frac{H(t)}{\bar{m} \int_0^t \phi(H(s))ds}\frac{\int_0^t \phi(H(s))ds}{\int_T^t \phi(H(s))ds}, \quad t \geq T.
	\]
	Hence
	\[
	\limsup_{t\to\infty}\frac{H_\epsilon(t)}{\bar{m} I_\epsilon(t)} \leq L_\phi(H) \limsup_{t\to\infty}\left\{\frac{H_\epsilon(t)}{H(t)}\,\frac{\int_0^t \phi(H(s))ds}{\int_T^t \phi(H(s))ds}\right\} = L_\phi(H).
	\]
	Thus $H_\epsilon(t) < \bar{m}L_\phi(H)(1+\epsilon)I_\epsilon(t)$ for $t \geq T' > T$. 
	Combine this estimate with \eqref{eq.Iprest} to obtain
	\begin{align*}
		I_\epsilon '(t) \leq \phi(H_\epsilon(t)+ \bar{m}(1+\epsilon)I_\epsilon(t)) \leq \phi((\bar{m}+\bar{m}L_\phi(H))(1+\epsilon)I_\epsilon(t)),\quad t \geq T'.
	\end{align*}
	Integrated the inequality above  reads
	\[
	\int_{T'}^t \frac{I_\epsilon '(s) ds}{\phi((\bar{m}+\bar{m}L_\phi(H))(1+\epsilon)I_\epsilon(s))} \leq t - T', \quad t \geq T'.
	\]
	Make the substitution $u = (\bar{m}+\bar{m}L_\phi(H))(1+\epsilon)I_\epsilon(s)$ to obtain
	\begin{multline*}
		\Phi((\bar{m}+\bar{m}L_\phi(H))(1+\epsilon)I_\epsilon(t)) - \Phi((\bar{m}+\bar{m}L_\phi(H))(1+\epsilon)I_\epsilon(T')) \leq (\bar{m}+\bar{m}L_\phi(H))(1+\epsilon)(t-T').
	\end{multline*}
	Define $\Phi_\epsilon = (\bar{m}+\bar{m}L_\phi(H))(1+\epsilon)I_\epsilon(T')$, so
	\[
	\bar{m}(1+L_\phi(H))(1+\epsilon)I_\epsilon(t) \leq \Phi^{-1}(\Phi_\epsilon + (\bar{m}+\bar{m}L_\phi(H))(1+\epsilon)(t-T')).
	\]
	Now combine equation \eqref{eq.x.est} with the inequality above as follows:
	\begin{multline*}
		x(t) \leq H_\epsilon(t) + \bar{m}(1+\epsilon)I_\epsilon(t)
		<\bar{m}(1+\epsilon)(1+L_\phi(H))I_\epsilon(t) \\
		<\Phi^{-1}(\Phi_\epsilon + \bar{m}(1+L_\phi(H))(1+\epsilon)(t-T')), 
	\end{multline*}
	for all $t \geq T'.$
	Thus 
	\[
	\Phi(x(t)) < \Phi_\epsilon + \bar{m}(1+L_\phi(H))(1+\epsilon)(t-T'), \quad t \geq T',
	\]
	and letting $t\to\infty$ yields
	$
	\limsup_{t\to\infty}\Phi(x(t))/\bar{m}t \leq (1+L_\phi(H))(1+\epsilon).
	$
	Recall Lemma \ref{capF_sim_capPhi} and let $\epsilon \to 0^+$ to obtain
	\[
	\limsup_{t\to\infty}\frac{F(x(t))}{\bar{m}t} \leq 1+L_f(H).
	\]
	Now assume that \eqref{fasym2} holds and show that $\liminf_{t\to\infty}x(t)/H(t) \geq 1 + 1/L_f(H)$. Since $t \mapsto m(t)$ is increasing there exists $T_2(\epsilon)>0$ such that $m(t) > (1-\epsilon)\bar{m}$ for all $t \geq T_2(\epsilon)$. Also, $f(x)>(1-\epsilon)\phi(x)$ for all $x \geq x_1(\epsilon)$ and owing to the divergence of $x(t)$ there exists $T_1(\epsilon)$ such that $x(t) > x_1(\epsilon)$ for all $t \geq T_1(\epsilon)$. Therefore, by positivity of $x(t)$,
	\begin{align*}
		x(t) &> H(t) + \bar{m}(1-\epsilon)^2\int_{T_1}^{t-T_2}\phi(x(s))\,ds, \quad t > T_1 + T_2.
	\end{align*} 
	Then, since $x(t) > H(t)$ for all $t \geq 0$,
	\[
	x(t) > H(t) + \bar{m}(1-\epsilon)^2 \int_{T_1}^{t-T_2}\phi(H(s))ds, \quad t \geq T_1+T_2,
	\]
	and it follows immediately that 
	\begin{align}\label{lower_x_over_H}
		\frac{x(t)}{H(t)} > 1 + \frac{1}{L_f(H)}\,\frac{L_f(H)\,\bar{m}(1-\epsilon)^2 \int_{T_1}^{t-T_2}\phi(H(s))ds}{H(t)}, \quad t \geq T_1 + T_2.
	\end{align}
	By hypothesis $H(t) \sim L_f(H)\,\bar{m}\,\int_0^t \phi(H(s))ds$ as $t\to\infty$ and consequently
	\[
	\max_{t-T_2 \leq s \leq t}H(s) \sim \max_{t-T_2 \leq s \leq t}L_f(H)\,\bar{m}\,\int_0^s \phi(H(u))du = L_f(H)\,\bar{m}\,\int_0^t \phi(H(s))ds.
	\]
	Furthermore, because $\phi$ preserves asymptotic equivalence (see Lemma \ref{asym_equiv} and note that it requires \eqref{fasym2}),
	\[
	\phi\left(\max_{t-T_2 \leq s \leq t}H(s)\right) \sim \phi \left( L_f(H)\,\bar{m}\,\int_0^t \phi(H(s))ds \right) \sim \phi(H(t)) \text{ as } t\to\infty.
	\]
	Hence 
	\[
	\limsup_{t\to\infty}\frac{\int_{t-T_2}^t \phi(H(s))ds }{\phi(H(t))} \leq \limsup_{t\to\infty}\frac{T_2\,\phi\left(\max_{t-T_2 \leq s \leq t}H(s)\right)}{\phi(H(t))} = T_2.
	\]
	Using the facts collected above compute as follows 
	\[
	\limsup_{t\to\infty}\frac{\int_{t-T_2}^t \phi(H(s))}{H(t)} = \limsup_{t\to\infty}\frac{\int_{t-T_2}^t \phi(H(s))ds}{\phi(H(t))}\,\frac{\phi(H(t))}{H(t)} \leq T_2\,\limsup_{t\to\infty}\frac{\phi(H(t))}{H(t)} = 0.
	\]
	Similarly, because $\lim_{t\to\infty}H(t)=\infty$, $\lim_{t\to\infty}\int_0^{T_1} \phi(H(s))ds/H(t)=0$. Thus
	\[
	\lim_{t\to\infty}\frac{L_f(H)\,\bar{m}\,\int_{T_1}^{t-T_2}\phi(H(s))ds }{H(t)} = 1.
	\]
	Returning to \eqref{lower_x_over_H} and using the limit above yields
	\[
	\liminf_{t\to\infty}\frac{x(t)}{H(t)} \geq 1 + \frac{(1-\epsilon)^2}{L_f(H)}.
	\]
	Finally, let $\epsilon \to 0^+$ to give the desired conclusion.
\end{proof}
\begin{proof}[Proof of Theorem \ref{Thm.det.L.H} $(a.)$]
	Case $(a.)$ follows from Theorem \ref{Thm.det.L.F} and by taking $\gamma=H$ in Theorem \ref{Thm.det.gamma}. Similarly, the first limit in \eqref{eq.thm25FxtinfxH1} is obtained by choosing $\gamma=H$ in Theorem \ref{Thm.det.gamma}. Note that $L_f(H) \in (1,\infty)$ and our positivity assumptions imply that $H$ is asymptotically increasing.
\end{proof}
\begin{proof}[Proof of Theorem \ref{Thm.det.L.H} $(b.)$]
	The first limit in \eqref{eq.thm25FxtinfxH1} follows from positivity of $H$, which implies\\ $\liminf_{t\to\infty}x(t)/H(t) \geq 1$ directly from \eqref{eq.xpert}, and setting $\gamma=H$ in case $(b.)$ of Theorem \ref{Thm.det.gamma}. The proof of the second limit in \eqref{eq.thm25FxtinfxH1} is straightforward. By hypothesis and Proposition \ref{PhiafterH},  
	$F(H(t))/t\to\infty$ as $t\to\infty$. Therefore, for every $N>1$ there is $T(N)>0$ such that $H(t)>F^{-1}(Nt)$ for $t\geq T(N)$. But $H$ positive implies $x(t)>H(t)$. Thus $x(t)>F^{-1}(Nt)$, or $F(x(t))/t> N$, for all $t\geq T(N)$. Hence $\liminf_{t\to\infty} F(x(t))/t\geq N$. Letting $N\to\infty$ gives the second part of  \eqref{eq.thm25FxtinfxH1}.
\end{proof}
\begin{proof}[Proof of Theorem \ref{Thm.det.gamma} (a.)]
	The hypotheses \eqref{asym_odd} and \eqref{fasym2} imply that there exists $\phi \in C^1(\mathbb{R}^+;\mathbb{R}^+)$ and $K(\epsilon)>0$ such that 
	\begin{align}\label{f_odd_est}
		|f(x)| < K(\epsilon) + (1+\epsilon)\phi(|x|), \mbox{ for all } x \in \mathbb{R}.
	\end{align}
	Now use equation \eqref{f_odd_est} to derive the following preliminary upper estimate on the size of the solution:
	\begin{align*}
		|x(t)| < |x(0)| + |H(t)| + \bar{m}K(\epsilon)\,t + \bar{m}(1+\epsilon)\int_0^t \phi(|x(s)|)ds, \quad t \geq 0.
	\end{align*}
	By L'H\^{o}pital's rule, $\lim_{x\to\infty}\Phi(x)/x = \lim_{x\to\infty}1/\phi(x) = 0$ and hence $\lim_{t\to\infty}\Phi(\gamma(t))/\gamma(t) = 0$. By Proposition \ref{props}, and since $L_f(\gamma) \in (1,\infty)$ by hypothesis, 
	\begin{align}\label{t_over_gamma}
		\lim_{t\to\infty}\frac{A + Bt}{\gamma(t)} = \lim_{t\to\infty}\frac{A+Bt}{\Phi(\gamma(t))}\frac{\Phi(\gamma(t))}{\gamma(t)} = 0,
	\end{align}
	for any nonnegative constants $A$ and $B$. Thus there exists $T(\epsilon)>0$ such that for all $t \geq T(\epsilon)$ we have $|x(0)|+\bar{m}\,K(\epsilon)\,t < \epsilon\,\gamma(t)$. By \eqref{gamma_H}, and the previous estimate, there exists $T_2(\epsilon)>T(\epsilon)$ such that for all $t \geq T_2(\epsilon)$, $|x(0)|+\bar{m}\,K(\epsilon)\,t+|H(t)| < (1+\epsilon)\gamma(t)$. Combining this with our initial estimate we obtain
	\begin{align*}
		|x(t)| &< (1+\epsilon)\gamma(t) + \bar{m}(1+\epsilon)\int_0^t \phi(|x(s)|)ds, \quad t \geq T_2(\epsilon).
	\end{align*}
	To ensure our comparison solution majorizes the true solution take $x^*= \max_{0\leq s \leq T_2}\phi(|x(s)|)$, so $\int_0^{T_2}\phi(|x(s)|)ds \leq T_2\,x^*$. Hence
	\begin{align*}
		|x(t)| &< T_2\,x^* + (1+\epsilon)\gamma(t) + \bar{m}(1+\epsilon)\int_{T_2}^t \phi(|x(s)|)ds, \quad t \geq T_2.
	\end{align*}
	Define the upper comparison solution, $x_+$, as follows:
	\begin{align}\label{x_gamma_initial}
		x_+(t) &= 1+T_2\,x^* + (1+\epsilon)\gamma(t) + \bar{m}(1+\epsilon)\int_{T_2}^t \phi(x_+(s))\,ds = \gamma_\epsilon(t) + \bar{m}(1+\epsilon)I_\epsilon(t), \quad t \geq T_2,
	\end{align}
	where $\gamma_\epsilon(t) =1 + T_2\,x^* + (1+\epsilon)\gamma(t)$ and $I_\epsilon(t) = \int_{T_2}^t \phi(x_+(s))\,ds$. By construction, $|x(t)| < x_+(t)$ for all $t\geq T_2$ (this follows immediately via a ``time of the first breakdown'' argument). Applying the same estimation procedures as in Theorems \ref{Thm.det.zero} and \ref{Thm.det.L.F} to $x_+$, and in particular to the quantity $I_\epsilon(t)$, we obtain an estimate analogous to \eqref{I'main_est}:
	\begin{align}\label{I'est_key_1}
		I_\epsilon '(t) < \phi(\gamma_\epsilon(t)) + \tilde{a}_\epsilon(t) I_\epsilon(t), \quad t \geq T_3(\epsilon),
	\end{align}
	where $\tilde{a}_\epsilon(t) = \bar{m}(1+\epsilon)^2 \phi(\gamma_\epsilon(t))/\gamma_\epsilon(t)$. Note once more that the hypothesis \eqref{fasym2} is needed to obtain the differential inequality \eqref{I'est_key_1}.  Before proceeding further with the line of argument from Theorem \ref{Thm.det.L.F} we need to refine the estimate above. $L_f(\gamma) \in (0,\infty)$ implies that $\lim_{t\to\infty}\gamma(t)=\infty$ and hence, by Lemma \ref{phi_props}, $\limsup_{t\to\infty}\phi(\gamma_\epsilon(t))/\phi(\gamma(t)) \leq (1+\epsilon)$. Therefore there exists a $T_4(\epsilon) > T_3(\epsilon)$ such that for all $t \geq T_4$ we have $\phi(\gamma_\epsilon(t)) < (1+\epsilon)^2\phi(\gamma(t))$. Hence
	\[
	I_\epsilon '(t) < (1+\epsilon)^2\phi(\gamma(t)) + \bar{m} (1+\epsilon)^4\frac{\phi(\gamma(t))}{\gamma_\epsilon(t)}I_\epsilon(t), \quad t \geq T_4.
	\] 
	$\gamma_\epsilon(t) \sim (1+\epsilon) \gamma(t)$ as $t\to\infty$ implies that there exists $T_5(\epsilon)>T_4(\epsilon)$ such that $\gamma_\epsilon(t) > (1-\epsilon)(1+\epsilon) \gamma(t)$ for all $t \geq T_5$. Taking reciprocals of the previous inequality and apply it to the previous estimate of $I_\epsilon '(t)$ to obtain
	\[
	I_\epsilon '(t) < (1+\epsilon)^2\phi(\gamma(t)) + \bar{m} (1+\epsilon)^3\frac{\phi(\gamma(t))}{(1-\epsilon)\gamma(t)}I_\epsilon(t), \quad t \geq T_5.
	\]
	Now let
	\[
	\alpha_\epsilon = (1+\epsilon)^2, \quad a_\epsilon(t) = \bar{m} (1+\epsilon)^3\frac{\phi(\gamma(t))}{(1-\epsilon)\gamma(t)},
	\]
	to obtain the consolidated estimate
	\begin{align}\label{I'_gamma_est}
		I_\epsilon '(t) \leq \alpha_\epsilon\, \phi(\gamma(t)) + a_\epsilon(t)\, I_\epsilon(t), \quad t\geq T_5.
	\end{align}
	Let $T'>T_5$ and solve the differential inequality above as follows
	\begin{align*}
		\frac{d}{dt}\left( I_\epsilon(t)e^{ -\int_{T'}^t a_\epsilon(s)ds } \right) &= I_\epsilon '(t) e^{-\int_{T'}^t a_\epsilon(s)ds} - a_\epsilon(t) I_\epsilon(t)e^{-\int_{T'}^t a_\epsilon(s)ds} \\
		&= e^{-\int_{T'}^t a_\epsilon(s)ds}\left\{ I_\epsilon '(t) - a_\epsilon(t) I_\epsilon(t) \right\}\leq \alpha_\epsilon\,\phi(\gamma(t))e^{-\int_{T'}^t a_\epsilon(s)ds}, \quad t \geq T'.
	\end{align*}
	Integration yields
	\[
	I_\epsilon(t) e^{-\int_{T'}^t a_\epsilon(s)ds} \leq I_\epsilon(T') + \alpha_\epsilon\int_{T'}^t \phi(\gamma(s))e^{-\int_{T'}^s a_\epsilon(u)du}ds, \quad t \geq T'.
	\]
	Hence
	\begin{align}\label{var_of_const}
		\frac{I_\epsilon(t)}{\int_{T'}^t \phi(\gamma(s))ds} \leq \frac{I_\epsilon(T')}{\int_{T'}^t \phi(\gamma(s))ds\,e^{-\int_{T'}^t a_\epsilon(s)ds}} + \frac{\alpha_\epsilon\int_{T'}^t \phi(\gamma(s))e^{-\int_{T'}^s a_\epsilon(u)du}ds}{\int_{T'}^t \phi(\gamma(s))ds\,e^{-\int_{T'}^t a_\epsilon(s)ds}}, \quad t \geq T'.
	\end{align}
	In the analysis which is required to show that the second term on the right-hand side of \eqref{var_of_const} is bounded it emerges that the first term on the right-hand side is also bounded so we immediately focus on the second term. Define 
	\begin{align*}
		C_\epsilon(t) = \alpha_\epsilon\int_{T'}^t \phi(\gamma(s))e^{-\int_{T'}^s a_\epsilon(u)du}\,ds, \quad
		B_\epsilon(t) = \int_{T'}^t \phi(\gamma(s))\,ds\, e^{-\int_{T'}^t a_\epsilon(s)ds},
	\end{align*}
	and restate \eqref{var_of_const} as
	\[
	\frac{I_\epsilon(t)}{\int_{T'}^t \phi(\gamma(s))ds} \leq \frac{I_\epsilon(T')}{B_\epsilon(t)} + \frac{C_\epsilon(t)}{B_\epsilon(t)}, \quad t \geq T'.
	\]
	By inspection $C_\epsilon '(t) > 0$, so either $\lim_{t\to\infty}C_\epsilon(t)=\infty$ or $\lim_{t\to\infty}C_\epsilon(t) = C(\epsilon) \in (0,\infty)$. Differentiating $B_\epsilon$ we obtain
	\begin{align*}
		B_\epsilon '(t) &= \phi(\gamma(t))e^{ -\int_{T'}^t a_\epsilon(s)ds } - a_\epsilon(t)e^{ -\int_{T'}^t a_\epsilon(s)ds } \int_{T'}^t \phi(\gamma(s))\,ds\nonumber\\
		&=e^{ -\int_{T'}^t a_\epsilon(s)ds } \left\{ \phi(\gamma(t)) - a_\epsilon(t)\int_{T'}^t \phi(\gamma(s))\,ds \right\} \nonumber\\
		&= C_\epsilon '(t) \left\{\frac{1}{\alpha_\epsilon} - a_\epsilon(t)\frac{\int_{T'}^t \phi(\gamma(s))\,ds }{\alpha_\epsilon\,\phi(\gamma(t))} \right\} = C_\epsilon '(t) \left\{ \frac{1}{\alpha_\epsilon} - \bar{m}\frac{(1+\epsilon)^4}{(1-\epsilon)} \frac{\int_{T'}^t \phi(\gamma(s))\,ds }{\alpha_\epsilon\,\gamma(t)} \right\}.
	\end{align*}
	Hence
	\begin{align}\label{B'overC'}
		\frac{B_\epsilon '(t)}{C_\epsilon '(t)} = \frac{1}{\alpha_\epsilon} - \frac{(1+\epsilon)^3}{(1-\epsilon)} \frac{\bar{m}\int_T^t \phi(\gamma(s))ds}{\alpha_\epsilon\,\gamma(t)}, \quad  t\geq T'.
	\end{align}
	Therefore, for $\epsilon$ sufficiently small,
	\begin{align}\label{c'_b'}
		\lim_{t\to\infty}\frac{B_\epsilon '(t)}{C_\epsilon '(t)} = \frac{1}{\alpha_\epsilon}- \frac{(1+\epsilon)^3}{(1-\epsilon)\,\alpha_\epsilon\, L_\phi(\gamma)} > 0.
	\end{align}
	\begin{remark}\label{remark_L>1}
		Note that the hypothesis $L_\phi(\gamma)>1$ implies that $B_\epsilon(t)$ is eventually increasing and hence has a limit $B(\epsilon) \in (0,\infty]$ at infinity. If $\lim_{t\to\infty}C_\epsilon(t) = \infty$ and $L_\phi(\gamma) \in (0,1]$, $B_\epsilon(t)$ is eventually decreasing and $\lim_{t\to\infty}B_\epsilon(t) \in [0,\infty)$. In this case $\lim_{t\to\infty}B_\epsilon(t) = 0$ for all $\epsilon \in (0,1)$ and we will be unable to obtain the required estimates to continue the proof.
	\end{remark}
	\noindent From \eqref{c'_b'}, by asymptotic integration, the convergence and divergence of $B_\epsilon$ and $C_\epsilon$ are equivalent. Hence
	\[\lim_{t\to\infty}\frac{C_\epsilon(t)}{B_\epsilon(t)}=
	\begin{cases}
		\left(1/\alpha_\epsilon- (1+\epsilon)^3/(1-\epsilon)\,\alpha_\epsilon\, L_\phi(\gamma)\right)^{-1}, \quad &\lim_{t\to\infty}C_\epsilon(t)=\infty,\\
		C_\epsilon/B_\epsilon, &\lim_{t\to\infty}C_\epsilon(t)=C(\epsilon).
	\end{cases}
	\]
	In both cases
	\[
	\limsup_{t\to\infty}\frac{I_\epsilon(t)}{\int_{T'}^t \phi(\gamma(s))ds} = K(\epsilon) < \infty.
	\]
	Therefore there exists $\bar{T} > T'$ such that $I_\epsilon(t) < K(\epsilon)(1+\epsilon)\int_{T'}^t \phi(\gamma(s))ds$ for all $t \geq \bar{T}$. Thus, recalling  \eqref{x_gamma_initial},
	\begin{align*}
		x_+(t) &= \gamma_\epsilon(t) + \bar{m}(1+\epsilon)I_\epsilon(t)
		\leq (1+2\epsilon)\gamma(t) + \bar{m}(1+\epsilon)^2 K(\epsilon)\int_{T'}^t \phi(\gamma(s))\,ds, \quad t \geq \bar{T}. 
	\end{align*}
	Hence
	\begin{align*}
		\limsup_{t\to\infty}\frac{x_+(t)}{\gamma(t)} \leq 1+2\epsilon + \bar{m}(1+\epsilon)^2 K(\epsilon) \limsup_{t\to\infty}\frac{\int_{T'}^t \phi(\gamma(s))ds}{\gamma(t)}
		= 1+2\epsilon + \frac{(1+\epsilon)^2 K(\epsilon)}{L_\phi(\gamma)} < \infty.
	\end{align*}
	Therefore, since $|x(t)| < x_+(t)$ for all $t \geq T_2$, $\limsup_{t\to\infty}|x(t)|/\gamma(t) < \infty$. Now let
	\begin{align}\label{limsup_finite_+}
		\limsup_{t\to\infty}\frac{|x(t)|}{\gamma(t)} =\lambda \in [0,\infty),
	\end{align}
	One can compute a definite upper bound on $\lambda$ in terms of the problem parameters as follows. Define $J(t) = \int_0^t m(t-s) f(x(s))ds$ and estimate as above
	\begin{align}\label{J_est}
		|J(t)| &\leq \bar{m}\, \int_0^t K(\epsilon)+ (1+\epsilon)\phi(|x(s)|)ds\nonumber\\
		&\leq \bar{m}\,K(\epsilon)\,t + \bar{m}\,T_2\,(1+\epsilon)\sup_{s \in [0,T_2]}\phi(|x(s)|)+ \bar{m}(1+\epsilon)\int_{T_2}^t \phi(|x(s)|)ds, \quad t \geq T_2.
	\end{align}
	Using \eqref{limsup_finite_+} there exists a $\bar{T}(\epsilon)>T_2$ such that
	\begin{align*}
		\limsup_{t\to\infty}\frac{|J(t)|}{\gamma(t)} \leq \bar{m}(1+\epsilon)\limsup_{t\to\infty}\frac{\int_{\bar{T}}^t \phi((\lambda+\epsilon) \gamma(s))ds}{\gamma(t)} \leq \frac{\max(1,\,\lambda+\epsilon)}{L_\phi(\gamma)}.
	\end{align*}
	Return to \eqref{eq.xpert}, take absolute values and apply the estimates above as follows
	\begin{align}\label{limsup_calc}
		\lambda = \limsup_{t\to\infty}\frac{|x(t)|}{\gamma(t)} &\leq \limsup_{t\to\infty}\frac{|x(0)|}{\gamma(t)} + \limsup_{t\to\infty}\frac{|H(t)|}{\gamma(t)} + \limsup_{t\to\infty}\frac{|J(t)|}{\gamma(t)}
		\leq 1 + \frac{\max(1,\,\lambda)}{L_f(\gamma)}.
	\end{align}
	Solving the inequalities above yields $\lambda \leq \max\left((1+L_f(\gamma))/L_f(\gamma),\, L_f(\gamma)/(L_f(\gamma)-1)\right)$. In fact the second quantity is always larger so $\limsup_{t\to\infty}|x(t)|/\gamma(t) \leq L_f(\gamma)/(L_f(\gamma)-1)$.
\end{proof}
\begin{proof}[Proof of Theorem \ref{Thm.det.gamma} (b.)]
	Follow the argument of Theorem \ref{Thm.det.gamma} $(a.)$ exactly to equation \eqref{B'overC'}, which we recall below. 
	\begin{align*}
		\frac{B_\epsilon '(t)}{C_\epsilon '(t)} = \frac{1}{\alpha_\epsilon} - \frac{(1+\epsilon)^3}{(1-\epsilon)} \frac{\bar{m}\int_T^t \phi(\gamma(s))ds}{\alpha_\epsilon\,\gamma(t)}, \quad  t\geq T'.
	\end{align*}
	Now $L_f(\gamma)=\infty$ implies 
	$\lim_{t\to\infty}B_\epsilon '(t)/C_\epsilon '(t) = 1/\alpha_\epsilon$. Thus $0 < C_\epsilon '(t) \sim \alpha_\epsilon B_\epsilon '(t)$ as $t\to\infty$. Recall equation \eqref{var_of_const}
	\[
	\frac{I_\epsilon(t)}{\int_{T'}^t \phi(\gamma(s))ds} \leq \frac{I_\epsilon(T')}{B_\epsilon(t)} + \frac{C_\epsilon(t)}{B_\epsilon(t)}, \quad t\geq T'.
	\]
	If $\lim_{t\to\infty}C_\epsilon(t) = \infty$, then $\lim_{t\to\infty}B_\epsilon(t) = \infty$ and $C_\epsilon(t) \sim \alpha_\epsilon B_\epsilon(t)$ as $t\to\infty$. Thus, when $C_\epsilon(t) \to \infty$ as $t\to\infty$,
	\[
	\limsup_{t\to\infty}\frac{I_\epsilon(t)}{\int_{T'}^t \phi(\gamma(s))ds} \leq \alpha_\epsilon.
	\]
	Alternatively, if $\lim_{t\to\infty}C_\epsilon(t) = C(\epsilon)$, $\lim_{t\to\infty}B_\epsilon(t) = B(\epsilon) \in (0,\infty)$, then
	\[
	\limsup_{t\to\infty}\frac{I_\epsilon(t)}{\int_{T'}^t \phi(\gamma(s))ds} \leq \frac{I_\epsilon(T')+ C(\epsilon)}{B(\epsilon)}.
	\]
	In both cases
	\[
	\limsup_{t\to\infty}\frac{I_\epsilon(t)}{\int_{T'}^t \phi(\gamma(s))ds} \leq K(\epsilon) < \infty.
	\]
	Once more we conclude that $\limsup_{t\to\infty}x_+(t)/\gamma(t) < \infty$ and hence that $\limsup_{t\to\infty}|x(t)|/\gamma(t) < \infty$. By an argument exactly analogous to that which completes the proof of Theorem \ref{Thm.det.gamma} case $(a.)$ we can show that $\lim_{t\to\infty}|J(t)|/\gamma(t) = 0$. Now write 
	\begin{align}\label{final_eqn}
		\frac{x(t)}{\gamma(t)} = \frac{x(0)}{\gamma(t)} + \frac{J(t)}{\gamma(t)} + \frac{H(t)}{\gamma(t)}, \quad t \geq 0.
	\end{align}
	Because $\limsup_{t\to\infty}|H(t)|/\gamma(t) = 1$, $\limsup_{t\to\infty}H(t)/\gamma(t) = 1$ or $\liminf_{t\to\infty}H(t)/\gamma(t) = -1$. Then, since $\lim_{t\to\infty}J(t)/\gamma(t)=0$, taking the limsup and liminf across \eqref{final_eqn} gives $\limsup_{t\to\infty}x(t)/\gamma(t) = 1$ or $\liminf_{t\to\infty}x(t)/\gamma(t) = -1$. In both cases $\limsup_{t\to\infty}|x(t)|/\gamma(t) = 1$. Noting that $J(t)/\gamma(t) \sim (x(t)-H(t))/\gamma(t)$ as $t\to\infty$ yields the second part of the conclusion. 
\end{proof}
\begin{proof}[Proof of Theorem \ref{Thm.det.gamma.pm} $(a.)$]
	The argument of Theorem \ref{Thm.det.gamma} $(a.)$ yields $\limsup_{t\to\infty}|x(t)|/\gamma_+(t) < \infty$. Let $\lambda_+ = \limsup_{t\to\infty}|x(t)|/\gamma_+(t) \in [0,\infty)$ and estimate as before to obtain $\limsup_{t\to\infty}|J(t)|/\gamma_+(t) \leq \max(1,\,\lambda_+)/L_f(\gamma_+)$. Calculating as in \eqref{limsup_calc} then yields $\lambda \leq \max(1,\,\lambda)/L_f(\gamma_+)$. Upon inspection we find that in all cases $\limsup_{t\to\infty}|x(t)|/\gamma(t) \in [0,1/L_f(\gamma_+)]$.
	
	For the second part of the claim suppose to the contrary that $\limsup_{t\to\infty}|x(t)|/\gamma_-(t) = \lambda_- < \infty$. Now argue, as in Theorem \ref{Thm.det.gamma}, that $\limsup_{t\to\infty}|J(t)|/\gamma_-(t) < \max(1,\lambda_-)/L_\phi(\gamma_-)$, where $J(t) = \int_0^t \bar{m}(t-s) f(x(s))\,ds$. However, by rearranging \eqref{eq.xpert} and taking absolute values
	\[
	|H(t)| \leq |x(0)| + |x(t)| + |J(t)|, \quad t \geq 0.
	\]
	Dividing across by $\gamma_-$ and taking the limsup immediately yields $\limsup_{t\to\infty}|H(t)|/\gamma_-(t) < \infty$, in contradiction to \eqref{gamma_pm}. Hence $\lambda_- = \infty$, as claimed.
\end{proof}
\begin{proof}[Proof of Theorem \ref{Thm.det.gamma.pm} $(b.)$]
	As with case $(a.)$, the proof is a consequence of Theorem \ref{Thm.det.gamma} and the stronger conclusion, $\lim_{t\to\infty}|x(t)|\gamma_+(t) = 0$, holds because in \eqref{limsup_calc} we have $\limsup_{t\to\infty}|H(t)|/\gamma_+(t)=0$ and $\limsup_{t\to\infty}|J(t)|/\gamma_+(t)=0$. The proof that $\limsup_{t\to\infty}|x(t)|/\gamma_-(t) = \infty$ is essentially unchanged.
\end{proof}
\section{Proofs of Results for Stochastic Volterra Equations}\label{sec_stoch_proofs}
\begin{proof}[Proof of Theorem \ref{Thm.stoch.epsilon}]
	The proof of this result follows directly from the argument used in the proof of Theorem \ref{Thm.stoch.zero} and the law of the iterated logarithm for continuous local martingales.
\end{proof}
\begin{proof}[Proof of Theorem \ref{Thm.stoch.zero}]
	We start by proving part (a), which covers the case when $\sigma\not \in L^2(0,\infty)$. Let $\epsilon, \, \eta \in (0,1)$ be arbitrary, rewrite \eqref{eq.introsfde} in integral form and estimate as follows
	\begin{align}\label{eq.abs_int_est}
		|X(t)| \leq |X(0)| + \int_0^t m(t-s) |f(X(s))|ds + \left|\int_0^t \sigma(s)dB(s) \right|, \quad t \geq 0.
	\end{align}
	Denote by $\Omega_1$ the a.s. event on which $t \mapsto X(t)(\omega)$ is continuous. We now recall the law of the iterated logarithm for continuous local martingales (see Revuz and Yor \cite[Ch. V, Ex. 1.15]{revuzandyor}) which states that if $N = \{N_t,\,t\geq 0\}$ is a continuous local martingale with $\langle N, N \rangle_\infty = \infty$, then
	\[
	\limsup_{t\to\infty}\frac{N_t}{\sqrt{ 2\langle N, N \rangle_t \log\log \langle N, N \rangle_t  }} = 1 \,\, a.s.,
	\]
	where $\langle N, N \rangle = \{\langle N, N \rangle_t,\, t \geq 0\}$ denotes the quadratic variation process of $N$. In our case 
	\[
	\left\langle \int_0^{\cdot}\sigma(s)dB(s), \int_0^{\cdot}\sigma(s)dB(s) \right\rangle_t = \int_0^t \sigma^2(s)ds
	\]
	and thus $\sigma \notin L^2(0,\infty)$ implies $\limsup_{t\to\infty}\left|\int_0^t \sigma(s)dB(s)\right|/ \Sigma(t) = 1$ a.s.
	
	Let $\eta>0$ be arbitrary. By hypothesis there exists $\phi \in C^1$ such that 
	\begin{equation} \label{eq.fglobalesteta}
		|f(x)| \leq K(\eta) + (1+\eta)\phi(|x|), \, x \in \mathbb{R}.
	\end{equation}
	Define $H_\eta(t) = \bar{m} K(\eta) t + (1+\eta)\Sigma(t)$ for $t \geq 0$. Note that $L_f(\Sigma)=0$ and Proposition \ref{PhiafterH} imply $
	\lim_{t\to\infty}\Phi(\Sigma(t))/t=0.
	$ 
	Therefore, for every $\epsilon\in (0,1)$ there exists $T_2(\epsilon)>0$ such that  
	\begin{align}\label{Sigma_est_2}
		\Sigma(t) < \Phi^{-1}(\epsilon t), \,\, t \geq T_2(\epsilon).
	\end{align}
	Similarly, by L'H\^{o}pital's rule,
	\[
	\lim_{t\to\infty}\frac{ \bar{m} K(\eta) t}{\int_0^t \phi( \bar{m} K(\eta) s)ds} = \lim_{t\to\infty}\frac{\bar{m} K(\eta)}{\phi( \bar{m} K(\eta) t)} = 0.
	\]
	Thus, again appealing to L'H\^{o}pital's rule, $\lim_{t\to\infty}\Phi(\bar{m} K(\eta) t)/t = 0$ and moreover, for any $\eta \in (0,1)$,
	$\lim_{t\to\infty}\Phi\left(\bar{m} K(\eta) t / \eta\right)/t=0$. Hence for every $\epsilon\in (0,1)$ there exists $T_3(\epsilon, \eta)$ such that 
	\begin{align}\label{lin_est}
		\bar{m} K(\eta) t < \eta\, \Phi^{-1}(\epsilon t), \quad t \geq T_3(\epsilon, \eta).
	\end{align} 
	Combining \eqref{Sigma_est_2} and \eqref{lin_est} yields
	\[
	H_\eta(t) =  \bar{m} K(\eta) t + (1+\eta)\Sigma(t) < (1+2\eta) \Phi^{-1}(\epsilon t), \quad t \geq T_4(\epsilon,\eta)= T_2 +T_3.
	\]
	Rearrange this inequality, let $t\to\infty$, and then let $\epsilon \to 0^+$ to obtain 
	$\lim_{t\to\infty}\Phi(H_\eta(t)/(1+2\eta))/\bar{m}t = 0$. Thus, by proceeding as above, for every $\epsilon\in (0,1)$ there is 
	$T_4'(\epsilon,\eta)>0$ such that 
	\begin{align}\label{H_eta_est}
		H_\eta(t) < (1+2\eta) \Phi^{-1}(\epsilon \bar{m} t), \quad t \geq T_4'(\epsilon,\eta).
	\end{align}
	Since $\Phi$ is concave, $\Phi^{-1}$ is convex and 
	$
	\Phi^{-1}(\epsilon \bar{m} t) \leq \epsilon \Phi^{-1}(\bar{m}t) + (1-\epsilon)\Phi^{-1}(0).
	$
	Therefore,
	\[
	\limsup_{t\to\infty}\Phi^{-1}(\epsilon \bar{m} t) / \Phi^{-1}(\bar{m}t) \leq \epsilon.
	\]
	Take limits in \eqref{H_eta_est} to give 
	\[
	\limsup_{t\to\infty} 
	\frac{H_\eta(t)}{\Phi^{-1}(\bar{m} t)} \leq  (1+2\eta) \epsilon,
	\]
	and then let $\epsilon\to 0$ to yield $\lim_{t\to\infty} H_\eta(t)/\Phi^{-1}(\bar{m} t) =0$. Therefore, for every $\epsilon\in (0,1)$ 
	there exists $T_5'(\epsilon,\eta)>0$ such that 
	$
	H_\eta(t) <\epsilon\, \Phi^{-1}(\bar{m} t)$ for $t\geq T_5'(\epsilon,\eta).
	$
	Now, let $T_5(\eta)=T_5'(\eta,\eta)$, so 
	\begin{align}\label{H_eta_est_2}
		H_\eta(t) < \eta\, \Phi^{-1}(\bar{m}t), \quad t \geq T_5(\eta).
	\end{align}
	On the other hand, because $\limsup_{t\to\infty}\left|\int_0^t \sigma(s)dB(s)\right|/ \Sigma(t) = 1$ a.s., 
	there exists an almost sure event $\Omega_2$ such that for all $\omega \in \Omega_2$
	\[
	\left|\int_0^t \sigma(s)dB(s)(\omega)\right| \leq (1+\eta)\Sigma(t), \quad t \geq T_1(\eta,\omega).
	\]
	Now let $T(\eta,\omega)=\max(T_1(\eta,\omega),T_5(\eta))$. 
	Thus for all $\omega \in \Omega^\ast= \Omega_1 \cap \Omega_2$ and $ t \geq T(\eta,\omega)$,
	\begin{align*}
		|X(t)| \leq |X(0)| + \int_0^t m(t-s) |f(X(s))|ds + (1+\eta)\Sigma(t).
	\end{align*} 
	Using the estimate \eqref{eq.fglobalesteta} on $f$ and the finiteness of $\lim_{t\to\infty}\bar{m}(t)$ we have
	\begin{align}\label{abs_est_1}
		|X(t)| &\leq |X(0)| + \bar{m} K(\eta) t + \bar{m}(1+\eta) \int_0^t \phi(|X(s)|)ds + (1+\eta)\Sigma(t) \nonumber\\
		&\leq X^*_0 + H_\eta(t) + \bar{m}(1+\eta) \int_T^t \phi(|X(s)|)ds, \quad t \geq T(\eta,\omega), \quad \omega \in \Omega^\ast,
	\end{align} 
	where $X(0)^* = |X(0)| + \bar{m} T \sup_{s \in [0,T]}\phi(|X(s)|)$. 
	
	Now since $t\geq T(\eta,\omega)\geq T_5(\eta)$, we have from \eqref{H_eta_est_2} that for all $\omega\in \Omega^\ast$ 
	\begin{align}\label{abs_est_2}
		|X(t)| \leq X(0)^\ast + \eta \,\Phi^{-1}(\bar{m}t) + \bar{m}(1+\eta) \int_T^t \phi(|X(s)|)ds, \quad t\geq T(\eta,\omega).
	\end{align}
	At this point we note that we are in the same position as in the proof of Theorem \ref{Thm.det.epsilon} at equation \eqref{x_est_upper}. From here a calculation exactly analogous to that which completes the proof of Theorem \ref{Thm.det.epsilon} will yield 
	\[
	\limsup_{t\to\infty}\frac{F(|X(t)|)}{\bar{m}t} \leq 1 \,\, a.s.
	\]
	
	To prove part (b), let $\epsilon, \, \eta \in (0,1)$ be arbitrary and rewrite \eqref{eq.introsfde} in integral form as before and take absolute values to obtain
	\begin{align*}
		|X(t)| \leq |X(0)| + \int_0^t m(t-s) |f(X(s))|ds + \left|\int_0^t \sigma(s)dB(s) \right|, \quad t \geq 0.
	\end{align*}
	Let $\Omega_1$ be as before. By the Martingale Convergence Theorem
	(see Revuz and Yor \cite[Ch. V, Prop. 1.8]{revuzandyor}), 
	if $N = \{N_t,\,t\geq 0\}$ is a continuous local martingale with $\langle N, N \rangle_\infty <+\infty$, then
	\[
	\lim_{t\to\infty} N_t\in (-\infty,\infty), \quad\text{a.s.}.
	\] 
	In our case,
	\[
	\left\langle \int_0^{\cdot}\sigma(s)dB(s), \int_0^{\cdot}\sigma(s)dB(s) \right\rangle_t = \int_0^t \sigma^2(s)ds
	\]
	and thus $\sigma \in L^2(0,\infty)$ implies that $\lim_{t\to\infty} N_t$ exists and is finite a.s.
	Therefore, as $t\mapsto N_t$ is a.s. continuous, there exists an almost sure event $\Omega_2$ such that for all $\omega \in \Omega_2$
	\[
	\sup_{t\geq 0}\left|\int_0^t \sigma(s)dB(s)(\omega)\right| \leq N^\ast(\omega)<+\infty.
	\]
	Thus for all $\omega \in \Omega^\ast= \Omega_1 \cap \Omega_2$ and $t \geq 0$,
	\begin{align*}
		|X(t)| \leq |X(0)| + N^\ast + \int_0^t m(t-s) |f(X(s))|ds.
	\end{align*} 
	Using the estimate \eqref{eq.fglobalesteta} on $f$ and the finiteness of $\lim_{t\to\infty}m(t)$, we have
	\begin{align*}
		|X(t)| &\leq |X(0)| + N^\ast +  \bar{m} K(\eta) t + \bar{m}(1+\eta)\int_0^t \phi(|X(s)|)) ds, \quad t\geq 0.
	\end{align*} 
	Lastly, define $X(0)^\ast=|X(0)| + N^\ast$ and $H_\eta(t)=\bar{m} K(\eta) t$ so that
	\[
	|X(t)| \leq X(0)^\ast + H_\eta(t) + \bar{m}(1+\eta)\int_0^t \phi(|X(s)|)) ds, \quad t\geq 0.
	\]
	Note that this estimate is in precisely the form of \eqref{abs_est_1}. It is easy to show, as above, that $H_\eta(t)=\bar{m} K(\eta) t$
	obeys an estimate of the form \eqref{H_eta_est_2} for all $t\geq T_5(\eta)$. Hence for all $t\geq T(\eta)=T_5(\eta)$ and 
	for all $\omega\in \Omega^\ast$, the estimate 
	\begin{align}\label{abs_est_3}
		|X(t)| \leq X(0)^\ast + \eta \,\Phi^{-1}(\bar{m}t) + \bar{m}(1+\eta) \int_T^t \phi(|X(s)|)ds, \quad t\geq T(\eta),
	\end{align}
	holds. At this point we note that we are in the same position as in the proof of part (a) after \eqref{abs_est_2}, 
	and exactly analogous calculations yield 
	\[
	\limsup_{t\to\infty}\frac{F(|X(t)|)}{\bar{m}t} \leq 1 \,\, a.s.
	\]
\end{proof}

\begin{proof}[Proof of Corollary \ref{sigma_const}]
	We first prove that $\limsup_{t\to\infty}|X(t)| = \infty$ a.s. by showing that $X$ cannot be bounded with positive probability. Suppose there exists an event $A$, with positive probability, such that $|X(t)| \leq N < \infty$ for all $t\geq 0$ on $A$. Now consider the linear SDE
	\[
	dY(t) = -Y(t) dt + \sigma dB(t), \quad t >0, \quad Y(0) = 0. 
	\]
	The solution to the SDE above is given by $Y(t) = \sigma\int_0^t e^{-(t-s)}dB(s)$. Furthermore, it can be shown that $Y$ obeys $\limsup_{t\to\infty}|Y(t)| = \infty$ a.s. and $\liminf_{t\to\infty}|Y(t)| = 0$ a.s. (see Appleby et al. \cite[Theorem 4.1]{appleby2011characterisation}). Write \eqref{eq.introsfde} as
	\[
	dX(t) = - X(t) dt + \{X(s) + \int_0^t \mu(ds)f(X_{t-s})\}dt + \sigma dB(t), \quad t >0.
	\]
	Applying the variation of constants formula we obtain
	\begin{align*}
		X(t) &= e^{-t}X(0) + \int_0^t e^{-(t-s)} \left\{ X(s) + \int_0^s \mu(du)f(X_{s-u}) \right\}ds + \sigma\int_0^t e^{-(t-s)}dB(s)\\
		&= e^{-t}X(0) + \int_0^t e^{-(t-s)} \left\{ X(s) + \int_0^s \mu(du)f(X_{s-u}) \right\}ds + Y(t), \quad t \geq 0.
	\end{align*}
	With some simple estimation it follows that, on $A$, $\limsup_{t\to\infty}X(t) = \infty$, a contradiction. To show that $\limsup_{t\to\infty}F(|X(t)|)/\bar{m}t \leq 1$ a.s. we check $\sigma(t) = \sigma \in \mathbb{R}/\{0\}$ obeys $L_f(\Sigma)=0$, so we can apply Theorem \ref{Thm.stoch.zero}. By L'H\^{o}pital's rule
	\begin{align*}
		\lim_{t\to\infty}\frac{\Sigma(t)}{\int_0^t f(\Sigma(s))ds} = \lim_{t\to\infty}\frac{\Sigma '(t)}{f(\Sigma(t))},
	\end{align*}
	assuming the limit on the right--hand side exists. In fact 
	\[
	\Sigma '(t) = \frac{\sigma^2}{\log(t \sigma^2)\sqrt{2t \sigma^2 \log\log(t \sigma^2)}} + \frac{\sigma^2 \log\log(t \sigma^2)}{\sqrt{2t \sigma^2 \log\log(t \sigma^2)}}.
	\]
	Hence $\lim_{t\to\infty}\Sigma '(t) = 0$ and $L_f(\Sigma)=0$.
\end{proof}
\begin{proof}[Proof of Theorem \ref{Thm.stoch.L.F}]
	Let $\epsilon \in (0,1)$ be arbitrary and follow the line of argument from the proof of Theorem \ref{Thm.stoch.L.H} to obtain
	\begin{align*}
		|X(t)| &\leq A_\epsilon + (1+2\epsilon)\Sigma(t) + \bar{m}(1+\epsilon) \int_{T}^t \phi(|X(s)|)ds , \quad t \geq T, \quad \omega \in \Omega,
	\end{align*}
	where $A_\epsilon = \bar{m}\,T\, \sup_{s \in [0,T_1]}|X(s)|$. We define the upper comparison solution $X_\epsilon$ as in \eqref{X_epsilon} by
	\[
	X_\epsilon(t) = 1+  A_\epsilon + (1+2\epsilon)\Sigma(t) + \bar{m} (1+\epsilon)\int_{T}^t \phi(X_\epsilon(s))ds, \quad t \geq T. 
	\]
	Now by \eqref{Sigma_dom} there exists $T_1(\epsilon) > T$ such that
	\begin{align}\label{X_epsilon_est}
		X_\epsilon(t) \leq (1+3\epsilon)\Sigma(t) + \bar{m} (1+\epsilon)\int_{T}^t \phi(X_\epsilon(s))ds, \quad t \geq T_1(\epsilon). 
	\end{align}
	Let $I_\epsilon(t) = \int_T^t \phi(X_\epsilon(s))ds$; monotonicity yields
	\[
	\lim_{t\to\infty}\frac{\Sigma(t)}{\bar{m}\,I_\epsilon(t)} \leq \lim_{t\to\infty}\frac{\Sigma(t)}{\bar{m}\,\int_T^t\phi(\Sigma(s))ds} = L_\phi(\Sigma) \in (0,\infty).
	\]
	Hence there exists $T_2(\epsilon) > T_1$ such that 
	\begin{align}\label{Sigma_est_3}
		\Sigma(t) \leq L_\phi(\Sigma)\bar{m}(1+\epsilon)I_\epsilon(t), \quad t \geq T_2.
	\end{align}
	For $t \geq T_2$, using \eqref{Sigma_est_3}, calculate as follows
	\begin{align*}
		I_\epsilon '(t) &= \phi(X_\epsilon(t)) \leq \phi\left( (1+3\epsilon)\Sigma(t) + \bar{m} (1+\epsilon)I_\epsilon(t) \right)\\ &\leq \phi\left( L_\phi(\Sigma)\bar{m}(1+3\epsilon)(1+\epsilon)I_\epsilon(t) + \bar{m}(1+\epsilon)I_\epsilon(t) \right) \leq \phi\left( (1+7\epsilon)(\bar{m}+L_\phi(\Sigma)\bar{m})I_\epsilon(t) \right).
	\end{align*}
	Integrating the previous inequality we obtain
	\[
	\int_{T_2}^t \frac{I_\epsilon '(s)ds}{\phi\left( (1+7\epsilon)(\bar{m}+L_\phi(\Sigma)\bar{m})I_\epsilon(s) \right)} \leq t - T_2, \quad t \geq T_2.
	\]
	Hence making the substitution $u = (1+7\epsilon)(\bar{m}+L_\phi(\Sigma)\bar{m})I_\epsilon(s)$ yields
	\[
	\Phi\left( (1+7\epsilon)(\bar{m}+L_\phi(\Sigma)\bar{m})I_\epsilon(t) \right) \leq (t-T_2)(1+7\epsilon)(\bar{m}+L_\phi(\Sigma)\bar{m}) + \Phi_\epsilon, \quad t \geq T_2,
	\]
	where $\Phi_\epsilon = \Phi\left( (1+7\epsilon)(\bar{m}+L_\phi(\Sigma)\bar{m})I_\epsilon(T_2) \right)$. Thus
	\[
	(1+7\epsilon)(\bar{m}+L_\phi(\Sigma)\bar{m})I_\epsilon(t) \leq \Phi^{-1}\left( (t-T_2)(1+7\epsilon)(\bar{m}+L_\phi(\Sigma)\bar{m}) + \Phi_\epsilon \right), \quad t \geq T_2.
	\]
	Returning to \eqref{X_epsilon_est} and using the estimate above we obtain, for $t \geq T_2$,
	\begin{align*}
		X_\epsilon(t) &\leq (1+3\epsilon)L_\phi(\Sigma)\bar{m}(1+\epsilon)I_\epsilon(t) + \bar{m}(1+\epsilon)I_\epsilon(t)
		\leq (1+7\epsilon)(\bar{m}+L_\phi(\Sigma)\bar{m})I_\epsilon(t)\\
		&\leq \Phi^{-1}\left( (t-T_2)(1+7\epsilon)(\bar{m}+L_\phi(\Sigma)\bar{m}) + \Phi_\epsilon \right).
	\end{align*}
	It immediately follows that 
	\[
	\limsup_{t\to\infty}\frac{\Phi(X_\epsilon(t))}{\bar{m}t} \leq (1+L_\phi(\Sigma))(1+7\epsilon).
	\]
	Let $\epsilon \to 0^+$ and note that by construction $|X(t)| \leq X_\epsilon(t)$ for all $t \geq T$. Therefore,
	\[
	\limsup_{t\to\infty}\frac{\Phi(|X(t)|)}{\bar{m}t} \leq 1+L_\phi(\Sigma) \,\, a.s.,
	\]
	as required.
\end{proof}
\begin{proof}[Proof of Theorem \ref{Thm.stoch.L.H}]
	By L'H\^{o}pital's rule, $\lim_{x\to\infty}\Phi(x)/x = \lim_{x\to\infty}1/\phi(x) = 0$ and hence\\ $\lim_{t\to\infty}\Phi(\Sigma(t))/\Sigma(t) = 0$. Therefore, using Proposition \ref{PhiafterH},
	\begin{align}\label{Sigma_dom}
		\lim_{t\to\infty}\frac{A + B t}{\Sigma(t)} = \lim_{t\to\infty}\frac{A + B t}{\Phi(\Sigma(t))}\,\frac{\Phi(\Sigma(t))}{\Sigma(t)}= 0,
	\end{align}
	for any nonnegative constants $A$ and $B$. Arguing as in the proof of Theorem \ref{Thm.stoch.zero}, with $T$ and $\Omega$ defined analogously, we have the initial estimate
	\begin{align*}
		|X(t)|		&\leq |X(0)|+ \bar{m} K(\epsilon)t + (1+\epsilon)\Sigma(t) + \bar{m}(1+\epsilon) \int_0^t \phi(|X(s)|)ds, \quad t \geq T(\epsilon,\omega),\quad\omega \in \Omega,
	\end{align*}
	where $\mathbb{P}[\Omega]=1$. By \eqref{Sigma_dom} there is $T_1(\epsilon,\omega)>T(\epsilon,\omega)$ such that for all $t \geq T_1(\epsilon,\omega)$ $|X(0)| + \bar{m} K(\epsilon) t < \epsilon\Sigma(t)$.
	Hence
	\begin{align}\label{stoch.brownian.upper}
		|X(t)| &\leq (1+2\epsilon)\Sigma(t) + \bar{m}(1+\epsilon) \int_0^t \phi(|X(s)|)ds\\
		&\leq A_\epsilon + (1+2\epsilon)\Sigma(t) + \bar{m}(1+\epsilon) \int_{T_1}^t \phi(|X(s)|)ds , \quad t \geq T_1, \quad \omega \in \Omega,
	\end{align}
	where $A_\epsilon = \bar{m}\,T_1\, \sup_{s \in [0,T_1]}\phi(|X(s)|)$. Now define the function $X_\epsilon(t)$ for $t \geq T_1$ by
	\begin{align}\label{X_epsilon}
		X_\epsilon(t) = 1+  A_\epsilon + (1+2\epsilon)\Sigma(t) + \bar{m} (1+\epsilon)\int_{T_1}^t \phi(X_\epsilon(s))ds. 
	\end{align}
	By construction $|X(t)| \leq X_\epsilon(t)$ for all $t \geq T_1(\epsilon)$. Let $I_\epsilon(t) = \int_{T_1}^t \phi(X_\epsilon(s))ds$, so 
	\[
	I_\epsilon '(t) = \phi(X_\epsilon(t)) = \phi(1+  A_\epsilon + (1+2\epsilon)\Sigma(t) + \bar{m} (1+\epsilon)I_\epsilon(t) ), \quad t \geq T_1(\epsilon).
	\]
	Since $\phi$ is increasing and there exists a $T_2(\epsilon)> T_1(\epsilon)$ such that $1 + A_\epsilon < \epsilon \Sigma(t)$ for all $t \geq T_2$ we have
	\[
	I_\epsilon '(t) \leq \phi((1+3\epsilon)\Sigma(t) + \bar{m} (1+\epsilon)I_\epsilon(t)), \quad t\geq T_2.
	\]
	By the Mean Value Theorem there exists $\theta_t \in [0,1]$ such that 
	\begin{align}\label{I'_est_1}
		I_\epsilon '(t) &= \phi((1+3\epsilon)\Sigma(t))+ \phi'((1+3\epsilon)\Sigma(t) + \theta_t \bar{m} (1+\epsilon)I_\epsilon(t))\,\bar{m} (1+\epsilon)I_\epsilon(t)\nonumber\\
		&\leq \phi((1+3\epsilon)\Sigma(t))+ \phi'((1+3\epsilon)\Sigma(t))\,\bar{m} (1+\epsilon)I_\epsilon(t)\nonumber\\
		&\leq \phi((1+3\epsilon)\Sigma(t))+ \bar{m}(1+\epsilon)^2 \frac{\phi((1+3\epsilon)\Sigma(t))}{(1+3\epsilon)\Sigma(t)}I_\epsilon(t), \quad t \geq T_2,
	\end{align}
	where the final inequality follows from Lemma \ref{phi_props}. Once more we exploit the Mean Value Theorem and the first part of Lemma \ref{phi_props} as follows
	\begin{align}\label{phi_mvt_est}
		\phi((1+3\epsilon)\Sigma(t)) &= \phi(\Sigma(t)) + \phi'(\Sigma(t)+ \rho_t 3\epsilon\Sigma(t))\,3\epsilon\Sigma(t), \quad \rho_t \in [0,1]\nonumber\\
		&\leq \phi(\Sigma(t)) + \phi'(\Sigma(t))\,3\epsilon\Sigma(t) = \phi(\Sigma(t))\left\{ 1 + 3\epsilon \frac{\phi'(\Sigma(t))\Sigma(t)}{\phi(\Sigma(t))} \right\}\nonumber\\ &\leq \phi(\Sigma(t))(1+4\epsilon), \quad t \geq T^* > T_2.
	\end{align}
	Hence \eqref{I'_est_1} becomes
	\[
	I_\epsilon '(t) \leq (1+4\epsilon)\phi(\Sigma(t))+ \bar{m}(1+\epsilon)^2\frac{(1+4\epsilon)}{(1+3\epsilon)} \frac{\phi(\Sigma(t))}{\Sigma(t)}I_\epsilon(t), \quad t \geq T^*.
	\]
	Let 
	\[
	a_\epsilon(t) = \bar{m}(1+\epsilon)^2\frac{(1+4\epsilon)}{(1+3\epsilon)} \frac{\phi(\Sigma(t))}{\Sigma(t)}\quad \mbox{ and }\quad H_\epsilon(t) = \Sigma(t).
	\] 
	Now apply the argument from the proof of Theorem \ref{Thm.det.L.H} beginning at \eqref{I'est_key_1}. Following this line of argument shows that
	\[
	\limsup_{t\to\infty}\frac{I_\epsilon(t)}{\int_{T_1}^t \phi(\Sigma(s))ds} \leq N(\epsilon) < \infty.
	\]
	Returning to \eqref{X_epsilon} this yields
	\begin{align*}
		X_\epsilon(t) &< 1 + A_\epsilon + (1+2\epsilon)\Sigma(t) + \bar{m} (1+\epsilon)^2 N(\epsilon) \int_{T_1}^t \phi(\Sigma(s))ds, \quad t \geq T^*.
	\end{align*}
	Therefore
	\[
	\frac{X_\epsilon(t)}{\Sigma(t)} < 1+2\epsilon + \frac{1 + A_\epsilon}{\Sigma(t)} + \frac{\bar{m} (1+\epsilon)^2 N(\epsilon) \int_{T_1}^t \phi(\Sigma(s))ds}{\Sigma(t)}, \quad t \geq T^*.
	\]
	Thus 
	\[
	\limsup_{t\to\infty}\frac{X_\epsilon(t)}{\Sigma(t)} \leq 1+2\epsilon + \frac{\bar{m} (1+\epsilon)^2 N(\epsilon)}{L_\phi(\Sigma)} < \infty.
	\]
	Hence we have that $\limsup_{t\to\infty}|X(t)|/\Sigma(t) < \infty$ a.s.
	
	Suppose that $\limsup_{t\to\infty}|X(t)|/\Sigma(t) = 0$ on an event $\Omega_p$ of positive probability, then there exists $\bar{T}(\epsilon)>0$ such that $|X(t)| < \epsilon\Sigma(t)$ for all $t \geq \bar{T}$, $\omega\in \Omega_p$. 
	Let $J(t) = \int_0^t \bar{m}(t-s) f(X(s))ds$ and estimate as before. For all $\omega\in\Omega_p$, we obtain
	\begin{align}\label{J_est_again}
		|J(t)| &\leq \bar{m}\, \int_0^t C(\epsilon)+ (1+\epsilon)\phi(|X(s)|)ds\nonumber\\
		&\leq \bar{m}\,C(\epsilon)\,t + \bar{m}\,\bar{T}\,(1+\epsilon)\sup_{s \in [0,\bar{T}]}\phi(|X(s)|)+ \bar{m}(1+\epsilon)\int_{\bar{T}}^t \phi(|X(s)|)ds, \quad t \geq \bar{T}.
	\end{align}
	Hence
	\begin{align*}
		\limsup_{t\to\infty}\frac{|J(t)|}{\Sigma(t)} \leq \bar{m}(1+\epsilon)\limsup_{t\to\infty}\frac{\int_{\bar{T}}^t \phi(\epsilon \Sigma(s))ds}{\Sigma(t)} \leq \frac{1+\epsilon}{L_\phi(\Sigma)},\quad \mbox{ for all }\omega\in\Omega_p \mbox{ and }\epsilon \in (0,1).
	\end{align*}
	Therefore, because $L_f(\Sigma)>1$, $\limsup_{t\to\infty}|J(t)|/\Sigma(t) = \lambda \in [0,1)$ on $\Omega_p$. It follows that there exists $T' > \bar{T}$ such that $J(t)/\Sigma(t) > - \lambda - \epsilon$ for all $t \geq T'$. Consider the stochastic integral equation
	\[
	X(t) = X(0) + \int_0^t \bar{m}(t-s) f(X(s))ds + \int_0^t \sigma(s)dB(s), \quad t \geq 0.
	\] 
	For all $t \geq T'$ and $\omega\in\Omega_p$,
	\begin{align*}
		\frac{X(t)}{\Sigma(t)} &= \frac{X(0)}{\Sigma(t)} + \frac{J(t)}{\Sigma(t)} + \frac{\int_0^t \sigma(s)dB(s)}{\Sigma(t)}\geq  \frac{X(0)}{\Sigma(t)} +\frac{\int_0^t \sigma(s)dB(s)}{\Sigma(t)} - \lambda - \epsilon.
	\end{align*}
	This implies that $\limsup_{t\to\infty}X(t)/\Sigma(t) \geq 1 - \lambda - \epsilon$ for all $\omega\in\Omega_p$ and for all $\epsilon \in (0,1)$. Hence $\limsup_{t\to\infty}X(t)/\Sigma(t)\geq 1 - 1/L_\phi(\Sigma)$ on $\Omega_p$ and similarly $\liminf_{t\to\infty}X(t)/\Sigma(t)\leq -1 + 1/L_\phi(\Sigma)$ on $\Omega_p$, a contradiction. Hence $\mathbb{P}[\Omega_p]=0$ and
	\[
	\limsup_{t\to\infty}\frac{|X(t)|}{\Sigma(t)} = \Lambda \in (0,\infty)\, a.s.
	\]
	From \eqref{J_est_again} we obtain the following a.s. estimate
	\[
	|J(t)| \leq \bar{m}\,C(\epsilon)\,t + \bar{m}\,\bar{T}\,(1+\epsilon)\sup_{s \in [0,\bar{T}]}\phi(|X(s)|)+ \bar{m}(1+\epsilon)\int_{\bar{T}}^t \phi((\Lambda+\epsilon)\Sigma(s))ds, \quad t \geq \bar{T}.
	\]
	If we have $\Lambda \in (0,1)$, then we can choose $\epsilon > 0$ sufficiently small that $\Lambda + \epsilon < 1$ and monotonicity of $\phi$ and $\Sigma$ will yield $\limsup_{t\to\infty} |J(t)|/ \Sigma(t) \leq \Lambda/L_\phi(\Sigma)$, as before. If $\Lambda \in [1,\infty)$, we can estimate via the second part of Lemma \ref{phi_props}. Suppose $\Lambda \in [1,\infty)$, then
	\[
	\limsup_{t\to\infty}\frac{|J(t)|}{\Sigma(t)} \leq \bar{m}(1+\epsilon)(\Lambda+\epsilon)\frac{\int_{\bar{T}}^t\phi(\Sigma(s))ds}{\Sigma(t)} = (1+\epsilon)\frac{\Lambda+\epsilon}{L_\phi(\Sigma)},
	\]
	and letting $\epsilon\to 0^+$ we obtain $\limsup_{t\to\infty} |J(t)|/ \Sigma(t) \leq \Lambda/L_\phi(\Sigma)$ a.s. Therefore
	\begin{align*}
		\limsup_{t\to\infty}\frac{X(t)}{\Sigma(t)} \leq \Lambda \leq \limsup_{t\to\infty}\frac{|X(0)|}{\Sigma(t)} + \limsup_{t\to\infty}\frac{|J(t)|}{\Sigma(t)} + \limsup_{t\to\infty}\frac{|\int_0^t \sigma(s)dB(s)|}{\Sigma(t)} \leq \frac{\Lambda}{L_\phi(\Sigma)} + 1\,\, a.s.
	\end{align*}
	Finally $\Lambda \leq L_f(\Sigma)/(L_f(\Sigma)-1)$. Thus, $\limsup_{t\to\infty}X(t)/\Sigma(t) \leq L_f(\Sigma)/(L_f(\Sigma)-1)$ a.s. and similarly $\liminf_{t\to\infty}X(t)/\Sigma(t) \geq -L_f(\Sigma)/(L_f(\Sigma)-1)$ a.s.
\end{proof}
\begin{proof}[Proof of Theorem \ref{Thm.stoch.infty}]
	We follow closely the line of argument from the proof of Theorem \ref{Thm.stoch.L.H}. First we establish the required analogue of \eqref{Sigma_dom}. $L_f(\Sigma)=\infty$, so by Proposition \ref{PhiafterH} $\lim_{t\to\infty}\Phi(\Sigma(t))/\Sigma(t) = \infty$. Hence, for any nonnegative constants $A$ and $B$,
	\[
	\lim_{t\to\infty}\frac{A + Bt}{\Sigma(t)} = \lim_{t\to\infty}\frac{A + Bt}{\int_0^t f(\Sigma(s))ds}\,\frac{\int_0^t f(\Sigma(s))ds}{\Sigma(t)} = 0.
	\]
	With this result in hand we can proceed with the argument from Theorem \ref{Thm.stoch.L.H} to obtain
	\begin{align*}
		|X(t)| &\leq A_\epsilon + (1+2\epsilon)\Sigma(t) + \bar{m}(1+\epsilon) \int_{T_1}^t \phi(|X(s)|)ds , \quad t \geq T_1, \quad \omega \in \Omega,
	\end{align*}
	where $A_\epsilon = \bar{m}\,T_1\, \sup_{s \in [0,T_1]}|X(s)|$. Define $X_\epsilon(t)$ as in \eqref{X_epsilon} and with the same estimation as before 
	$
	\limsup_{t\to\infty}\int_T^t \phi(X_\epsilon(s))ds/\int_T^t \phi(\Sigma(s))ds < N(\epsilon) < \infty.
	$
	Therefore, since $L_f(\Sigma)=\infty$,
	\[
	\limsup_{t\to\infty}\frac{X_\epsilon(t)}{\Sigma(t)} \leq 1 + 2\epsilon + \bar{m}(1+\epsilon)^2N(\epsilon) \limsup_{t\to\infty}\frac{\int_T^t \phi(\Sigma(s))ds}{\Sigma(t)} = 1 + 2\epsilon.
	\]
	Note that $|X(t)| \leq X_\epsilon(t)$ a.s for all $t \geq T$ and let $\epsilon \to 0^+$ to conclude that 
	\[
	\limsup_{t\to\infty}\frac{|X(t)|}{\Sigma(t)} \leq 1\,\, a.s.
	\]
	The event on which $\limsup_{t\to\infty}|X(t)|/\Sigma(t) = 0$ is shown to have probability zero by exactly the line of argument which concludes the proof of Theorem \ref{Thm.stoch.L.H}. Hence $\limsup_{t\to\infty}|X(t)|/\Sigma(t) = \lambda \in (0,1]$ a.s. and $|X(t)| \leq (\lambda+\epsilon)\Sigma(t)$ for all $t \geq T(\epsilon)$ on an event of probability one. Once more using the notation that $J(t) = \int_T^t m(t-s)f(X(s))\,ds$ we recall the a.s. estimate \eqref{J_est_again}
	\begin{align*}
		|J(t)| &\leq \bar{m}\,C(\epsilon)\,t + \bar{m}\,\bar{T}\,(1+\epsilon)\sup_{s \in [0,T]}\phi(|X(s)|)+ \bar{m}(1+\epsilon)\int_{T}^t \phi(|X(s)|)ds, \quad t \geq T.
	\end{align*}
	Using the monotonicity of $\phi$, an estimate of the form \eqref{phi_mvt_est} and the hypothesis that $L_{\phi}(\Sigma)=\infty$,
	\begin{align*}
		\limsup_{t\to\infty}\frac{|J(t)|}{\Sigma(t)} &\leq \bar{m}(1+\epsilon)\limsup_{t\to\infty}\frac{\int_T^t \phi((\lambda+\epsilon)\Sigma(s))ds}{\Sigma(t)}\leq \bar{m}(1+\epsilon)\limsup_{t\to\infty}\frac{\int_T^t \phi((1+\epsilon)\Sigma(s))ds}{\Sigma(t)}\\ &\leq \bar{m}(1+\epsilon)(1+2\epsilon)\limsup_{t\to\infty}\frac{\int_T^t \phi(\Sigma(s))ds}{\Sigma(t)} = 0 \,\, a.s.
	\end{align*}
	Hence $\lim_{t\to\infty}J(t)/\Sigma(t)=0$ a.s. and the claim \eqref{X_over_sigma_to_zero} is proven. Now compute $\limsup_{t\to\infty}X(t)/\Sigma(t)$ as follows
	\[
	\limsup_{t\to\infty}\frac{X(t)}{\Sigma(t)} = \limsup_{t\to\infty}\left\{\frac{X(0)}{\Sigma(t)} + \frac{J(t)}{\Sigma(t)} + \frac{\int_0^t \sigma(s)dB(s)}{\Sigma(t)}\right\} = 1 \,\, a.s.
	\]
	Taking the liminf, rather than the limsup, in the equation above yields $\liminf_{t\to\infty}X(t)/\Sigma(t) = -1$ a.s., concluding the proof.
\end{proof}
\begin{proof}[Proof of Theorem \ref{Thm.Levy0}]
	First note that $\int_0^\infty {\gamma(s)}^{-\alpha}ds < \infty$ implies $\limsup_{t\to\infty}|Z(t)|/\gamma(t)=0$ a.s. (see Bertoin \cite[Theorem 5, Ch. VIII]{bertoin1998levy}). This proof follows by applying the argument used to establish Theorem \ref{Thm.stoch.zero} with $\Sigma$ replaced by $\gamma$ as appropriate.
\end{proof}
\begin{proof}[Proof of Theorem \ref{Thm.Levy1}]
	Suppose that $\gamma_+$ and $\gamma_-$ both satisfy the hypotheses on $\gamma$ with $\int_0^\infty {\gamma_+(s)}^{-\alpha}ds < \infty$ and $\int_0^\infty {\gamma_-(s)}^{-\alpha}ds = \infty$. It follows that
	\begin{equation}\label{levy_limits}
		\limsup_{t\to\infty}\frac{|Z(t)|}{\gamma_+(t)} = 0 \mbox{ a.s. and } \limsup_{t\to\infty}\frac{|Z(t)|}{\gamma_-(t)}=\infty \mbox{ a.s.}.
	\end{equation}
	Relevant properties of $\alpha$--stable processes can be found in Bertoin \cite[Theorem 5, Ch. VIII]{bertoin1998levy}. 
	
	We first deal with the claim that $\limsup_{t\to\infty}|X(t)|/\gamma_+(t) \leq 1/L_f(\gamma_+)$ a.s. when $L_f(\gamma_+) \in (1,\infty)$. Analogous to the beginning of the proof of Theorem \ref{Thm.stoch.L.H} use Proposition \ref{PhiafterH} to show that
	\[
	\lim_{t\to\infty}\frac{A+Bt}{\gamma_+(t)} = \lim_{t\to\infty}\frac{A+Bt}{\Phi(\gamma_+(t))}\frac{\Phi(\gamma_+(t))}{\gamma_+(t)} = 0,
	\]
	for any nonnegative constants $A$ and $B$. With the estimate above in hand and the proof proceeds as in that of Theorem \ref{Thm.stoch.L.H} but we arrive at a slightly different initial upper estimate to that derived in equation \eqref{stoch.brownian.upper} since we employ \eqref{levy_limits} for the asymptotics of $Z$. In this case
	\begin{align}
		|X(t)| &\leq A_\epsilon + 3\epsilon\, \gamma_+(t) + \bar{m}(1+\epsilon) \int_{T_1}^t \phi(|X(s)|)ds , \quad t \geq T_1, \quad \omega \in \Omega_1,
	\end{align}
	where $A_\epsilon= \bar{m}\,T_1\, \sup_{s \in [0,T_1]}\phi(|X(s)|)$. Now we are free to define the comparison solution
	\begin{align}\label{X_epsilon_levy}
		X_\epsilon(t) = 1+  A_\epsilon + 3\epsilon\,\gamma_+(t) + \bar{m} (1+\epsilon)\int_{T_1}^t \phi(X_\epsilon(s))ds, \quad t \geq T_1.
	\end{align}
	By following exactly the steps from the proof of Theorem \ref{Thm.stoch.L.H} we obtain $\limsup_{t\to\infty}|X_\epsilon(t)|/\gamma_+(t) < \infty$ with probability one and hence 
	\begin{equation}\label{levy_limsup_est}
		\limsup_{t\to\infty}\frac{|X(t)|}{\gamma_+(t)} < \infty \mbox{ a.s.}
	\end{equation}
	With the usual notation that $J(t) = \int_0^t \bar{m}(t-s)f(X(s))\,ds$ write
	\[
	\frac{|X(t)|}{\gamma_+(t)} \leq \frac{|X(0)|}{\gamma_+(t)} + \frac{|J(t)|}{\gamma_+(t)} + \frac{|Z(t)|}{\gamma_+(t)}.
	\]
	To finally derive the required bound on $\limsup_{t\to\infty}|X(t)|/\gamma_+(t)$ estimate $|J(t)|$ using \eqref{levy_limsup_est} (as was done in the proof of Theorem \ref{Thm.stoch.L.H}, for example); conclude by plugging in this estimate above and using \eqref{levy_limits}.
	
	The proof is essentially the same when $L_f(\gamma_+)=\infty$. To show that $\limsup_{t\to\infty}|X(t)|/\gamma_+(t) = 0$ a.s. proceed as before in applying the argument of Theorem \ref{Thm.stoch.L.H} but note now that this will give $\limsup_{t\to\infty}X_\epsilon(t)/\gamma_+(t) \leq 3\epsilon$ for the comparison solution. The conclusion now follows readily. 
	
	It remains to show that $\limsup_{t\to\infty}|X(t)|/\gamma_-(t) = \infty$ a.s. Begin by assuming to the contrary that there exists an event $\Omega_2$ with positive probability on which $\limsup_{t\to\infty}|X(t)|/\gamma_-(t) =: L \in [0,\infty)$. We first show that $\limsup_{t\to\infty}|J(t)|/\gamma_-(t) < \infty$ on an event of positive probability; work on $\Omega_2$ and estimate as follows
	\begin{align}\label{J_est_gamma_-}
		|J(t)| &\leq \bar{m} \int_0^t \left\{K + (1+\epsilon)\phi(|X(s))| \right\}\,ds  \nonumber\\
		&\leq \bar{m} K t + \bar{m}(1+\epsilon)T \sup_{s\in[0,T]}\phi(|X(s)|) + \bar{m}(1+\epsilon)\int_T^t \phi((1+\epsilon)L\gamma_-(s))\,ds \nonumber\\
		& \leq \bar{m} K t + \bar{m}(1+\epsilon)T \sup_{s\in[0,T]}\phi(|X(s)|) + \bar{m}(1+\epsilon)^2 \max(1,L)\int_T^t \phi(\gamma_-(s))\,ds,
	\end{align}  
	for $T$ sufficiently large and $t \geq T$ (the last inequality uses Lemma \ref{phi_props}). Divide by $\gamma_-$ and take the limsup across \eqref{J_est_gamma_-}; the final term on the right--hand side can be dealt with using the hypothesis $L_f(\gamma_-) \in (1,\infty]$, the first two terms are $o(\gamma_-)$ and hence we obtain
	\[
	\limsup_{t\to\infty}\frac{|J(t)|}{\gamma_-(t)} < \infty \mbox{ with positive probability}.
	\]
	Therefore the following holds on an event of positive probability 
	\[
	\limsup_{t\to\infty}\frac{|Z(t)|}{\gamma_-(t)} \leq \limsup_{t\to\infty} \left\{ \frac{|X(0)|}{\gamma_-(t)} + \frac{|X(t)|}{\gamma_-(t)} + \frac{|J(t)|}{\gamma_-(t)} \right\} < \infty,
	\]
	in contradiction of the fact that $\limsup_{t\to\infty}|Z(t)|/\gamma_-(t) = \infty$ a.s.
\end{proof}
\bibliography{perturbed_volterra_refs}

\begin{thebibliography}{10}

\bibitem{appleby2011characterisation}
J.~A. Appleby, J.~Cheng, and A.~Rodkina.
\newblock Characterisation of the asymptotic behaviour of scalar linear
  differential equations with respect to a fading stochastic perturbation.
\newblock {\em Discrete \& Contininuous Dynamical Systems Supplement}, pages
  79--90, 2011.

\bibitem{appleby2005mean}
J.~A. Appleby, S.~Devin, and D.~W. Reynolds.
\newblock Mean square convergence of solutions of linear stochastic {V}olterra
  equations to non--equilibrium limits.
\newblock {\em Dynam. Con. Disc. Imp. Sys. Ser A Math Anal. B}, 13:515--534,
  2005.

\bibitem{appleby2007almost}
J.~A. Appleby, S.~Devin, and D.~W. Reynolds.
\newblock Almost sure convergence of solutions of linear stochastic {V}olterra
  equations to nonequilibrium limits.
\newblock {\em Journal of Integral Equations and Applications}, pages 405--437,
  2007.

\bibitem{appleby2017growth}
J.~A. Appleby and D.~D. Patterson.
\newblock Growth rates of solutions of superlinear ordinary differential
  equations.
\newblock {\em Applied Mathematics Letters}, 71:30--37, 2017.

\bibitem{appleby2017memory}
J.~A. Appleby and D.~D. Patterson.
\newblock Memory dependent growth in sublinear {V}olterra differential
  equations.
\newblock {\em Journal of Integral Equations and Applications}, 29(4):531--584,
  2017.

\bibitem{appleby2003non}
J.~A. Appleby and D.~W. Reynolds.
\newblock Non-exponential stability of scalar stochastic {V}olterra equations.
\newblock {\em Statistics \& Probability letters}, 62(4):335--343, 2003.

\bibitem{appleby2004asymptotic}
J.~A. Appleby and A.~Rodkina.
\newblock Asymptotic stability of polynomial stochastic delay differential
  equations with damped perturbations.
\newblock {\em Functional Differential Equations}, 12(1--2):35--66, 2004.

\bibitem{sublinear2015}
J.~A.~D. Appleby and D.~D. Patterson.
\newblock Growth rates of sublinear functional and {V}olterra differential
  equations.
\newblock {\em SIAM Journal on Mathematical Analysis}, 50(2):2086--2110, 2018.

\bibitem{baker2000numerical}
C.~T. Baker and E.~Buckwar.
\newblock Numerical analysis of explicit one-step methods for stochastic delay
  differential equations.
\newblock {\em LMS Journal of Computation and Mathematics}, 3:315--335, 2000.

\bibitem{benhabib1991vintage}
J.~Benhabib and A.~Rustichini.
\newblock Vintage capital, investment, and growth.
\newblock {\em Journal of Economic theory}, 55(2):323--339, 1991.

\bibitem{berger1980volterra}
M.~A. Berger and V.~J. Mizel.
\newblock Volterra equations with {I}t{\^o} integrals — {I}.
\newblock {\em The Journal of Integral Equations}, pages 187--245, 1980.

\bibitem{bertoin1998levy}
J.~Bertoin.
\newblock {\em L{\'e}vy processes}, volume 121.
\newblock Cambridge University Press, 1998.

\bibitem{bharucha1972random}
A.~T. Bharucha-Reid.
\newblock {\em Random integral equations}.
\newblock Academic press, 1972.

\bibitem{bihari1956generalization}
I.~Bihari.
\newblock A generalization of a lemma of {B}ellman and its application to
  uniqueness problems of differential equations.
\newblock {\em Acta Mathematica Academiae Scientiarum Hungarica}, 7(1):81--94,
  1956.

\bibitem{BGT}
N.~H. Bingham, C.~M. Goldie, and J.~L. Teugels.
\newblock {\em Regular variation}, volume~27.
\newblock Cambridge University Press, 1989.

\bibitem{brockwell2009existence}
P.~J. Brockwell and A.~Lindner.
\newblock Existence and uniqueness of stationary {L}{\'e}vy-driven {CARMA}
  processes.
\newblock {\em Stochastic Process. Appl.}, 119(8):2660--2681, 2009.

\bibitem{corduneanu1991integral}
C.~Corduneanu.
\newblock {\em Integral equations and applications}.
\newblock Cambridge University Press, 1991.

\bibitem{dieudonne1971infinitesimal}
J.~A. Dieudonn{\'e}.
\newblock {\em Infinitesimal calculus}.
\newblock Hermann, 1971.

\bibitem{GLS}
G.~Gripenberg, S.-O. Londen, and O.~Staffans.
\newblock {\em Volterra integral and functional equations}, volume~34.
\newblock Cambridge University Press, 1990.

\bibitem{hartman2002ordinary}
P.~Hartman.
\newblock {\em Ordinary Differential Equations}.
\newblock SIAM, 2nd edition edition, 2002.

\bibitem{kolmanovskii2012applied}
V.~Kolmanovskii and A.~Myshkis.
\newblock {\em Applied Theory of Functional Differential Equations}, volume~85.
\newblock Springer Science \& Business Media, 2012.

\bibitem{kolmanovskii2013introduction}
V.~Kolmanovskii and A.~Myshkis.
\newblock {\em Introduction to the theory and applications of functional
  differential equations}, volume 463.
\newblock Springer Science \& Business Media, 2013.

\bibitem{lipovan2006}
O.~Lipovan.
\newblock Integral inequalities for retarded {V}olterra equations.
\newblock {\em Journal of Mathematical Analysis and Applications},
  322(1):349--358, 2006.

\bibitem{mao2007exponential}
X.~Mao.
\newblock Exponential stability of equidistant {E}uler--{M}aruyama
  approximations of stochastic differential delay equations.
\newblock {\em Journal of Computational and Applied Mathematics},
  200(1):297--316, 2007.

\bibitem{mao2007stochastic}
X.~Mao.
\newblock {\em Stochastic Differential Equations and Applications}.
\newblock Horwood Publishing, 2nd edition, 2007.

\bibitem{mao2006mean}
X.~Mao and M.~Riedle.
\newblock Mean square stability of stochastic {V}olterra integro-differential
  equations.
\newblock {\em Systems \& Control Letters}, 55(6):459--465, 2006.

\bibitem{marquardt2007multivariate}
T.~Marquardt and R.~Stelzer.
\newblock Multivariate {CARMA} processes.
\newblock {\em Stochastic Process and their Applications}, 117(1):96--120,
  2007.

\bibitem{metivier1980stopped}
M.~M{\'e}tivier and J.~Pellaumail.
\newblock On a stopped {D}oob's inequality and general stochastic equations.
\newblock {\em The Annals of Probability}, pages 96--114, 1980.

\bibitem{mohammed1996lyapunov}
S.-E.~A. Mohammed and M.~K. Scheutzow.
\newblock Lyapunov exponents of linear stochastic functional differential
  equations driven by semimartingales. {P}art {I}: the multiplicative ergodic
  theory.
\newblock In {\em Annales de l'IHP Probabilit{\'e}s et statistiques},
  volume~32, pages 69--105, 1996.

\bibitem{mohammed2003stable}
S.-E.~A. Mohammed and M.~K. Scheutzow.
\newblock The stable manifold theorem for non-linear stochastic systems with
  memory. {I}. existence of the semiflow.
\newblock {\em Journal of Functional Analysis}, 205(2):271--305, 2003.

\bibitem{mohammed2004stable}
S.-E.~A. Mohammed and M.~K. Scheutzow.
\newblock The stable manifold theorem for non-linear stochastic systems with
  memory: Ii. the local stable manifold theorem.
\newblock {\em Journal of Functional Analysis}, 206(2):253--306, 2004.

\bibitem{nualart2000large}
D.~Nualart and C.~Rovira.
\newblock Large deviations for stochastic {V}olterra equations.
\newblock {\em Bernoulli}, pages 339--355, 2000.

\bibitem{pituk1999hartman}
M.~Pituk.
\newblock The {H}artman--{W}intner theorem for functional differential
  equations.
\newblock {\em Journal of Differential Equations}, 155(1):1--16, 1999.

\bibitem{protter2004stochastic}
P.~E. Protter.
\newblock Stochastic integration and differential equation.
\newblock {\em Stochastic Modeling and Applied Probability}, 21, 2004.

\bibitem{revuzandyor}
D.~Revuz and M.~Yor.
\newblock {\em Continuous martingales and Brownian motion}, volume 293.
\newblock Springer Science \& Business Media, 1999.

\bibitem{reynolds2008decay}
D.~Reynolds and J.~Appleby.
\newblock Decay rates of solutions of linear stochastic {V}olterra equations.
\newblock {\em Electronic Journal of Probability}, 13:922--943, 2008.

\bibitem{shaikhet2013lyapunov}
L.~Shaikhet.
\newblock {\em Lyapunov functionals and stability of stochastic functional
  differential equations}.
\newblock Springer Science \& Business Media, 2013.

\bibitem{wu2012lasalle}
F.~Wu and S.~Hu.
\newblock The {L}a{S}alle-type theorem for neutral stochastic functional
  differential equations with infinite delay.
\newblock {\em Discrete \& Continuous Dynamical Systems-A}, 32(3):1065, 2012.

\bibitem{wu2010robustness}
F.~Wu, S.~Hu, and C.~Huang.
\newblock Robustness of general decay stability of nonlinear neutral stochastic
  functional differential equations with infinite delay.
\newblock {\em Systems \& Control Letters}, 59(3-4):195--202, 2010.

\bibitem{zhang2008euler}
X.~Zhang.
\newblock Euler schemes and large deviations for stochastic {V}olterra
  equations with singular kernels.
\newblock {\em Journal of Differential Equations}, 244(9):2226--2250, 2008.

\end{thebibliography}
\bibliographystyle{abbrv}
\newpage
\appendix
\section{Examples \& Numerical Experiments}\label{sec_examples}
\begin{proof}[Example \ref{example_L=1}]
	$x(t) = \exp\left( \lambda(t) + \sqrt{2(t+1)} \right)-e = \exp(P(t))-e$ for $t \geq 0$, with $\lambda(t) = (1+t)^\alpha$ for some $\alpha \in \left(0,1/2\right)$. We first show that  $\lim_{t\to\infty}H(t)/x(t) = 0$ which, combined with positivity of $H$, yields $\lim_{t\to\infty}x(t)/H(t) = \infty$. 
	\begin{align*}
		\lim_{t\to\infty}\frac{x(t)-x(0)-\int_0^t f(x(s))ds}{x(t)} &= 1 - \lim_{t\to\infty}\frac{\int_0^t f(x(s))ds}{x(t)} = 1 - \lim_{t\to\infty}\frac{f(x(t))}{x'(t)} \\
		&= 1 - \lim_{t\to\infty}\frac{\left( \alpha(1+t)^{\alpha-1}+[2(t+1)]^{-1/2} \right)^{-1}}{(1+t)^\alpha + \sqrt{2(t+1)}} = 0.
	\end{align*}
	Similarly,
	\[
	\lim_{t\to\infty}\frac{\int_0^t e^{-(t-s)}f(x(s))ds}{x(t)} = \lim_{t\to\infty}\frac{f(x(t))}{x(t)} = 0,
	\]
	and it then follows from \eqref{example_1} that $\lim_{t\to\infty}H(t)/x(t) = 0.$ Thus
	\begin{align}\label{asym_est_1}
		H(t) &= x(t) - x(0) - \int_0^t f(x(s))ds + \int_0^t e^{-(t-s)}f(x(s))ds\nonumber\\
		&\sim e^{P(t)} - \int_0^t \frac{e^{P(s)}}{P(s)}ds + f(x(t)) 
		\sim e^{P(t)} - \int_0^t \frac{e^{P(s)}}{P(s)}ds + \frac{e^{P(t)}}{P(t)}, \mbox{ as }t\to\infty.
	\end{align}
	We make the substitution $u = P(s)$ in the integral term and define $Q(s) = P(s)P'(s)$. Now estimate as follows
	\begin{align*}
		\int_0^t \frac{e^{P(s)}}{P(s)}ds &= \int_{P(0)}^{P(t)} \frac{e^u}{Q(P^{-1}(u))}du 
		= \int_{P(0)}^{P(t)} \frac{Q(P^{-1}(u))-1}{Q(P^{-1}(u))}e^u du + \int_{P(0)}^{P(t)}\, e^u\, du\\
		&= e^{P(t)} - e^{P(0)} + \int_{P(0)}^{P(t)} \frac{Q(P^{-1}(u))-1}{Q(P^{-1}(u))}\,e^u\, du.
	\end{align*}
	Combining this with \eqref{asym_est_1} we obtain
	\begin{align}\label{asym_est_2}
		H(t) \sim \frac{e^{P(t)}}{P(t)} + \int_{P(0)}^{P(t)} \frac{Q(P^{-1}(u))-1}{Q(P^{-1}(u))}e^u du, \mbox{ as }t\to\infty.
	\end{align}
	It remains to make an asymptotic estimate of the integral term on the right--hand side of equation \eqref{asym_est_2}. Expanding $Q$ yields
	\begin{align*}
		Q(s) &= \lambda(s)\lambda '(s) + \lambda(s) [2(s+1)]^{-1/2} + \lambda '(s) [2(s+1)]^{1/2} + 1\\
		&\sim 1 + \left( 2^{-1/2} + \alpha\sqrt{2} \right)s^{\alpha - 1/2} + o\left(s^{\alpha - 1/2}\right), \mbox{ as }s\to\infty.
	\end{align*}
	Hence $\lim_{s\to\infty}Q(s) =  1 + \lim_{s\to\infty}\left( 2^{-1/2} + \alpha\sqrt{2} \right)s^{\alpha - 1/2} = 1$ and since $P(s) \sim \sqrt{2s}$, $P^{-1}(s) \sim s^2 /2$, as $s\to\infty$. Therefore 
	\[
	\lim_{u\to\infty}\frac{Q(P^{-1}(u))-1}{Q(P^{-1}(u))} \sim \left( 2^{-1/2} + \alpha\sqrt{2} \right)\left(\frac{1}{2} \right)^{\alpha-1/2} u^{2\alpha - 1} = 0.
	\]
	It is straightforward to show that $\int_1^x u^{2\alpha - 1} e^u du \sim x^{2\alpha-1}e^x$ as $x\to\infty$ and thus
	\[
	\int_{P(0)}^{P(t)} \frac{Q(P^{-1}(u))-1}{Q(P^{-1}(u))}\,e^u\, du \sim K\,P(t)^{2\alpha-1}e^{P(t)}, \mbox{ as }t\to\infty,
	\]
	with $K$ a positive constant. Combining this with \eqref{asym_est_2} yields
	\[
	H(t) \sim \frac{e^{P(t)}}{P(t)} + K\,P(t)^{2\alpha-1}e^{P(t)} \sim K\,P(t)^{2\alpha-1}e^{P(t)}.
	\]
	Before calculating $\lim_{t\to\infty}H(t)/\int_0^t f(H(s))\,ds$ we note that
	\[
	f(H(t)) \sim \frac{K\,P(t)^{2\alpha-1}e^{P(t)}}{\log\left( K\,P(t)^{2\alpha-1}e^{P(t)}\right)} \sim K\,P(t)^{2\alpha-2}e^{P(t)}, \mbox{ as }t\to\infty.
	\]
	Hence, by L'H\^{o}pital's rule,
	\begin{align*}
		L_F(H) &= \lim_{t\to\infty}\frac{H(t)}{\int_0^t f(H(s))ds} = \lim_{t\to\infty}\frac{(2\alpha-1)P'(t)P(t)^{2\alpha-2}e^{P(t)} + P(t)^{\alpha-1}P'(t) e^{P(t)}}{P(t)^{2\alpha-2}e^{P(t)}}\\ &= \lim_{t\to\infty}(2\alpha-1)P'(t) + P'(t)P(t) = \lim_{t\to\infty}Q(t) = 1.
	\end{align*}
	Thus $L_f(H) = 1$ and $\lim_{t\to\infty}x(t)/H(t)=\infty$, as claimed.
\end{proof}
\begin{proof}[Example \ref{golden}]
	Write 
	\[
	H(t) = x(t) - x(0) - \int_0^t f(x(s))ds + \int_0^t e^{-(t-s)}f(x(s))ds, \quad t \geq 0,
	\]
	in order to work out the asymptotics of $H$. Firstly, 
	\[
	\int_0^t f(x(s))ds = A^\beta\,(1-\beta)^{\tfrac{\beta}{1-\beta}}\, \int_0^t s^{\tfrac{\beta}{1-\beta}}ds = A^\beta \, [(1-\beta)t]^{\tfrac{1}{1-\beta}}.
	\]
	Next
	\[
	\int_0^t e^{-(t-s)}f(x(s))ds \sim f(x(t)) = A^\beta \, [(1-\beta)t]^{\tfrac{\beta}{1-\beta}},\mbox{ as }t\to\infty,
	\]
	and hence this term will not affect the asymptotics of $H$. Thus
	\[
	H(t) \sim x(t) - f(x(t)) = (A - A^\beta)[(1-\beta)t]^{\tfrac{1}{1-\beta}}, \mbox{ as }t\to\infty.
	\]
	Now suppose instead that we had
	\[
	H(t) = [L_f(H)(1-\beta)t]^{\tfrac{1}{1-\beta}}, \quad t \geq 0.
	\]
	In this case
	\[
	\lim_{t\to\infty}\frac{H(t)}{\int_0^t f(H(s))ds} = \lim_{t\to\infty}\frac{H'(t)}{f(H(t))} = \lim_{t\to\infty}\frac{L_f(H)^{\tfrac{1}{1-\beta}}\,[(1-\beta)t]^{\tfrac{\beta}{1-\beta}} }{L_f(H)^{\tfrac{\beta}{1-\beta}}\,[(1-\beta)t]^{\tfrac{\beta}{1-\beta}}} = L_f(H).
	\]
	In order to choose $L_f(H)$ freely in this example we must solve $A - A^\beta = L_f(H)^{1/(1-\beta)}$ for $A \in [1,\infty)$, for a given $L_f(H) \in (0,\infty)$. To simplify the calculation choose $L_f(H) = 1$ and $\beta = 1/2$, so we must solve $A - A^{1/2} = 1$, or equivalently, $x^2 - x - 1 =0$, where $x^2 = A$. It follows that $x = (1 \pm \sqrt{5})/2$ and since we require $A \geq 1 $ select the solution $x = (1+\sqrt{5})/2$ (the so--called golden ratio). This yields $A \approx 2.618$. Finally, 
	\[
	\lim_{t\to\infty}\frac{F(x(t))}{\mu(\mathbb{R}^+)t} = A^{1-\beta} = A^{1/2} = \frac{1+\sqrt{5}}{2} \approx 1.618.
	\]  
\end{proof}
\begin{proof}[Example \ref{example_L>1}]
	Suppose $L_f(H) \in (1,\infty)$ and let $x(t) = \exp\left( \sqrt{2L_f(H)(t+1)} \right) - e$ for $t \geq 0$, we have
	\[
	f(x(t)) = [2(t+1)]^{-1/2}\,e^{\sqrt{2L_f(H)(t+1)}}.
	\]
	Integrating we obtain
	\begin{align*}
		\int_0^t f(x(s))ds &= \frac{1}{L_f(H)}\int_{\sqrt{2L_f(H)}}^{\sqrt{2L_f(H)(t+1)}}\, e^u\, du = \frac{1}{L_f(H)} e^{\sqrt{2L_f(H)(t+1)}} - \frac{1}{L_f(H)} e^{\sqrt{2L_f(H)}}.
	\end{align*}
	Therefore,
	\[
	x(t)-x(0)-\int_0^t f(x(s))\,ds \sim \left(\frac{L_f(H)-1}{L_f(H)}\right)e^{\sqrt{2L_f(H)(t+1)}}, \mbox{ as }t\to\infty.
	\]
	Using the fact that $f \circ x$ is sub-exponential and increasing we have
	\[
	\int_0^t e^{-(t-s)}f(x(s))ds \sim f(x(t)) = \frac{e^{\sqrt{2L_f(H)(t+1)}}}{\sqrt{2L_f(H)(t+1)}}, \mbox{ as }t\to\infty.
	\]
	Now, from \eqref{example_1}, we have $H(t) \sim \left((L_f(H)-1)/L_f(H)\right)e^{\sqrt{2L_f(H)(t+1)}}.$ It follows  that $\lim_{t\to\infty}x(t)/H(t) = L_f(H)/(L_f(H)-1)$. Finally,
	\begin{align*}
		\lim_{t\to\infty}\frac{H(t)}{\int_0^t f(H(s))ds} &= \lim_{t\to\infty}\frac{H'(t)}{f(H(t))}= \lim_{t\to\infty}\frac{\left(\tfrac{L_f(H)-1}{L_f(H)}\right)e^{\sqrt{2L_f(H)(t+1)}}L_f(H)}{\sqrt{2L_f(H)(t+1)}f(H(t))} 
		\\ &= \lim_{t\to\infty}\left(\frac{L_f(H)-1}{L_f(H)}\right)\,L_f(H)\,e^{\sqrt{2L_f(H)(t+1)}}
		\frac{\left(\tfrac{L_f(H)}{L_f(H)-1}\right)\sqrt{2L_f(H)(t+1)}}{\sqrt{2L_f(H)(t+1)}e^{\sqrt{2L_f(H)(t+1)}}}= L_f(H).
	\end{align*}
\end{proof}
\begin{proof}[Example \ref{example_L=infty}]
	With $x(t) = \exp\left([2(t+1)]^\alpha \right)-e$, $\alpha \in \left(\tfrac{1}{2},1\right)$, $t \geq 0$, we have
	\[
	f(x(t)) = [2(t+1)]^{-\alpha}\exp\left([2(t+1)]^\alpha \right).
	\]
	Hence 
	\[
	\lim_{t\to\infty}\frac{\int_0^t f(x(s))ds}{x(t)} = \lim_{t\to\infty}\frac{f(x(t))}{x'(t)} = \lim_{t\to\infty}\frac{1}{2\alpha[2(t+1)]^{2\alpha-1}} = 0.
	\]
	Similarly,
	\[
	\lim_{t\to\infty}\frac{\int_0^t e^{-(t-s)} f(x(s))ds}{x(t)} = \lim_{t\to\infty}\frac{f(x(t))}{x(t)} = 0,
	\]
	since $f \circ x$ is sub-exponential and $f$ is sublinear. It follows from \eqref{example_1} that $x \sim H$ and hence
	\[
	\lim_{t\to\infty}\frac{H(t)}{\int_0^t f(H(s))ds} = \lim_{t\to\infty}\frac{x(t)}{\int_0^t f(x(s))ds} = \infty,
	\]
	by the argument above for the limit of the reciprocal.
\end{proof}

\noindent \emph{Figure \ref{fig_sharpness}: Numerical Experiments}

To investigate further the results from Section \ref{stochastic} numerically we consider the following test equation
\begin{align}\label{eq.test}
	dX(t) = \left(\int_0^t e^{-(t-s)}f(X(s))ds\right)dt + \sigma(t)dB(t), \quad t >0,\quad X(0) = \psi \in \mathbb{R}.
\end{align}
Hence $\bar{m} = 1$ throughout this section. Let $I(t) = \int_0^t e^{-(t-s)}f(X(s))\,ds$ and rewrite \eqref{eq.test} as a coupled system as follows
\begin{align*}
	dI(t) &= \left( -I(t) + f(X(t))\right)dt, \quad t >0, \quad I(0) = 0,\\
	dX(t) &= I(t)dt + \sigma(t)dB(t),\quad t >0,\quad X(0) = \psi.
\end{align*}
Various authors have shown that an explicit one step Euler--Maruyama discretisation reliably approximates solutions to equations of the type \eqref{eq.test} in a mean--square sense for both fixed time lags and Volterra equations~\cite{baker2000numerical,mao2007exponential}. In particular, we can have confidence in such a numerical scheme if we assume global Lipschitz and linear growth conditions on both $f$ and $\sigma$. Using an Euler--Maruyama discretisation of the coupled form of \eqref{eq.test} we obtain, for $h>0$,
\begin{align*}
	I_{n+1} &= I_n + h \left( I_n - f(X_n) \right), \,\, n \geq 1, \,\, I_0 = 0,\nonumber\\
	X_{n+1} &= X_n +h \, I_n + \sigma(nh)\Delta W_n, \,\, n \geq 1, \,\, X_0 = \psi,
\end{align*}
where $\Delta W_n$ is a normal random variable with mean zero and variance $h$ for each $n \geq 1$. We take
\[
f(x) = \sign(x)|x|^\beta, \quad x\in \mathbb{R}, \quad \beta \in (0,1),
\]
so that $f$ obeys a global linear bound and $\lim_{|x|\to\infty}|f(x)|/f(|x|)=1$. We choose
\[
\sigma(t) = \frac{L_f(\Sigma)^{1/(1-\beta)} \, \left[ (1-\beta)t \right]^{(1+\beta)/(2-2\beta)}}{\sqrt{\log\log(t+e)}}, \quad t \geq 0,
\]
so that $\sigma \notin L^2(0,\infty)$ and $L_f(\Sigma) \in (0,\infty)$ is a free parameter. It is straightforward to show that 
\[
\Sigma(t) \sim \left(\, L_f(\Sigma) (1-\beta) t\, \right)^{1/1-\beta}, \quad \mbox{ as } t \to\infty.
\]

In Figure \ref{fig_sharpness} (left) we take $s(t) = (1+t)^{-3}$, $\beta=0.5$, $L_f(\Sigma) = 2$ and plot $s(t) X(t)/\Sigma(t)$. The scaled process $s(t)\,X(t)/\Sigma(t)$ fluctuates between $\pm s(t)L_f(\Sigma)/(L_f(\Sigma)-1)$ as it tends to zero a.s. As predicted by Theorem \ref{Thm.stoch.L.H}, the largest values of the scaled process are well approximated by $\pm s(t)L_f(\Sigma)/(L_f(\Sigma)-1)$. We also make a nonlinear transformation of the coordinates so that the convergence does not take place too quickly to observe. The co-ordinate transformation is given by $(x,y) \mapsto (\sign(x)|x|^\epsilon, \sign(y)|y|^\epsilon)$ with $\epsilon=0.08$; this nonlinear transformation is intended to be analogous to a log-log plot, typically used with non-negative data.

In Figure \ref{fig_sharpness} (right) we plot the quantity $F(|X(t)|)/\mu(\mathbb{R^+})t$ for various values of $\beta$ and $L_f(\Sigma)=1$ (without any scaling or coordinate transformations). We observe that $F(|X(t)|)/\mu(\mathbb{R^+})t$ appears to be pathwise bounded by $1+L_f(\Sigma)$ for large $t$, as guaranteed by Theorem \ref{Thm.stoch.L.F}.
\end{document}